\documentclass[aps,preprint,groupedaddress,showkeys]{revtex4-1}
\usepackage{graphicx}
\usepackage{dcolumn}
\usepackage{bm}
\usepackage{amssymb}
\usepackage{amsthm}
\usepackage{amsmath}
\usepackage{graphicx}
\usepackage{epstopdf}
\usepackage{mathrsfs}

\newtheorem{corollary}{\hskip\parindent\bf Corollary}
\newtheorem{theorem}{\hskip\parindent\bf Theorem}
\newtheorem{lemma}{\hskip\parindent\bf Lemma}
\newtheorem{remark}{\hskip\parindent\bf Remark}
\newtheorem{definition}{\hskip\parindent\bf Definition}
\begin{document}

\preprint{}

\title[Two delays induce  double Hopf bifurcation in a   predator-prey model]{Two delays induce Hopf bifurcation and  double Hopf bifurcation in a  diffusive Leslie-Gower  predator-prey system}

\author{Yanfei Du$^{1,2}$} \author{Ben Niu$^{2}$} \email{ niu@hit.edu.cn; author to whom correspondence should be addressed.}    \author{Junjie Wei$^{2}$}

\affiliation{$^1$College of Arts and Sciences,   Shaanxi University of Science and Technology, Xi'an 710021,  China.}
\affiliation{$^2$Department of Mathematics, Harbin  Institute of Technology at Weihai, Weihai 264209, China.}

\date{\today}

\begin{abstract}
In this paper, the dynamics of a modified Leslie-Gower predator-prey system with two delays and diffusion is considered. By calculating stability switching curves, the stability of positive equilibrium and the existence of  Hopf bifurcation and double Hopf bifurcation are investigated   on the  parametric plane of two delays. Taking two time delays as bifurcation parameters, the normal form on the center manifold near the double Hopf bifurcation point is derived, and  the unfoldings near the critical points are given.  Finally, we obtain the complex dynamics near the double Hopf bifurcation point, including the existence of  quasi-periodic solutions on a 2-torus,  quasi-periodic solutions on a 3-torus, and  strange attractors.

\end{abstract}

                           \maketitle
\textbf{Diffusive predator-prey models  with delays have been investigated widely, and the delay induced Hopf  bifurcation analysis has been well studied.  However, the study about bifurcation analysis of  predator-prey models with two simultaneously varying delays  has not been well established. Neither the  Hopf bifurcation theorem with two parameters nor the derivation process of normal form for two delays induced  double Hopf bifurcation has been   proposed in literatures. In this paper, we investigate a diffusive Leslie-Gower model with two delays, and carry out Hopf and double Hopf bifurcation analysis of the model.  Applying the method of studying characteristic equation with two delays, we get the stability switching curves and the crossing direction, after which we give the Hopf bifurcation theorem in two-parameter plane for the first time. Under some condition, the intersections of two stability switching curves are double Hopf bifurcation points. To figure out the dynamics near the double Hopf bifurcation point, we calculate the normal form on the center manifold. The derivation process of normal form we use in this paper  can be extended to other models   with two delays, one delay, or without delay.}
\section{Introduction}
\label{}

 The Leslie-Gower model, one of the most widely used  predator-prey models,  was proposed by  Leslie and Gower \cite{P. Leslie,PH Leslie}
   \begin{equation}
   \begin{aligned}
    &\dot{u}=r_1u(1-\frac{u}{K})-a uv,\\
      &\dot{v}=r_2v(1-\frac{v}{\gamma u}),
   \end{aligned}
     \end{equation}
     where $u(t)$ and $v(t)$ represent the populations of the prey and the predator at time $t$, respectively.
     $r_1$ and $r_2$ are the  intrinsic growth
     rates for prey and predator, respectively. $K$ is the environmental carrying capacity
     for prey population. $a$ is the per capita capturing rate of prey by a predator during unit  time. $\frac{v}{\gamma u}$ is Leslie-Gower term with carrying capacity of the predator $\gamma u$, which means that the carrying
     is proportional to the population size of the prey, and $\gamma$ is referred to as a measure of the quality of the prey as the food for the predator.  Since then, various researches on this model and modified models have been carried out. \cite{M. Aziz,J. Collings,P. Feng,Y. Ma,J. Zhou,Yuan S. L.1,Yuan S. L.2}

Refuges have important effects on the coexistence of predator and prey, and reducing the chance of extinction due to predation.  Chen et al. \cite{F. Chen} incorporated a  refuge protecting $mu$ of the prey into  Leslie-Gower system, which means that the remaining $(1-m)u$ of the prey is available to the predator. They considered the following Leslie-Gower predator-prey model
\begin{equation*}
\begin{aligned}
&\dot{u}=(r_1-b_1u)u-a_1(1-m)vu,\\
&\dot{v}=[ r_2-a_2\frac{v}{(1-m)u}] v,
\end{aligned}
\end{equation*}
 where $m\in \left[ 0,1\right) $   is a refuge protecting rate of the prey.

  Time delays are ubiquitous  in predator-prey systems.   It seems that time delays play an important role in  the stability of species densities. Researches are carried out to figure out the effect of delays on predator-prey systems.
 May \cite{R.M. May} considered the feedback time delay in prey growth, and the term $r_1u(t)(1-\frac{u(t-\tau)}{K})$ is the well-known delayed logistic equation.   Another type of time delay was introduced to  the negative feedback of the predator's density in  Leslie-Gower model in Refs. \cite{Yuan R.,A.F. Nindjin},  which denoted the time taken for digestion of the prey.       Liu et al. \cite{Leslie} considered both  delays mentioned above, and investigated  a modified Leslie-Gower predator-prey system with two delays described by the following system:
\begin{equation}
\label{odepredator}
     \left\lbrace            \begin{array}{lll}
\dot{u}(t) &=& r_1u(t)[ 1-\frac{u(t-\tau_1)}{K}] -a(1-m)u(t)v(t), \\
\dot{v}(t) &=&r_2v(t)[ 1-\frac{v(t-\tau_2)}{\gamma(1-m)u(t-\tau_2)}] , \\
\end{array}
\right.\end{equation}
with the initial conditions
\begin{equation}
\label{initial}
\begin{array}{l}
(\varphi_{1},\varphi_2)\in\textbf{C}([-\tau,0],\mathbb{R}_+^2),\varphi_i(0)>0,i=1,2,
\end{array}
\end{equation}
where  $\tau_1$ is  the feedback time delay in prey growth,  $\tau_2$  is the feedback time delay in predator growth, and we define $\tau={\rm max}\{\tau_1,\tau_2\}$.

 Related to systems with  two  delays,   general approach is to fix one delay, and vary another, or let  $\tau_1+\tau_2=\tau$. \cite{Song,K. Li,S. Ruan,C. Xu,L. Deng,Y. Ma}  However, sometimes, we want to investigate the dynamics of system when  two delays  vary simultaneously. To discuss systems with two delays, Gu et al. \cite{Gu K} analyzed the characteristic quasipolynomial
\begin{equation}
p(s)=p_0(s)+p_1(s)e^{-\tau_1s}+p_2(s)e^{-\tau_2s},
\end{equation}
where
\begin{equation*}
p_l(s)=\sum_{k=0}^np_{lk}s^k,
\end{equation*}
and provided a detailed study on the stability crossing curves such that the characteristic quasipolynomial has at least one imaginary zero and the crossing direction. Lin and Wang \cite{Lin X} considered the following  characteristic  functions
\begin{equation*}
D(\lambda;\tau_1,\tau_2)= P_{0}(\lambda)+P_{1}(\lambda)e^{-\lambda\tau_1}+P_{2}(\lambda)e^{-\lambda\tau_2}+P_{3}(\lambda)e^{-\lambda(\tau_1+\tau_2)}.
\end{equation*}
They derived an explicit expression for the stability switching curves in the $(\tau_1, \tau_2)$ plane, and gave a criterion to determine switching directions.

Since the preys and predators distribute inhomogenously in different locations, the diffusion should be taken into account  in more realistic ecological models. To  reveal new phenomena caused by the introduction of inhomogenous spatial environment, Du and Hsu \cite{Y. Du} considered a diffusive predator-prey model
\begin{equation}
\label{diffusion du}
\left\{
\begin{array}{ll}
 \dfrac{\partial u(x,t)} {\partial t}= d_1\Delta u(x,t)+\lambda u(x,t)-\alpha u(x,t)^2-\beta u(x,t)v(x,t),&x\in \Omega, t>0,\\
\dfrac{\partial v(x,t)}{\partial t }= d_2\Delta v(x,t)+\mu v(x,t)[1-\delta \frac{v(x,t)}{u(x,t)}],&x\in \Omega, ~t>0, \\
\dfrac{\partial u(x,t)} {\partial v}= 0,~~\dfrac{\partial v(x,t)} {\partial v}=0, & x\in \partial \Omega, t>0.\\
\end{array}
  \right.
\end{equation}
The
Neumann boundary condition means that no species can pass across the boundary
of $\Omega$. They showed the existence of steady-state solutions with certain prescribed spatial patterns.

 Motivated by the previous work, we consider the following modified Leslie-Gower predator-prey model with diffusion and Neumann boundary conditions
\begin{equation}
\label{diffusion predator}
\left\{
\begin{array}{l}
\begin{array}{l}
 \dfrac{\partial u(x,t)} {\partial t}= d_1\Delta u(x,t)+r_1u(x,t)[ 1-\frac{u(x,t-\tau_1)}{K}] -a(1-m)u(x,t)v(x,t),\\
\dfrac{\partial v(x,t)}{\partial t }= d_2\Delta v(x,t)+r_2v(x,t)[1-\frac{v(x,t-\tau_2)}{\gamma(1-m)u(x,t-\tau_2)}],~~~ \\
\end{array}
x\in [0,l\pi],~t>0\\
\dfrac{\partial u(x,t)} {\partial x}= 0,~~\dfrac{\partial v(x,t)} {\partial x}=0, ~at~ x=0~and~l\pi,\\
\end{array}
  \right.
\end{equation}
where $d_1,d_2>0$ are the diffusion coefficients characterizing the rate of the spatial dispersion of the prey and predator population, respectively.  $m\in [0,1)$ is a refuge
protecting rate of the prey. The spatial interval has been normalized as $[0,l\pi]$.

In fact, there are many literatures on    predator-prey model with delays, we refer to Refs. \cite{SHANSHAN CHEN,I. Al-Darabsah,Tian C,Yang R,Yuan R.} and the references therein.  Among them, effects of one  delay  on systems have been discussed widely. In this paper, we focus on the joint effect of two delays on  system (\ref{diffusion predator}). Adjusting the method  given in Ref. \cite{Lin X}, which was proposed to solve the analysis of stability switching of delayed differential equations, we can apply the method on reaction-diffusion systems with two delays, and obtain  the stability switching results when $(\tau_1,\tau_2)$ varies simultaneously.  To perform bifurcation analysis, we will extend the normal form method given by Faria and Magalh\~aes \cite{Faria,FariaJDE} to the double Hopf bifurcation analysis of the  reaction-diffusion system (\ref{diffusion predator}),
since the normal forms theory, which is  an efficient  method  for bifurcation analysis,    can be used to transform the original system to a qualitatively equivalent  equation with the simplest form by using  near-identity nonlinear transformations.

There are many  realistic problems with two delays,  such as epidemic model, \cite{K. L. Cooke,Jackson M.}   population interactions, \cite{H. I. Freedman,SHANSHAN CHEN,Song,K. Li,S. Ruan,C. Xu,L. Deng,Y. Ma}   neural networks, \cite{J. Wei and S. Ruan}  coupled
oscillators, \cite{Nguimdo} and so on. As we know, delay usually  destabilizes the equilibrium and induces Hopf bifurcation, which  gives rise to the periodic activities. For these systems, when we take two delays as parameters,  two Hopf bifurcation curves may intersect, and thus double Hopf bifurcation may occur, which is a source of  complicated dynamical behaviors.  The dynamics near double Hopf bifurcation is much more complicated than those near  Hopf bifurcation.  We can usually  observe rich dynamical behaviors, such as periodic and quasi-periodic oscillations,  coexisting of several oscillations, two- or three-dimensional invariant torus, and even chaos. \cite{Kuznetsov,Guckenheimer}      The analysis of double Hopf bifurcation provides a qualitative classification of bifurcating solutions arisen from  double Hopf bifurcation, which can help us to figure out the  dynamics  corresponding to different values of two delays near the double Hopf bifurcation singularity.
In fact,  the results obtained and the methods used in this paper can  be applied to the realistic problems with two delays mentioned above.   On one hand, as  was pointed out by Lin and Wang, \cite{Lin X}   the method of stability switching curves can be used to find the curves where the stability switches, and determine the crossing direction,   as long as  the characteristic equation of the  system we consider has the form of (\ref{character}).  On the other hand,  the normal form  derivation proposed in this paper can also be used to  reaction-diffusion systems with two delays and  Neumann boundary conditions,  by which the dynamics near the double Hopf singularity can be obtained.

This paper is organized as follows. In section \ref{existence}, we investigate the stability of the positive equilibrium  and the existence of Hopf bifurcation by the method of stability switching curves given in Ref. \cite{Lin X}.  In section  \ref{normal form}, taking two delays as bifurcation parameters, we derive the normal form on the center manifold near the double Hopf bifurcation point, and give the unfoldings near the critical points.  In section \ref{simulations}, we carry out some numerical simulations to support our analytical results.

\section{Stability switching curves and existence of double Hopf bifurcation}
\label{existence}
In this section, we  perform  bifurcation analysis  near the positive equilibrium  $E^*$ of system (\ref{diffusion predator}).  In order to discuss the joint effect of two delays $\tau_1$ and $\tau_2$ on system (\ref{diffusion predator}),  we will apply the method of the stability switching curves which is given in Ref.  \cite{Lin X}.

\subsection{Stability switching curves}
\label{Stability switching curves}
Clearly, system (\ref{diffusion predator}) has a unique positive constant equilibrium $E^*(u^*,v^*)$   with $u^*=\frac{Kr_1}{r_1+aK\gamma (1-m)^2}$ and   $v^*=\gamma(1-m)u^*$.

 The linearization of system (\ref{diffusion predator}) at the equilibrium $E^*$ is
     \begin{equation}
     \label{linear e2}
   \frac{\partial }{\partial t} \left( \begin{array}{l}
     u(x,t) \\
       v(x,t) \\
           \end{array}\right)
     =
    (D\Delta      +
     A) \left(  \begin{array}{l}
           u(x,t)\\
            v(x,t)\\
                     \end{array}
          \right)
    + B \left(  \begin{array}{l}
                u(x,t-\tau_1)\\
                 v(x,t-\tau_2)\\
                               \end{array}
               \right)
    +C\left(  \begin{array}{l}
            u(x,t-\tau_2)\\
              v(x,t-\tau_2)\\
              \end{array}
     \right),
          \end{equation}
where
  \begin{equation*}
   D=\left( \begin{array}{cc}
      d_1& 0\\
      0& d_2\\
         \end{array}\right),~
  A=\left( \begin{array}{cc}
    0& -a(1-m)u^*\\
    0& 0\\
       \end{array}\right),~
  B=\left( \begin{array}{cc}
     -\frac{r_1u^*}{K}& 0\\
     0 & 0\\
         \end{array}\right),~
    C=\left( \begin{array}{cc}
     0& 0\\
    \gamma(1-m)r_2 & -r_2\\
          \end{array}\right),
  \end{equation*}
  and $u(x,t),v(x,t)$ satisfy the homogeneous Neumann boundary condition.

 The characteristic equation of (\ref{linear e2}) is
 \begin{equation}\label{characterAG}
 {\rm det}(\lambda I_2-M_n-A-Be^{-\lambda\tau_1}-Ce^{-\lambda\tau_2})=0,
 \end{equation}
 where $I_2$ is the $2\times 2$ identity matrix and $M_n=-\frac{n^2}{l^2} D$, $n\in \mathbb{N}_0$.
The characteristic equation (\ref{characterAG}) is equivalent to \begin{equation}
\label{character}
D_n(\lambda;\tau_1,\tau_2)= P_{0,n}(\lambda)+P_{1,n}(\lambda)e^{-\lambda\tau_1}+P_{2,n}(\lambda)e^{-\lambda\tau_2}+P_{3,n}(\lambda)e^{-\lambda(\tau_1+\tau_2)}=0,
\end{equation}
where
\begin{equation*}
\begin{array}{l}
P_{0,n}(\lambda)=(\lambda+d_1\frac{n^2}{l^2})(\lambda+d_2\frac{n^2}{l^2}),\\P_{1,n}(\lambda)=\frac{r_1}{K}u^*(\lambda+d_2\frac{n^2}{l^2}),\\P_{2,n}(\lambda)=r_2(\lambda+d_1\frac{n^2}{l^2})+a(1-m)^2\gamma r_2u^*,\\P_{3,n}(\lambda)=\frac{r_1}{K}u^*r_2.
\end{array}
\end{equation*}

When $\tau_1=\tau_2=0$,  Eq. (\ref{character}) becomes
\begin{equation}
\label{charactertau120}
 \lambda^2+A\lambda+B=0,                                                                       \end{equation}
where
\begin{equation*}
\begin{aligned}
&A=d_1\frac{n^2}{l^2}+d_2\frac{n^2}{l^2}+\frac{r_1}{K}u^*+r_2>0,\\
&B=d_1d_2\frac{n^4}{l^4}+\frac{r_1}{K}u^*d_2\frac{n^2}{l^2}+r_2d_1\frac{n^2}{l^2}+a(1-m)^2\gamma r_2u^*+\frac{r_1}{K}u^*r_2>0.
\end{aligned}
\end{equation*}
It is clear that all roots of (\ref{charactertau120}) have negative real parts. Thus,
  when $\tau_1=\tau_2=0$, the positive equilibrium $E^*(u^*,v^*)$ is locally asymptotically stable.

\begin{remark}\label{Du}
    Du and Hsu  \cite{Y. Du} have proved that  the
positive equilibrium  of  the  diffusive Leslie-Gower predator-prey system (\ref{diffusion du}) is globally asymptotically stable under certain conditions.  We will also mention that Chen et al. \cite{SHANSHAN CHEN} have given a global stability result for  a  diffusive Leslie-Gower predator-prey system with two delays (the delay terms are different from these proposed in this paper). However, delays in this paper will destabilize the equilibrium. Thus, we only give the global stability result for $\tau_1=\tau_2=0$, which is a direct application of Proposition 2.1 in  Ref. \cite{Y. Du}.
Applying the globally asymptotical stability result in Ref. \cite{Y. Du} directly, we  have the following result:
 if $\frac{r_1}{K}>a(1-m)$,  the positive equilibrium $E^*(u^*,v^*)$ of (\ref{diffusion predator})  is globally asymptotically stable  when $\tau_1=\tau_2=0$.
\end{remark}

 In order to  apply the method of the stability switching curves, \cite{Lin X}
 we first verify the assumptions (i)-(iv) in Ref.  \cite{Lin X} are all true for any fixed $n$.
  \begin{itemize}

 \item[(i)] Finite number of characteristic roots on $\mathbb{C}_+=\{\lambda\in\mathbb{C}:{\rm Re}\lambda>0\}$ under the condition $${\rm deg} (P_{0,n}(\lambda))\geq {\rm max}\{{\rm deg}(P_{1,n}(\lambda)),{\rm deg}(P_{2,n}(\lambda)),{\rm deg}(P_{3,n}(\lambda))\}.$$

 \item[ (ii)] $P_{0,n}(0)+P_{1,n}(0)+P_{2,n}(0)+P_{3,n}(0)\neq 0$.
 \item[(iii)]
 $P_{0,n}(\lambda),P_{1,n}(\lambda),P_{2,n}(\lambda),P_{3,n}(\lambda)$ are coprime polynomials.
  \item[(iv)] $\lim\limits_{\lambda\rightarrow\infty} \left(\left| \frac{P_{1,n}(\lambda)}{P_{0,n}(\lambda)}\right|+  \left| \frac{P_{2,n}(\lambda)}{P_{0,n}(\lambda)}\right|+\left| \frac{P_{3,n}(\lambda)}{P_{0,n}(\lambda)}\right|\right) <1$.
 \end{itemize}

 In fact, condition (ii)-(iv) are obviously satisfied and (i) follows from  Ref. \cite{J. Hale}. Thus, similar to Ref. \cite{K.L. Cooke}, we have the following lemma.
    \begin{lemma}
 As the delays $(\tau_1,\tau_2)$ vary continuously in $\mathbb{R}_+^2$, the number of zeros (counting multiplicity) of $D_n(\lambda; \tau_1,\tau_2)$ on $\mathbb{C}_+$ can change only if a zero appears on or cross the imaginary axis.
 \end{lemma}

  To find the stability switching curves,  we  should seek  all the points $(\tau_1,\tau_2)$ such that $D_n(\lambda; \tau_1,\tau_2)$ has at least one zero on the imaginary axis.    Substituting $\lambda =i\omega~(\omega>0)$ into (\ref{character}), we obtain
  \begin{equation*}
  \label{characteriomega}
   (P_{0,n}(i\omega)+P_{1,n}(i\omega)e^{-i\omega\tau_1})+(P_{2,n}(i\omega)+P_{3,n}(i\omega)e^{-i\omega\tau_1})e^{-i\omega\tau_2}=0.
  \end{equation*}
 From $|e^{-i\omega\tau_2}|=1$, we have
  \begin{equation*}
  |P_{0,n}(i\omega)+P_{1,n}(i\omega)e^{-i\omega\tau_1}|=|P_{2,n}(i\omega)+P_{3,n}(i\omega)e^{-i\omega\tau_1}|.
  \end{equation*}
  Thus,  we have
  \begin{equation}\label{omegatau1}
  |P_{0,n}(i\omega)|^2+|P_{1,n}(i\omega)|^2-|P_{2,n}(i\omega)|^2-|P_{3,n}(i\omega)|^2=2A_{1,n}(\omega)\cos(\omega\tau_1)-2B_{1,n}(\omega)\sin(\omega\tau_1),
  \end{equation}
  with
  \begin{equation*}
  \begin{aligned}
  A_{1,n}(\omega)={\rm Re}(P_{2,n}(i\omega)\overline{P}_{3,n}(i\omega)-P_{0,n}(i\omega)\overline{P}_{1,n}(i\omega)),\\
   B_{1,n}(\omega)={\rm Im}(P_{2,n}(i\omega)\overline{P}_{3,n}(i\omega)-P_{0,n}(i\omega)\overline{P}_{1,n}(i\omega)).\\
  \end{aligned}
  \end{equation*}


If $A_{1,n}(\omega)^2+B_{1,n}(\omega)^2>0$,  there exists a function $\varphi_{1,n}(\omega)$ such that
 \begin{equation*}
    \begin{aligned}
        A_{1,n}(\omega)=\sqrt{A_{1,n}(\omega)^2+B_{1,n}(\omega)^2}\cos(\varphi_{1,n}(\omega)),\\
     B_{1,n}(\omega)=\sqrt{A_{1,n}(\omega)^2+B_{1,n}(\omega)^2}\sin(\varphi_{1,n}(\omega)),\\
    \end{aligned}
    \end{equation*}
where, $\varphi_{1,n}(\omega)={\rm arg}\{P_{2,n}(i\omega)\overline{P}_{3,n}(i\omega)-P_{0,n}(i\omega)\overline{P}_{1,n}(i\omega)\}\in(-\pi,\pi]$. Thus, (\ref{omegatau1}) can be written as
   \begin{equation}\label{omegatau1xin}
    |P_{0,n}(i\omega)|^2+|P_{1,n}(i\omega)|^2-|P_{2,n}(i\omega)|^2-|P_{3,n}(i\omega)|^2 =2\sqrt{A_{1,n}(\omega)^2+B_{1,n}(\omega)^2}\cos(\varphi_{1,n}(\omega)+\omega\tau_1).
    \end{equation}
 It is obvious that there exists  $\tau_1\in \mathbb{R}_+$ satisfying (\ref{omegatau1xin}) if and only if \begin{equation}\label{conditiontau1}
  \left(  (P_{0,n}(i\omega)|^2+|P_{1,n}(i\omega)|^2-|P_{2,n}(i\omega)|^2-|P_{3,n}(i\omega)|^2\right)^2 \leq 4(A_{1,n}(\omega)^2+B_{1,n}(\omega)^2).
 \end{equation}
 Denote  the set of $\omega\in\mathbb{R}_+$ which  satisfies (\ref{conditiontau1}) as $\Sigma^1_n$. We notice that (\ref{conditiontau1}) also includes the case $A_{1,n}^2(\omega)+B_{1,n}^2(\omega)=0$.

Denote
  \begin{equation*}
  \cos(\theta_{1,n}(\omega))=\dfrac{ |P_{0,n}(i\omega)|^2+|P_{1,n}(i\omega)|^2-|P_{2,n}(i\omega)|^2-|P_{3,n}(i\omega)|^2}{2\sqrt{A_{1,n}(\omega)^2+B_{1,n}(\omega)^2}},~~~~~\theta_{1,n}\in[0,\pi],
  \end{equation*}
  which leads to
  \begin{equation}\label{tau1}
  \tau_{1,j_1,n}^{\pm}(\omega)=\dfrac{\pm\theta_{1,n}(\omega)-\varphi_{1,n}(\omega)+2j_1\pi}{\omega},~~~j_1\in\mathbb{Z}.
  \end{equation}
  Similarly, we have
   \begin{equation}\label{tau2}
    \tau_{2,j_2,n}^{\pm}(\omega)=\dfrac{\pm\theta_{2,n}(\omega)-\varphi_{2,n}(\omega)+2j_2\pi}{\omega},~~~j_2\in\mathbb{Z},
    \end{equation}
   where
   \begin{equation*}
    \begin{aligned}
     \cos(\theta_{2,n}(\omega))=\dfrac{|P_{0,n}(i\omega)|^2-|P_{1,n}(i\omega)|^2+|P_{2,n}(i\omega)|^2-|P_{3,n}(i\omega)|^2}{2\sqrt{A_{2,n}(\omega)^2+B_{2,n}(\omega)^2}},~~~~~\theta_{2,n}\in[0,\pi],\\
    A_{2,n}(\omega)={\rm Re}(P_{1,n}(i\omega)\overline{P}_{3,n}(i\omega)-P_{0,n}(i\omega)\overline{P}_{2,n}(i\omega))=2\sqrt{A_{2,n}(\omega)^2+B_{2,n}(\omega)^2}\cos(\varphi_{2,n}(\omega)),\\
     B_{2,n}(\omega)={\rm Im}(P_{1,n}(i\omega)\overline{P}_{3,n}(i\omega))-P_{0,n}(i\omega)\overline{P}_{2,n}(i\omega))=2\sqrt{A_{2,n}(\omega)^2+B_{2,n}(\omega)^2}\sin(\varphi_{2,n}(\omega)).\\
    \end{aligned}
    \end{equation*}
Here the condition on $\omega$ is as follows
 \begin{equation}\label{conditiontau2}
  \left(  |P_{0,n}(i\omega)|^2-|P_{1,n}(i\omega)|^2+|P_{2,n}(i\omega)|^2-|P_{3,n}(i\omega)|^2\right)^2 \leq 4(A_{2,n}(\omega)^2+B_{2,n}(\omega)^2).
 \end{equation}
Denote the set of $\omega\in\mathbb{R}_+$ which  satisfies (\ref{conditiontau2}) as $\Sigma^2_n$.  In fact, we can easily show that (\ref{conditiontau1}) is equivalent  to (\ref{conditiontau2}) by squaring both sides of the two conditions (\ref{conditiontau1})  and   (\ref{conditiontau2}). Thus,  $\Sigma^1_n=\Sigma^2_n\stackrel{\vartriangle}{=}\Omega_n$.

\begin{definition}
The set
\begin{equation*}\begin{array}{r}
\Omega_n=\Bigg\{     \omega\in \mathbb{R}_+:F_n(\omega)\stackrel{\vartriangle}{=} ( |P_{0,n}(i\omega)|^2+|P_{1,n}(i\omega)|^2-|P_{2,n}(i\omega)|^2-|P_{3,n}(i\omega)|^2)^2\\ \left.-4(A_{1,n}(\omega)^2+B_{1,n}(\omega)^2)\leq 0  \Bigg\}\right.
\end{array}
\end{equation*}
is called the crossing set of $D_n(\lambda;\tau_1,\tau_2)=0$.
\end{definition}

Obviously, when $\omega\in \Omega_n$, both (\ref{conditiontau1}) and (\ref{conditiontau2}) hold. Now we consider the composition of set $\Omega_n$.

\begin{lemma}
The crossing set $\Omega_n$ consists of a finite number of intervals of finite length.
\end{lemma}
\noindent\textbf{Proof.} We follow the similar method in Ref. \cite{Lin X} to show this result.
 Since $F_n(\omega)$ is an eighth degree polynomial, and $F_n(+\infty)=+\infty$, $F_n(\omega)$ has a finite number of roots on $\mathbb{R}_+$.

If $F_n(0)>0$, Denote the  roots of $F_n(\omega)=0$ as $0<a_{1,n}<b_{1,n}\leq a_{2,n}<b_{2,n}<\cdots\leq a_{N,n}<b_{N,n}<+\infty$, and we $\Omega_n=\bigcup\limits_{j=1}^N\Omega_{j,n},~~ \Omega_{j,n}=[a_{j,n},b_{j,n}].$

 If $F_n(0)\leq 0$, denote the roots of $F_n(\omega)$ as $0<b_{1,n}\leq a_{2,n}<b_{2,n}<\cdots\leq a_{N,n}<b_{N,n}<+\infty$, and we have $\Omega_n=\bigcup\limits_{j=1}^N\Omega_{j,n},~~ \Omega_{1,n}=\left( 0,b_{1,n}\right] , \Omega_{j,n}=[a_{j,n},b_{j,n}]~~  (j\geq 2).       ~~~~~~~~~~\Box$

In fact, we can verify that when $\tau_1=\tau_{1,j_1,n}^+(\omega)$, we have $\tau_2=\tau_{2,j_2,n}^-(\omega)$, and when $\tau_1=\tau_{1,j_1,n}^-(\omega)$, we have $\tau_2=\tau_{2,j_2,n}^+(\omega)$. Denote
\begin{equation}\label{Tzfk}
\begin{aligned}
\mathcal{T}_{j_1,j_2,n}^{\pm j}&=\left\lbrace \left( \tau_{1,j_1,n}^{\pm}(\omega),\tau_{2,j_2,n}^{\mp}(\omega)\right):\omega\in\Omega_{j,n} \right\rbrace \\&=\left\lbrace \left( \dfrac{\pm\theta_{1,n}(\omega)-\varphi_{1,n}(\omega)+2j_1\pi}{\omega},\dfrac{\mp\theta_{2,n}(\omega)-\varphi_{2,n}(\omega)+2j_2\pi}{\omega}\right):\omega\in\Omega_{j,n} \right\rbrace,
\end{aligned}
\end{equation}
\begin{equation}\label{Tk}
\mathcal{T}^{j}_n=\bigcup_{j_1=-\infty}^{\infty}\bigcup_{j_2=-\infty}^{\infty}(\mathcal{T}_{j_1,j_2,n}^{+j}\cup\mathcal{T}_{j_1,j_2,n}^{-j})\cap \mathbb{R}_+^2,
\end{equation}
and
\[\mathcal{T}_n=\bigcup_{j=1}^N\mathcal{T}^j_n.\]

 \begin{definition}
  Any $(\tau_1,\tau_2)\in \mathcal{T}_n$ is called a crossing point, which makes $D_n(\lambda;\tau_1,\tau_2)=0$ have at least one root $i\omega$ with $\omega$ belongs to the crossing set $\Omega_n$. The set $\mathcal{T}_n$, which is the collection of all the crossing points, is called stability switching curves.
 \end{definition}

Since $F_n(a_{j,n})=F_n(b_{j,n})=0$, we have
\begin{equation*}
\theta_{i,n}(a_{j,n})=\delta_i^a\pi,~~\theta_{i,n}(b_{j,n})=\delta_i^b\pi,
\end{equation*}
where $\delta_i^a,\delta_i^b=0,1,i=1,2$. By (\ref{tau1}) and (\ref{tau2}), we can easily  confirm  that
\begin{equation}\label{connect}
\begin{aligned}
(\tau_{1,j_1,n}^{+j}(a_{j,n}),\tau_{2,j_2,n}^{-j}(a_{j,n}))=(\tau_{1,j_1+\delta_1^a,n}^{-j}(a_{j,n}),\tau_{2,j_2-\delta_2^a,n}^{+j}(a_{j,n})),\\
(\tau_{1,j_1,n}^{+j}(b_{j,n}),\tau_{2,j_2,n}^{-j}(b_{j,n}))=(\tau_{1,j_1+\delta_1^b,n}^{-j}(b_{j,n}),\tau_{2,j_2-\delta_2^b,n}^{+j}(b_{j,n})).\\
\end{aligned}
\end{equation}
Thus, for the stability switching curves corresponding to $\Omega_{j,n}$, $\mathcal{T}_{j_1,j_2,n}^{+j}$ is  connected to  $\mathcal{T}_{j_1+\delta_1^a,j_2-\delta_2^a,n}^{-j}$ at one end $a_{j,n}$,  and connected to  $\mathcal{T}_{j_1+\delta_1^b,j_2-\delta_2^b,n}^{-j}$ at the other end $b_{j,n}$.

\subsection{Crossing directions}
In the following, in order to identify the existence of Hopf bifurcation, we consider the  direction in which the root of (\ref{character}) cross the imaginary axis as $(\tau_1,\tau_2)$ deviates from a  stability switching curve $\mathcal{T}_n^j$ by the method given by Lin and Wang. \cite{Lin X}

Let $\lambda=\sigma+i\omega$. By (\ref{character}), and the implicit function theorem, $\tau_1$, $\tau_2$ can be expressed as function of $\sigma$ and $\omega$.
From (\ref{character}), we have
\begin{equation*}
\begin{aligned}
&\dfrac{\partial {\rm Re} D_n(\lambda;\tau_1,\tau_2)}{\partial \sigma}|_{\lambda=i\omega}
=R_0,
\dfrac{\partial {\rm Im} D_n(\lambda;\tau_1,\tau_2)}{\partial \sigma}|_{\lambda=i\omega} 
=I_0,
\end{aligned}
\end{equation*}
\begin{equation*}
\dfrac{\partial {\rm Re} D_n(\lambda;\tau_1,\tau_2)}{\partial \omega}|_{\lambda=i\omega}=-I_0,
\dfrac{\partial {\rm Im} D_n(\lambda;\tau_1,\tau_2)}{\partial \omega}|_{\lambda=i\omega}=R_0,
\end{equation*}
\begin{equation*}
\begin{aligned}
&\dfrac{\partial {\rm Re} D_n(\lambda;\tau_1,\tau_2)}{\partial \tau_l}|_{\lambda=i\omega} 
=R_l,
\dfrac{\partial {\rm Im} D_n(\lambda;\tau_1,\tau_2)}{\partial \tau_l}|_{\lambda=i\omega} 
=I_l,
\end{aligned}
\end{equation*}
where $l=1,2$.
By the implicit function theory, if $
{\rm det}\left(\begin{array}{cc}
R_1& R_2\\I_1&I_2
\end{array}\right)=R_1I_2-R_2I_1\neq 0.
$
 we have
\begin{equation}\label{juzhen}
\Delta(\omega):=\left( \begin{array}{cc}
\frac{\partial \tau_1}{\partial \sigma}&\frac{\partial \tau_1}{\partial \omega}\\\frac{\partial \tau_2}{\partial \sigma}&\frac{\partial \tau_2}{\partial \omega}
\end{array}\right) \arrowvert_{\sigma=0,\omega\in \Omega_n}=-\left( \begin{array}{cc}
R_1& R_2\\I_1&I_2
\end{array}\right)^{-1}\left( \begin{array}{cc}
R_0&-I_0\\I_0&R_0
\end{array}\right) .
\end{equation}

\begin{figure}[h]
                \centering
          \includegraphics[width=0.7\textwidth,height=0.4\textwidth]{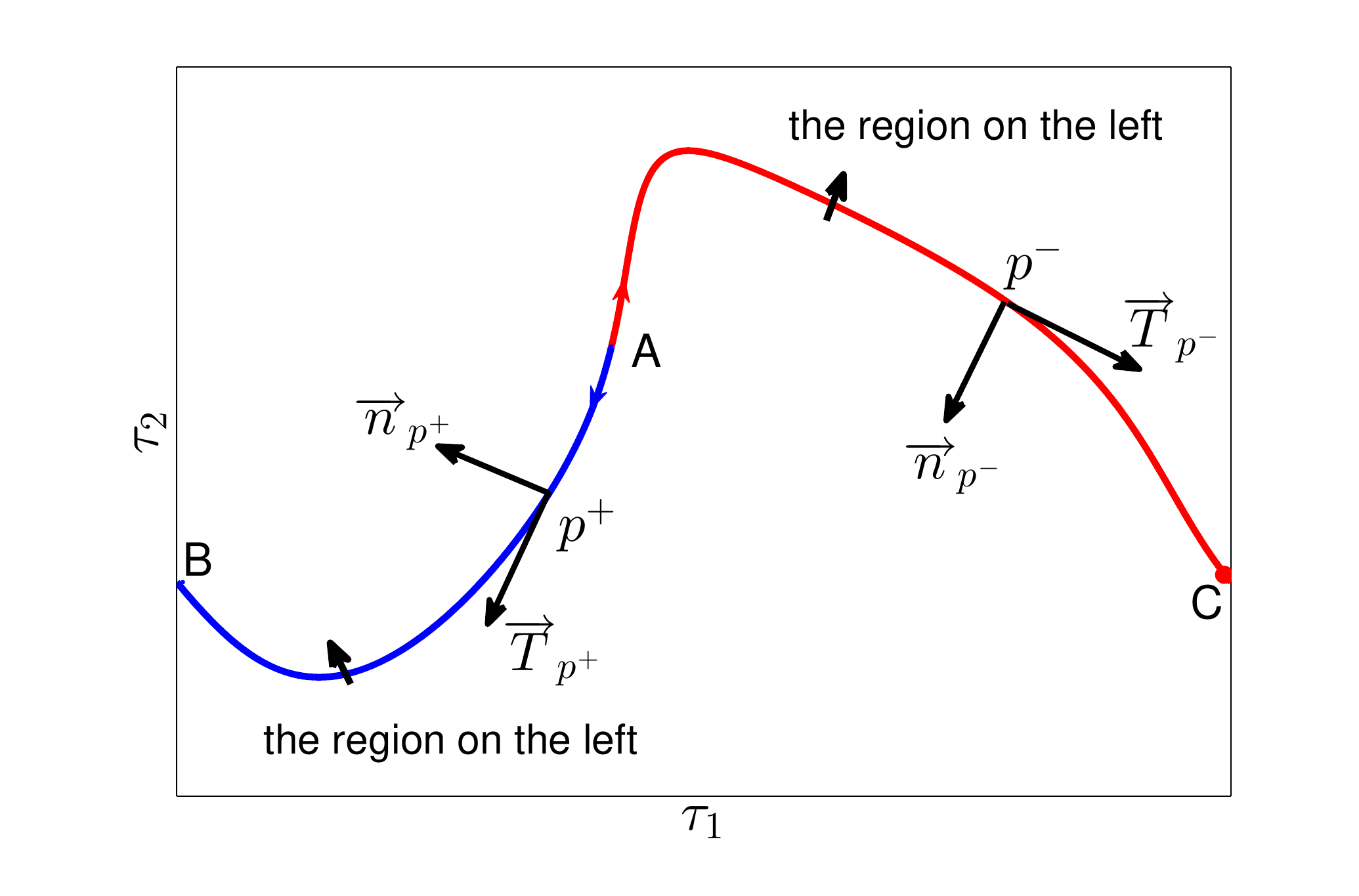}
      \caption  {This is a part of stability switching curves corresponding to $\Omega_{j,n}=[a_{j,n},b_{j,n}]$. The blue curve stands for $\mathcal{T}_{j_1,j_2,n}^{+ j}$, with two ends $A(\tau_{1,j_1,n}^{+j}(a_{j,n}),\tau_{2,j_2,n}^{-j}(a_{j,n}))$, and $B(\tau_{1,j_1,n}^{+j}(b_{j,n}),\tau_{2,j_2,n}^{-j}(b_{j,n}))$.   The red curve denotes $\mathcal{T}_{j_1+\delta_1^a,j_2-\delta_2^a,n}^{- j}$, which is   connected to $\mathcal{T}_{j_1,j_2,n}^{+ j}$ at $A$ corresponding to $a_{j,n}$, with the positive direction from $A$ to $C$.   }
       \label{fig:neighbor}
      \end{figure}

For any stability switching curves $\mathcal{T}_{j_1,j_2,n}^{\pm j}$, the direction of the curve corresponding to increasing $\omega\in \Omega_{j,n}$ is called the positive direction, i.e. from $(\tau_{1,j_1,n}^{\pm j}(a_{j,n}),\tau_{2,j_2,n}^{\mp j}(a_{j,n}))$ to $(\tau_{1,j_1,n}^{\pm j }(b_{j,n}),\tau_{2,j_2,n}^{\mp j}(b_{j,n}))$. The region on the left-hand (right-hand) side as we head in the positive directions of the curve is called the region on the left (right). As we have mentioned  in the previous section,  $\mathcal{T}_{j_1,j_2,n}^{+j}$ is  connected to  $\mathcal{T}_{j_1+\delta_1^a,j_2-\delta_2^a,n}^{-j}$ at  $a_{j,n}$,  then the positive direction of the two curves are opposite. Since the tangent vector of $\mathcal{T}_{j_1,j_2,n}^{\pm j}$ at $p^{\pm}(\tau_{1,j_1,n}^{\pm},\tau_{2,j_2,n}^{\mp})$ along the positive direction is $(\frac{\partial \tau_1}{\partial \omega},\frac{\partial \tau_2}{\partial \omega})\mid_{p^{\pm}}\stackrel{\vartriangle}{=}\overrightarrow{T}_{p^{\pm}}$, the normal vector of $\mathcal{T}_{j_1,j_2,n}^{\pm j}$  pointing to the right region is $(\frac{\partial \tau_2}{\partial \omega},-\frac{\partial \tau_1}{\partial \omega})\mid_{p^{\pm}}\stackrel{\vartriangle}{=}\overrightarrow{n}_{p^{\pm}}$ (see Fig. \ref{fig:neighbor}). On the other hand, as a pair of complex characteristic roots cross the imaginary axis to the right half plane, $(\tau_1,\tau_2)$ moves along the direction
$(\frac{\partial \tau_1}{\partial \sigma},\frac{\partial \tau_2}{\partial \sigma})\mid_{p^{\pm}}$.  We can conclude that if the inner product of these two vectors are positive, i.e.,
\begin{equation}\label{delta}
\delta(\omega)\mid_{p^{\pm}}:=\frac{\partial \tau_1}{\partial \sigma}\frac{\partial \tau_2}{\partial \omega}-\frac{\partial \tau_2}{\partial \sigma}\frac{\partial \tau_1}{\partial \omega}\mid_{p^{\pm}}>0,
\end{equation}
 Eq. (\ref{character}) has two more characteristic roots with positive real parts in the region on the right of  $\mathcal{T}_{j_1,j_2,n}^{\pm j}$.  If the inequality (\ref{delta}) is reversed, then the region on the left (\ref{character}) has two more characteristic roots with positive real parts.

It is easy to see that $\delta(\omega)={\rm det} \Delta(\omega)$.  Since $
{\rm det}\left(\begin{array}{cc}
-R_0& I_0\\-I_0&-R_0
\end{array}\right)=R_0^2+I_0^2\geq 0,$
(\ref{delta}) can be written as  $R_1I_2-R_2I_1>0$, if either $R_0\neq 0$ or $I_0\neq 0$, which is satisfied   since we do not consider the case that $i\omega$ is the multiple root of $D_n(\lambda;\tau_1,\tau_2)=0$,  i.e.,  $\frac{d D_n(\lambda;\tau_1,\tau_2)}{d \lambda}\mid_{\lambda=i\omega}=R_0+iI_0\neq 0$.

 We can verify that
\begin{equation}\begin{array}{l}
 R_1I_2-R_2I_1\mid_{p^{\pm}}\\={\rm Im}\{\overline{-i\omega (P_{1,n}e^{-i\omega\tau_{1,j_1,n}^{\pm}}+P_{3,n}e^{-i\omega(\tau_{1,j_1,n}^{\pm}+\tau_{2,j_2,n}^{\mp})})}(-i\omega)(P_{2,n}e^{-i\omega\tau_{2,j_2,n}^{\mp}}+P_{3,n}e^{-i\omega(\tau_{1,j_1,n}^{\pm}+\tau_{2,j_2,n}^{\mp})})\}\\
 =\pm\omega^2 |P_{2,n}\overline{P}_{3,n}-P_{0,n}\overline{P}_{1,n}|\sin \theta_{1,n}.\\
\end{array}
\end{equation}
Hence,
\begin{equation}\label{pone}
\begin{array}{l}
\delta(\omega\in \mathring{\Omega}_{j,n})\mid_{p^+}>0,~~\forall p^+\in\mathcal{T}_{j_1,j_2,n}^{+ j},~ and ~
\delta(\omega\in \mathring{\Omega}_{j,n})\mid_{p^-}<0,~~\forall p^-\in\mathcal{T}_{j_1,j_2,n}^{- j}\\
\end{array}
\end{equation}
since $\theta_{1,n}( \mathring{\Omega}_{j,n})\subset (0,\pi)$. Here, $\mathring{\Omega}_{j,n}$ denotes the interior of $\Omega_{j,n}$.

 We have the following conclusion.
\begin{lemma}\label{direction}
For any $j=1,2,\cdots,N$, we have
\begin{equation*}
\delta(\omega\in \mathring{\Omega}_{j,n}) >0(<0),~~\forall (\tau_1(\omega),\tau_2(\omega))\in\mathcal{T}_{j_1,j_2,n}^{+ j}((\tau_1(\omega),\tau_2(\omega))\in\mathcal{T}_{j_1,j_2,n}^{- j}).
\end{equation*}
Therefore, the region on the right  of $\mathcal{T}_{j_1,j_2,n}^{+j}$ $(\mathcal{T}_{j_1,j_2,n}^{- j})$ has two more (less) characteristic roots with positive real parts.
\end{lemma}

  We can see from  Fig. \ref{fig:neighbor} that the region on the right of $\mathcal{T}_{j_1,j_2,n}^{+ j}$, and the region on the left of $\mathcal{T}_{j_1+\delta_1^a,j_2-\delta_2^a,n}^{- j}$, with the black arrows pointing to,  have two more characteristic roots with positive real part.
Thus, as we move along these curves, stability crossing directions are consistent (see Fig.\ref{fig:neighbor}, the region with two more characteristic roots with positive real part is on the same side of $\mathcal{T}_{j_1,j_2,n}^{+j}$ and $\mathcal{T}_{j_1+\delta_1^a,j_2-\delta_2^a,n}^{-j}$).

Any given direction, $\overrightarrow{l}=(l_1,l_2)$, is pointing to the right region of the curve $\mathcal{T}_{j_1,j_2,n}^{\pm j}$,  if its inner product with the right-hand side  normal $(\frac{\partial \tau_2}{\partial \omega},-\frac{\partial \tau_1}{\partial \omega})$ is positive, i.e.,
\begin{equation}\label{l1l2}
l_1\frac{\partial \tau_2}{\partial \omega}-l_2\frac{\partial \tau_1}{\partial \omega}>0.
\end{equation}
And it is  pointing to the left region of the curve $\mathcal{T}_{j_1,j_2,n}^{\pm j}$  if its inner product with the right-hand side  normal  is negative.

We have the following  result
\begin{corollary}\label{coro}
As $(\tau_1,\tau_2)$ crosses the curve $\mathcal{T}_{j_1,j_2,n}^{\pm j}$ along the direction $\overrightarrow{l}=(l_1,l_2)$, there are two more (less) characteristic roots with positive real parts  if
\begin{equation}
-l_1(I_0I_1+R_0R_1)-l_2(I_0I_2+R_0R_2)>0 (<0).
\end{equation}
\end{corollary}
\noindent Proof. From (\ref{juzhen}),  the left side of (\ref{l1l2}) becomes\begin{equation}\label{cc}
[-l_1(I_0I_1+R_0R_1)-l_2(I_0I_2+R_0R_2)]/[R_1I_2-I_1R_2].
\end{equation}
If $-l_1(I_0I_1+R_0R_1)-l_2(I_0I_2+R_0R_2)>0$, $(l_1,l_2)$ is in the same (opposite) side as the right-hand side normal of $\mathcal{T}_{j_1,j_2,n}^{+j}$ ($\mathcal{T}_{j_1,j_2,n}^{-j}$). From Theorem \ref{direction}, we can conclude that there are two more characteristic roots with positive real parts  of (\ref{character}) as $(\tau_1,\tau_2)$ crosses the curve along the direction $\overrightarrow{l}=(l_1,l_2)$.  we can prove the result similarly when the inequality is reversed. ~~~~~~~~~~~~~ $\Box$
\subsection{Theorem of Hopf bifurcation}
   From the previous discussion, we have the following conclusion about Hopf bifurcation.
    \begin{figure}[h]
                    \centering
              \includegraphics[width=0.65\textwidth,height=0.4\textwidth]{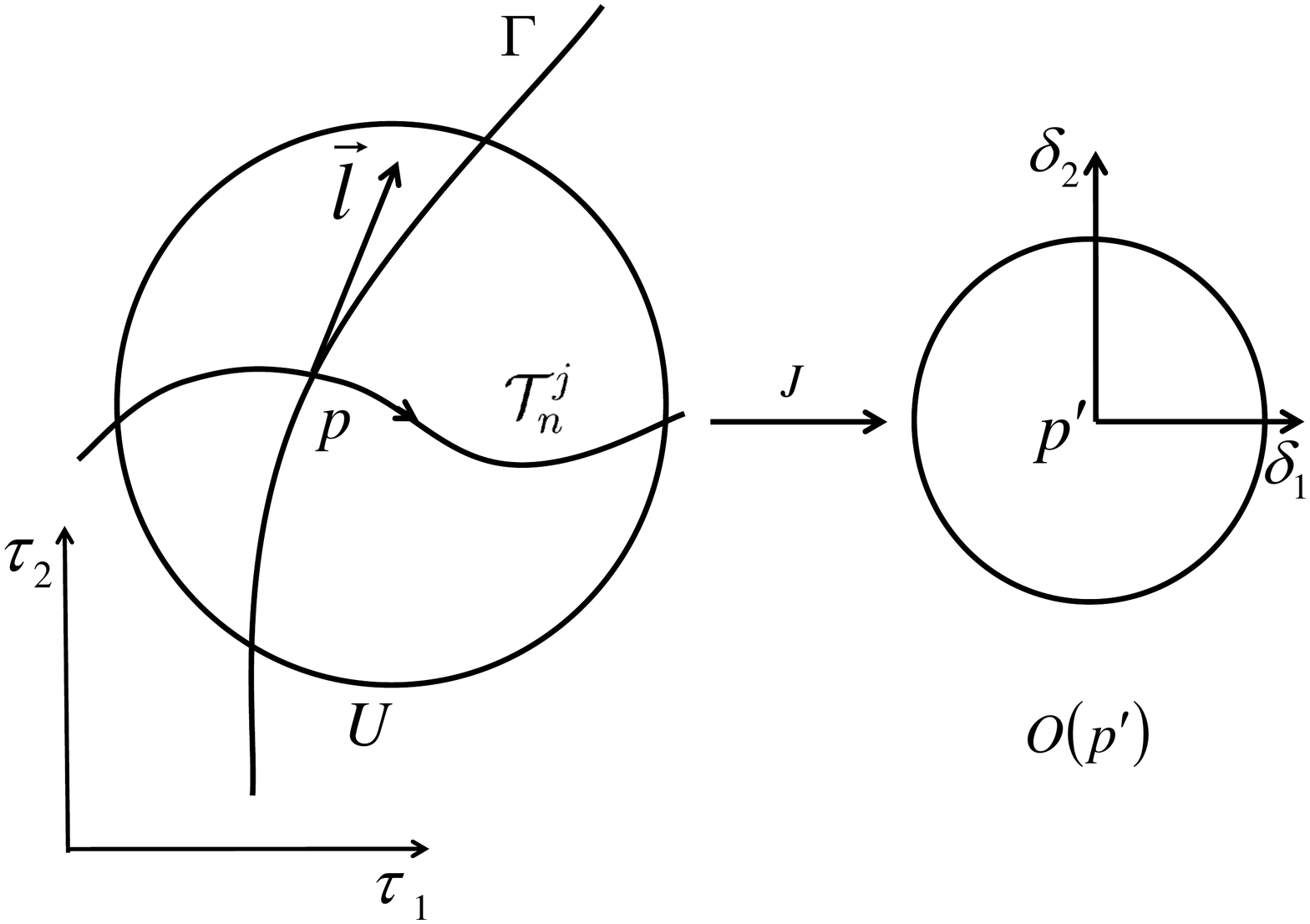}
          \caption  {A sketch of the transformation  form $(\tau_1,\tau_2)$ plane to $(\delta_1,\delta_2)$ plane.}
           \label{fig:transform}
          \end{figure}
\begin{theorem}\label{bifur}
For any $j=1,2,\cdots,N$, $\mathcal{T}_n^j$ is a Hopf bifurcation curve in the following sense:
for any $p\in \mathcal{T}_n^j$ and for any smooth curve $\Gamma$ intersecting with $\mathcal{T}_n^j$ transversely at $p$, we define the tangent of $\Gamma$ at $p$ by $\overrightarrow{l}$. If $
\frac{\partial {\rm Re}\lambda}{\partial \overrightarrow{l}}\mid_p\neq 0
$, and the other eigenvalues of (\ref{character}) at $p$ have non-zero real parts, then system (\ref{diffusion predator}) undergoes a Hopf bifurcation at $p$ when parameters $(\tau_1,\tau_2)$ cross $\mathcal{T}_n^j$ at $p$ along $\Gamma$.
\end{theorem}
\noindent Proof. Denote $p(\tau_{1}^0,\tau_2^0)$. Let $U$ be a neighbourhood of $p$. Suppose that the equation of curve $\Gamma$ is $\Gamma(\tau_1,\tau_2)=0$.   Introduce a mapping $J:U\rightarrow \mathbb{R}^2$, whose coordinate component function is expressed by $\left\lbrace \begin{array}{c}
\delta_1=\delta_1(\tau_1,\tau_2)\\
\delta_2=\delta_2(\tau_1,\tau_2)
\end{array}\right.$. Suppose that $J$ locally maps $p(\tau_1^0,\tau_2^0)$, $\mathcal{T}_n^j$ and $\Gamma$ to $p'(0,0)$, $\delta_1$ axis and $\delta_2$ axis, respectively (shown in Fig. \ref{fig:transform}), and the Jacobian determinant $\frac{\partial (\delta_1,\delta_2)}{\partial (\tau_1,\tau_2)}\arrowvert _p$ of mapping $J$ is not zero, then by  inverse function group theorem, there exists a  neighborhood of $p'$,  $O(p')$,  such that  there is a unique inverse mapping $J^{-1}$ of $J$, $$\left\lbrace \begin{array}{c}
\tau_1=\tau_1(\delta_1,\delta_2),\\
\tau_2=\tau_2(\delta_1,\delta_2),
\end{array}\right. ~~~~(\delta_1,\delta_2)\in O(p').$$

Now the characteristic equation of system (\ref{diffusion predator}) with $\delta_1=0$ has purely imaginary root $i\omega$ at $\delta_2=0$. We only need to further verify that $\frac{d {\rm Re} \lambda}{d \delta_2}\neq 0$.

In fact, we have $\delta_1=\Gamma(\tau_1,\tau_2)$, and the tangent vector  of curve $\Gamma$ is $\overrightarrow{l}=(-\frac{\partial \Gamma}{\partial \tau_2},\frac{\partial \Gamma}{\partial \tau_1})^T=(-\frac{\partial \delta_1}{\partial \tau_2},\frac{\partial \delta_1}{\partial \tau_1})^T$.   For convenience, denote  $J_1=\left(\begin{array}{cc}
  \frac{\partial \delta_1}{\partial \tau_1}&\frac{\partial \delta_1}{\partial \tau_2}\\
  \frac{\partial \delta_2}{\partial \tau_1}&\frac{\partial \delta_2}{\partial \tau_2}
  \end{array} \right) $.  Obviously, $\overrightarrow{e}_{\delta_2}=(0,1)^T=\frac{1}{{\rm det} J_1} J_1 ~\overrightarrow{l}$.

Since
$\frac{d {\rm Re} \lambda}{d \delta_2}=\frac{d \sigma}{d \delta_2}=\frac{\partial \sigma}{\partial  \overrightarrow{e}_{\delta_2}}=(\frac{\partial \sigma}{\partial \delta_1},\frac{\partial \sigma}{\partial \delta_2})^T  \cdotp \overrightarrow{e}_{\delta_2}$, we further treat  the inner product as  matrix multiplication  $(\frac{\partial \sigma}{\partial \delta_1},\frac{\partial \sigma}{\partial \delta_2}) \overrightarrow{e}_{\delta_2}=(\frac{\partial \sigma}{\partial \tau_1},\frac{\partial \sigma}{\partial \tau_2})J_1^{-1}  \frac{1}{{\rm det}J_1}J_1 \overrightarrow{l}=\frac{1}{{\rm det}J_1}(\frac{\partial \sigma}{\partial \tau_1},\frac{\partial \sigma}{\partial \tau_2})^T\cdotp   \overrightarrow{l}=\frac{1}{{\rm det}J_1}\frac{\partial {\rm \sigma}}{\partial \overrightarrow{l}} $.
 Thus, the transversality condition $\frac{d {\rm Re} \lambda}{d \delta_2}\neq 0\mid_{\delta_2=0}$ holds if $\frac{\partial {\rm Re \lambda}}{\partial \overrightarrow{l}}\mid_p\neq 0$.   According to  Corollary 2.4 in Ref.  \cite{JWu},   the conclusion follows. ~~~~~~~~~~~~~~$\Box$

\begin{remark}
Suppose  that there exists  $\omega_{j_1,k_1}\in \Omega_{j_1,k_1}$ and $\omega_{j_2,k_2}\in\Omega_{j_2,k_2}$, such that $\mathcal{T}_{k_1}^{j_1}$ and  $\mathcal{T}_{k_2}^{j_2}$ intersect. Then there are two pairs of pure imaginary roots of (\ref{character}) at the intersection.
\end{remark}

Thus, system (\ref{diffusion predator}) may undergoes double Hopf bifurcations near the positive equilibrium $E^*$ at the intersection of  two stability switching curves.

\section{Normal form on the center manifold  for  double Hopf bifurcation}
\label{normal form}
From the previous section, when two stability switching curves intersect, system (\ref{diffusion predator})  may undergoes double Hopf bifurcations near the positive equilibrium $E^*$.  In order to investigate the dynamical behavior of  system (\ref{diffusion predator}) near the  double Hopf bifurcation point,  we will calculate  the  normal forms of double Hopf bifurcation,   by applying  the normal form method of partial functional differential equations. \cite{Faria}

Without loss of generality, we always assume  $\tau_1>\tau_2$ in this section. Otherwise, all the derivations can be proceeded in  similar forms. Let $\overline{u}(x,t)=u(x,\tau_1t)-u^*,\overline{v}(x,t)=v(x,\tau_1t)-v^*$, and drop the bars, system (\ref{diffusion predator})  can be written as
    \begin{equation}
     \label{diffusion predator fourier}
   \frac{\partial }{\partial t} \left( \begin{array}{l}
     u(x,t) \\
       v(x,t) \\
           \end{array}\right)
     =
    \tau_1(D\Delta      +
     A) \left(  \begin{array}{l}
           u(x,t)\\
            v(x,t)\\
                     \end{array}
          \right)
    +\tau_1 B \left(  \begin{array}{l}
                u(x,t-1)\\
                 v(x,t-1)\\
                               \end{array}
               \right)
    +\tau_1C\left(  \begin{array}{l}
            u(x,t-\tau_2/\tau_1)\\
              v(x,t-\tau_2/\tau_1)\\
              \end{array}
     \right)+\tau_1\left(  \begin{array}{l}
                 f_1\\
                  f_2\\
                   \end{array}
          \right),
          \end{equation}
where
\begin{equation*}
\begin{array}{l}
f_1=-\frac{r_1}{K}u(x,t)u(x,t-1)-a(1-m)u(x,t)v(x,t),\\
f_2=-\frac{r_2}{\gamma(1-m)u^*}v(x,t)v(x,t-\tau_2/\tau_1)+\frac{r_2}{u^*}v(x,t)u(x,t-\tau_2/\tau_1)+\frac{r_2}{u^*}v(x,t-\tau_2/\tau_1)u(x,t-\tau_2/\tau_1)\\-\frac{r_2\gamma(1-m)}{u^*}u^2(x,t-\tau_2/\tau_1)+\frac{r_2}{\gamma(1-m)u^{*2}}v(x,t)u(x,t-\tau_2/\tau_1)v(x,t-\tau_2/\tau_1)\\-\frac{r_2}{u^{*2}}v(x,t)u^2(x,t-\tau_2/\tau_1)-\frac{r_2}{u^{*2}}v(x,t-\tau_2/\tau_1)u^2(x,t-\tau_2/\tau_1)+\frac{r_2\gamma(1-m)}{u^{*2}}u^3(x,t-\tau_2/\tau_1).
\end{array}
\end{equation*}

For the Neumann boundary condition, we
 define the real-valued Hilbert space \[X=\left\lbrace (u,v)^T\in H^2(0,l\pi)\times H^2(0,l\pi):\frac{\partial u}{\partial x}=\frac{\partial v}{\partial x}=0~ at ~ x=0,l\pi\right\rbrace,\]
  and  the  corresponding complexification space of $X$ by
    $X_{\mathbb{C}}:=X\oplus iX=\{U_1+iU_2:U_1,U_2\in X\}$,
   with  the general complex-value   $L^2$ inner product
   $\langle U,V\rangle=\int_0^{l\pi}(\overline{u}_1v_1+\overline{u}_2v_2)dx,$
    for $U=(u_1,u_2)^T$, $V=(v_1,v_2)^T\in X_{\mathbb{C}}$.
Let $\mathscr{C}:=C([-1,0],X_{\mathbb{C}})$ denote the phase space with the sup norm. We write $u^t\in\mathscr{C}$ for $u^t(\theta)=u(t+\theta)$, $-1\leq \theta\leq 0$.

Denote the double Hopf bifurcation point by $(\tau_1^*,\tau_2^*)$. Introduce two bifurcation parameters $\sigma=(\sigma_1,\sigma_2)$ by setting $\sigma_1=\tau_1-\tau_1^*$,  $\sigma_2=\tau_2-\tau_2^*$, and denote $U(t)=(u(t),v(t))^T$, then (\ref{diffusion predator fourier}) can be written as
   \begin{equation}
     \label{dudt}
   \dfrac{dU(t)}{dt}=D(\tau_1^*+\sigma_1,\tau_2^*+\sigma_2)\Delta U(t)+L(\tau_1^*+\sigma_1,\tau_2^*+\sigma_2)(U^t)+F(\tau_1^*+\sigma_1,\tau_2^*+\sigma_2,U^t),
          \end{equation}
where
 \begin{equation*}
      \begin{array}{l}
     D(\tau_1^*+\sigma_1,\tau_2^*+\sigma_2)=(\tau_1^*+\sigma_1)D=\tau_1^*D+\sigma_1 D,\\
    L(\tau_1^*+\sigma_1,\tau_2^*+\sigma_2)U^t\\=(\tau_1^*+\sigma_1)AU^t(0)+(\tau_1^*+\sigma_1)BU^t(-1)+(\tau_1^*+\sigma_1)CU^t(-(\tau_2^*+\sigma_2)/(\tau_1^*+\sigma_1))\\=\tau_1^*(AU^t(0)+BU^t(-1)+CU^t(-\tau_2^*/\tau_1^*))+\sigma_1(AU^t(0)+BU^t(-1)+CU^t(-\tau_2^*/\tau_1^*))\\+\tau_1^*C(U^t(-(\tau_2^*+\sigma_2)/(\tau_1^*+\sigma_1))-U^t(-\tau_2^*/\tau_1^*))+\sigma_1C(U^t(-(\tau_2^*+\sigma_2)/(\tau_1^*+\sigma_1))-U^t(-\tau_2^*/\tau_1^*)),\\
    F(\tau_1^*+\sigma_1,\tau_2^*+\sigma_2,U^t)=   (\tau_1^*+\sigma_1) (f_1,
    f_2)^T.
      \end{array}
         \end{equation*}
Consider the linearized system of (\ref{dudt})
   \begin{equation}
     \label{dudtlinear}
   \dfrac{dU(t)}{dt}=\tau_1^*D\Delta U(t)+\tau_1^*(AU^t(0)+BU^t(-1)+CU^t(-\tau_2^*/\tau_1^*))\stackrel{\vartriangle}{=}D_0\Delta U(t)+L_0(U^t).
          \end{equation}

 It is well known that the eigenvalues of $D\Delta$ on $X$ are $-d_1\frac{n^2}{l^2}$ and $-d_2\frac{n^2}{l^2}$, $n\in \mathbb{N}_0=\{0,1,2,\cdots\}$,
 with corresponding normalized eigenfunctions  $\beta_n^{1}(x)=\gamma_n(x)(1,0)^T$ and $\beta_n^{2}(x)=\gamma_n(x)(0,1)^T$, where
 $\gamma_n(x)=\dfrac{\cos\frac{n}{l}x}{\parallel\cos\frac{n}{l}x\parallel_{L^2}}$.  Define $\mathscr{B}_n$ the subspace of $\mathscr{C}$, by
$\mathscr{B}_n:={\rm span} \left\lbrace \langle v(\cdot),\beta_n^{j}\rangle\beta_n^{j}~\arrowvert  ~v\in \mathscr{C},j=1,2\right\rbrace$, satisfying $L(\mathscr{B}_n)(\tau_1,\tau_2)\subset {\rm span}\{\beta_n^{1},\beta_n^{2}\}$.
For simplification of notations,   we write
 $\left\langle v(\cdot),\beta_n \right\rangle=\left(
\langle v(\cdot),\beta_n^{1} \rangle,
\langle v(\cdot),\beta_n^{2} \rangle
\right)^T. $

System (\ref{dudt}) can be written as
          \begin{equation}
               \label{dudt2}
             \dfrac{dU(t)}{dt}=D_0\Delta U(t)+L_0(U^t)+G(\sigma,U^t).
                    \end{equation}
where $$\begin{aligned}
G(\sigma,U^t)&=\sigma_1(D\Delta U^t(0)+AU^t(0)+BU^t(-1)+CU^t(-\tau_2^*/\tau_1^*))\\&+(\tau_1^*+\sigma_1)C(U^t(-(\tau_2^*+\sigma_2)/(\tau_1^*+\sigma_1))-U^t(-\tau_2^*/\tau_1^*))\\&+F(\tau_1^*+\sigma_1,\tau_2^*+\sigma_2,U^t).
\end{aligned}$$
Rewrite (\ref{dudt2}) as an abstract ordinary differential equation on $\mathscr{C}$  \cite{Faria}
\begin{equation}\label{ode}
\frac{d}{dt}U^t=AU^t+X_0G(\sigma,U^t),
\end{equation}
where $A$ is the infinitesimal generator of the $C_0$-semigroup of solution maps of  the linear equation (\ref{dudt}), defined by
$
  \label{A}
  A:\mathscr{C}_0^1\cap\mathscr{C}\rightarrow \mathscr{C}, ~A\varphi=\dot{\varphi}+X_0[D_0\Delta\varphi(0)+L_0(\varphi)-\dot{\varphi}(0)],
  $
  with ${\rm dom}(A)=\{\varphi\in\mathscr{C}:\dot{\varphi}\in\mathscr{C},\varphi(0)\in {\rm dom}(\Delta)\}$,   and $X_0$ is given by
  $X_0(\theta)=0$ for  $\theta\in[-1,0)$ and $X_0(0)=I$.

   Then on $\mathscr{B}_n$, the linear equation $\frac{d}{dt}U(t)=D_0\Delta U(t)+L_0(U^t)$ is equivalent to the retarded functional differential equation on $\mathbb{C}^2$:
   $
   \label{RFDE}
   \dot{z}(t)=-\frac{n^2}{l^2}D_0z(t)+L_0z^t.
  $
   Define functions of bounded variation $\eta_m\in BV([-1,0],\mathbb{R})$ such that
   \begin{equation*}
   -\frac{k_m^2}{l^2}D_0\varphi(0)+L_0(\varphi)=\int_{-1}^0d\eta_m(\theta)\varphi(\theta), \varphi\in \mathscr{C}.
   \end{equation*}
   Let $A_m$ ($m=1,2$) denote the infinitesimal generator of the semigroup generated by (\ref{RFDE}),  and $A_m^*$ denote the formal adjoint of $A_m$ under the bilinear form
   \begin{equation*}
   (\alpha,\beta)_m=\alpha(0)\beta(0)-\int_{-1}^0\int_0^\theta\alpha(\xi-\theta)d\eta_m(\theta)\beta(\xi)d\xi.
   \end{equation*}

From the previous section, we know that system (\ref{dudtlinear}) has two pairs of pure imaginary eigenvalues $ \pm i\omega_{j_1,k_1}\tau_1^*,$ $\pm i\omega_{j_2,k_2}\tau_1^* $ at the double Hopf bifurcation point and the other eigenvalues with nonzero real  parts. For simplicity of notation, we denote them by $\{\pm i\omega_1\tau_1^*,\pm i\omega_2\tau_1^*\}$. Suppose that  $\omega_1:\omega_2\neq m:n$ for $m,n\in\mathbb{N}$ and $1\leq m,n\leq 3$, i.e., we do not consider the strongly resonant cases.
   Using the formal adjiont theory, we decompose $\mathscr{C}$ by $\{\pm i\omega_1\tau_1^*,\pm i\omega_2\tau_1^*\}$ as $\mathscr{C}=P_m\oplus Q_m$, where $Q_m=\{\phi\in \mathscr{C}:(\psi,\phi)_m=0,~for ~\psi~\in P_m^*\}$, $m=1,2$.
    We choose the basis
    $\varPhi_1(\theta)=(\phi_1(\theta),\overline{\phi}_1(\theta)), \varPhi_2(\theta)=(\phi_3(\theta),\overline{\phi}_3(\theta)$,
    $\varPsi_1(s)=(\psi_1(s),\overline{\psi}_1(s))^T$, $\varPsi_2(s)=(\psi_3(s),\overline{\psi}_3(s))^T$, in $P_1$, $P_1^*$, $P_2$, $P_2^*$, satisfying $(\Psi_m,\Phi_m)_{m}=I$, and
    \[A_m\Phi_m=\Phi_mB_m,~A_m^*\Psi_m=B_m \Psi_m,~ m=1,2,\]
    with $B_1={ \rm diag} (i\omega_1\tau_1^*,-i\omega_1\tau_1^*)$, $B_2={ \rm diag} (i\omega_2\tau_1^*,-i\omega_2\tau_1^*)$. Denote $\Phi(\theta)=(\Phi_1(\theta),\Phi_2(\theta))$, and $ \Psi(s)={ (\Psi_1(s),\Psi_2(s))^T}.$  By a few calculations, we have
     \[
    \begin{aligned}
    \phi_1(\theta)=(1,r_{12})^Te^{i\omega_1\tau_1^*\theta},\phi_3(\theta)=(1,r_{32})^Te^{i\omega_2\tau_1^*\theta},\\
    \psi_1(s)=D_1(1,r_{12}^*)e^{-i\omega_1\tau_1^*s}, \psi_3(s)=D_2(1,r_{32}^*)e^{-i\omega_2\tau_1^*s},
    \end{aligned}\]
     where
      \begin{equation}
            \begin{array}{ll}
    r_{12}=\dfrac{\gamma(1-m)r_2e^{-i\omega_1\tau_2^*}}{r_2e^{-i\omega_1\tau_2^*}+d_2\frac{n^2}{l^2}+i\omega_1},~~~~
    r_{32}=\dfrac{\gamma(1-m)r_2e^{-i\omega_2\tau_2^*}}{r_2e^{-i\omega_2\tau_2^*}+d_2\frac{n^2}{l^2}+i\omega_2},\\
    r_{12}^*=\dfrac{-a(1-m)u^*}{r_2e^{-i\omega_1\tau_2^*}+d_2\frac{n^2}{l^2}+i\omega_1},~~~~
   r_{32}^*=\dfrac{-a(1-m)u^*}{r_2e^{-i\omega_2\tau_2^*}+d_2\frac{n^2}{l^2}+i\omega_2},\\
    D_1=\frac{1 }{1+r_{12}^*r_{12}-\tau_1^*\frac{r_1}{K}u^*e^{-i\omega_1\tau_1^*}+r_{12}^*\gamma(1-m)r_2\tau_2^*e^{-i\omega_1\tau_2^*}-r_2\tau_2^*e^{-i\omega_1\tau_2^*}r_{12}^*r_{12} },\\D_2=\frac{1 }{1+r_{32}^*r_{32}-\tau_1^*\frac{r_1}{K}u^*e^{-i\omega_2\tau_1^*}+r_{32}^*\gamma(1-m)r_2\tau_2^*e^{-i\omega_2\tau_2^*}-r_2\tau_2^*e^{-i\omega_2\tau_2^*}r_{32}^*r_{32} }.
              \end{array}
         \end{equation}

 Now, we can decompose $\mathscr{C}$ into a center subspace and its orthocomplement, i.e.,
       \begin{equation}\label{ker}
       \mathscr{C}=\mathcal{P}\oplus{\rm Ker}\pi,
       \end{equation}
       where $\pi:\mathscr{C}\rightarrow\mathcal{P}$ is the projection defined by
       $
       \pi(\varphi)=\sum_{m=1}^2\Phi_m(\Psi_m,\langle\varphi(\cdot),\beta_{k_m}\rangle)_m\cdot\beta_{k_m},
       $
      with $\beta_{k_m}=\left( \beta_{k_m}^{1}, \beta_{k_m}^{2}\right) $, $m=1,2$.

Define the enlarged phase space \cite{Faria}
$\mathscr{BC}:=\{ \psi:[-1,0]\rightarrow X_{\mathbb{C}}:\psi {\rm ~is~ continuous~ on~} [-1,0), $ $\exists \lim_{\theta\rightarrow 0^-}\psi(\theta)\in X_{\mathbb{C} }\} .$

 According  to (\ref{ker}),   $U^t\in X$  can be composed as
 \begin{equation}\label{Ut}
 U^t(\theta)=\phi_1(\theta)z_1\gamma_{k_1}+\overline{\phi}_1(\theta)z_2\gamma_{k_1}+\phi_3(\theta)z_3\gamma_{k_2}+\overline{\phi}_3(\theta)z_4\gamma_{k_2}+w(\theta)\stackrel{\vartriangle}{=}\Phi(\theta) z_x+w(\theta),
 \end{equation}
where $w \in \mathscr{C}^1\bigcap{\rm Ker}\pi:=\mathcal{Q}^1$ for any $t$.
 Then in $\mathscr{BC}$ the system (\ref{ode}) is equivalent to the system
  \begin{equation}\label{zdoty}
  \begin{array}{l}
   \dot{z_1}=i\omega_1\tau_1^*z_1+\psi_1(0) \langle G(\sigma,\Phi(\theta) z_x+w(\theta)),\beta_{k_1}\rangle,\\
    \dot{z_2}=-i\omega_1\tau_1^*z_2+\overline{\psi}_1(0)\langle G(\sigma,\Phi(\theta) z_x+w(\theta)),\beta_{k_1}\rangle,\\
     \dot{z_3}=i\omega_2\tau_1^*z_3+\psi_3(0) \langle G(\sigma,\Phi(\theta) z_x+w(\theta)),\beta_{k_2}\rangle,\\
      \dot{z_4}=-i\omega_2\tau_1^*z_1+\overline{\psi}_3(0) \langle G(\sigma,\Phi(\theta) z_x+w(\theta)),\beta_{k_2}\rangle,\\
  \frac{dw}{dt} =A_1w+(I-\pi)X_0 G(\sigma,\Phi(\theta) z_x+w(\theta)),
  \end{array}
  \end{equation}
  where  $A_1$ is  the restriction of $A$ on  $\mathcal{Q}^1\subset{\rm Ker \pi\rightarrow{\rm Ker}\pi}$, $A_1\varphi=A\varphi$ for $\varphi\in \mathcal{Q}^1$.

 Consider the formal Taylor expansion
 \begin{equation*}
 G(\sigma,\varphi)=\frac{1}{2!}G_2(\sigma,\varphi)+\frac{1}{3!}G_3(\sigma,\varphi),
 \end{equation*}
 where $G_j$ is the $j{\rm th}$ Fr\'{e}chet derivation of $G$, which we calculate  in section  1 of supplement material. Then (\ref{zdoty}) can be written  as
 \begin{equation}\label{zdotTaylor}
 \begin{array}{l}
 \dot{z}=Bz+\sum\limits_{j\geq 2}\frac{1}{j!}f_j^1(z,w,\sigma),\\
 \frac{d}{dt}w =A_1w+\sum\limits_{j\geq 2}\frac{1}{j!}f_j^2(z,w,\sigma),
 \end{array}
 \end{equation}
 where  $z=(z_1,z_2,z_3,z_4)\in \mathbb{C}^4, w\in \mathcal{Q}^1$,  $B={\rm diag}$,  $ (B_1,B_2)={\rm diag (i\omega_1\tau_1^*,-i\omega_1\tau_1^*,i\omega_2\tau_1^*,-i\omega_2\tau_1^*)}$, and $f_j=(f_j^1,f_j^2), j\geq 2$, are defined by
 \begin{equation}
 \label{fj1fj2}
 \begin{array}{l}
 f_j^1(z,w,\sigma)=\left(\begin{array}{c}
 \psi_1(0) \langle G_j(\sigma,\Phi(\theta) z_x+w(\theta)),\beta_{k_1}\rangle\\
     \overline{\psi}_1(0)\langle G_j(\sigma,\Phi(\theta) z_x+w(\theta)),\beta_{k_1}\rangle\\
      \psi_3(0) \langle G_j(\sigma,\Phi(\theta) z_x+w(\theta)),\beta_{k_2}\rangle\\
       \overline{\psi}_3(0) \langle G_j(\sigma,\Phi(\theta) z_x+w(\theta)),\beta_{k_2}\rangle\\
 \end{array} \right),\\
 f_j^2=(I-\pi)X_0 G_j(\sigma,\Phi(\theta) z_x+w(\theta)),
 \end{array}
 \end{equation}
 Similar as  in Ref. \cite{Faria}, and using the notations in it,
 define the operator $M_j=(M_j^1,M_j^2)$, $j\geq 2$ by
 \begin{equation}\label{mjU}
 \begin{array}{ll}
 M_j^1:V_j^{4+2}(\mathbb{C}^4)\rightarrow V_j^{4+2}(\mathbb{C}^4),&(M_j^1p)(z,\sigma)=D_zp(z,\sigma)Bz-Bp(z,\sigma),\\
 M_j^2:V_j^{4+2}(\mathcal{Q}_1)\subset V_j^{4+2}({\rm Ker\pi})\rightarrow  V_j^{4+2}({\rm Ker\pi}),&(M_j^2h)(z,\sigma)=D_zh(z,\sigma)Bz-A_1h(z,\sigma).\\
 \end{array}
 \end{equation}
When  $\omega_1:\omega_2\neq m:n$ for $m,n\in\mathbb{N}$ and $1\leq m,n\leq 3$ (i.e., the case of  nonresonant double Hopf bifurcation)  it is easy to verify that
 \begin{equation*}
 \begin{aligned}
   &{\rm Im}(M_2^1)^c={\rm span}\left\lbrace \sigma_1z_1e_1,\sigma_2z_1e_1,\sigma_1z_2e_2,\sigma_2z_2e_2,\sigma_1z_3e_3,\sigma_2z_3e_3,\sigma_1z_4e_4,\sigma_2z_4e_4
   \right\rbrace.\\
 &{\rm Im}(M_3^1)^c={\rm span}\left\{
 z_1^2z_2e_1 ,
 z_1z_3z_4e_1 ,
 z_1z_2^2e_2 ,
 z_2z_3z_4e_2,
 z_3^2z_4e_3
 z_1z_2z_3e_3 ,
 z_3z_4^2e_4,
 z_1z_2z_4e_4
 \right\}.
 \end{aligned}
\end{equation*}
where $e_1=(1,0,0,0)^T,e_2=(0,1,0,0)^T,e_3=(0,0,1,0)^T,e_4=(0,0,0,1)^T$.
 According to Ref. \cite{Faria},  by a recursive transformations of variables of the following form
 \begin{equation}
 \label{zyalpha}
 (z,w,\sigma)=(\widehat{z},\widehat{w},\sigma)+\frac{1}{j!}(U_j^1(\widehat{z},\sigma),U_j^2(\widehat{z},\sigma),0),
 \end{equation}
 with $U_j=(U_j^1,U_j^2)\in V_j^{4+2}(\mathbb{C}^4)\times V_j^{4+2}(\mathcal{Q}_1)$,  \cite{Faria}
 which  transforms (\ref{zdotTaylor}) into the normal form
 \begin{equation}\label{zdotnormalform}
 \begin{array}{l}
 \dot{z}=Bz+\sum\limits_{j\geq 2}\frac{1}{j!}g_j^1(z,w,\sigma),\\
 \frac{dw}{dt} =A_1w+\sum\limits_{j\geq 2}\frac{1}{j!}g_j^2(z,w,\sigma),
 \end{array}
 \end{equation}
 where $g_j=(g_j^1,g_j^2), j\geq 2$, are  given by
 $g_j(z,w,\sigma)=\overline{f}_j(z,w,\sigma)-M_jU_j(z,\sigma)$,
 and $U_j\in V_j^{4+2}(\mathbb{C}^4)\times V_j^{4+2}(\mathcal{Q}_1)$
 are expressed as
 \begin{equation}\label{Ujzalpha}
 U_j(z,\sigma)=(M_j)^{-1}{\rm Proj}_{{\rm Im}(M_j^1)\times{\rm Im}(M_j^2)}\circ\overline{f}_j(z,0,\sigma),
 \end{equation}
 where $\overline{f}_j=(\overline{f}_j^1,\overline{f}_j^2)$ stand for the terms of order $j$ in $(z,w)$, which are obtained after the computation of normal forms up to order $j-1$.
 From Ref. \cite{Faria}, the normal form truncated to the third order has the following form
 \begin{equation}
 \dot{z}=Bz+\frac{1}{2!}g_2^1(z,0,\sigma)+\frac{1}{3!}g_3^1(z,0,0)+h.o.t..
 \end{equation}
Here $g_3^1(z,0,0)={\rm Proj}_{{\rm Ker}(M_3^1)}\overline{f}_3^1(z,0,0)$, where
 \begin{equation}
 \label{f31bar}
 \begin{aligned}
 \overline{f}_3^1(z,0,0)=f_3^1(z,0,0)+\frac{3}{2}[D_zf_2^1(z,0,0)U_2^1(z,0)\\+D_wf_2^1(z,0,0)U_2^2(z,0)-D_zU_2^1(z,0)g_2^1(z,0,0)],
 \end{aligned}
 \end{equation}
and $(U_2^1(z,\sigma),U_2^2(z,\sigma))\in V_j^{4+2}(\mathbb{C}^4)\times V_j^{4+2}(\mathcal{Q}_1)$ given by (\ref{Ujzalpha}).

 The calculations of $g_2^1(z,0,\sigma)$ and $g_3^1(z,0,0)$ heavily depends on tedious mathematical derivations such as solve the center manifold function and some projections, thus we  leave the details  in section 2 and 3
  of   the supplement materials.

\begin{figure}
  \centering
  \includegraphics[height=0.44\textwidth]{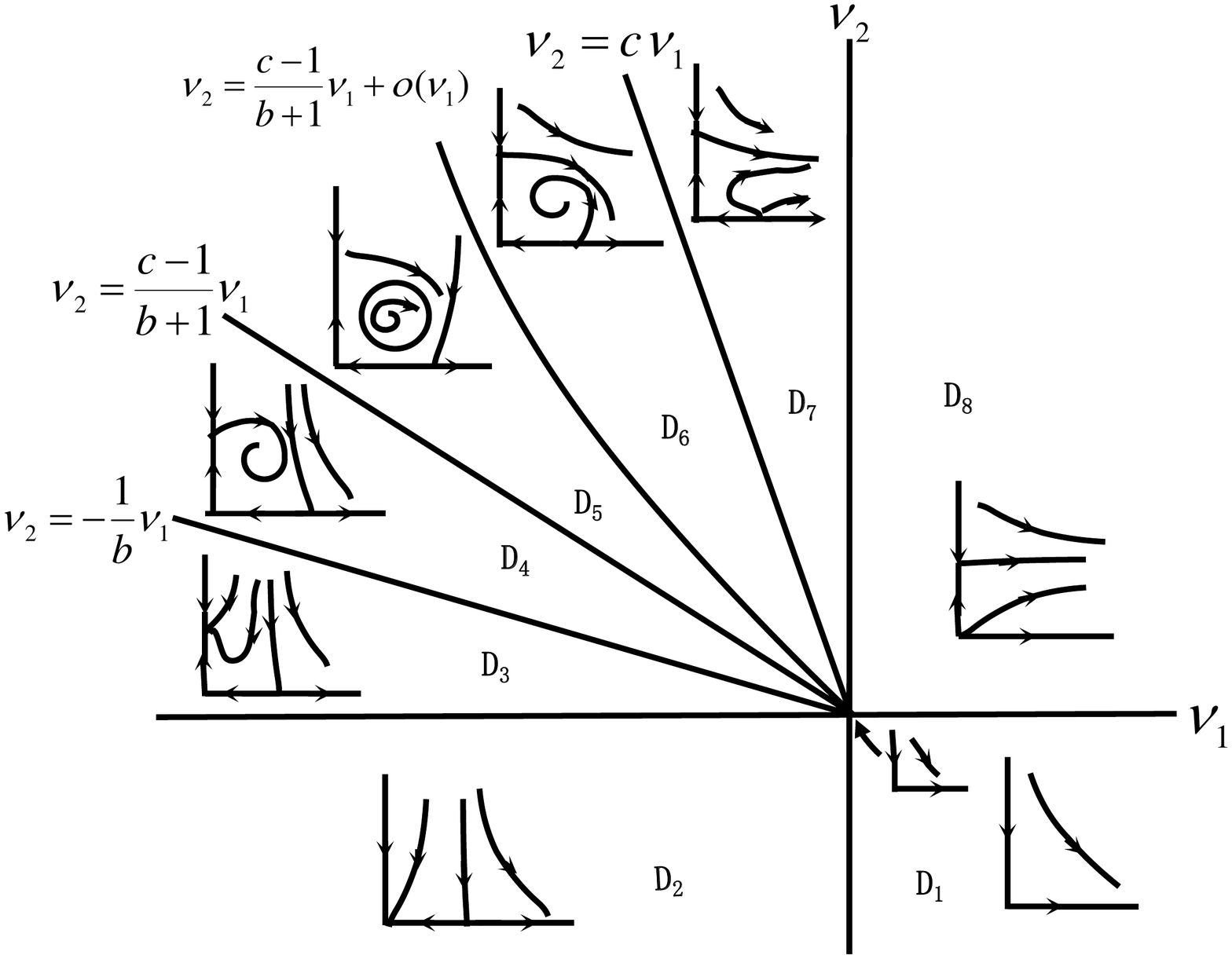}
  \caption{Phase portraits for  the unfoldings of case VIa with $\epsilon_1=1$.}
  \label{fig:VIa}
   \end{figure}

Finally, the normal form truncated to the third order on the center manifold for double Hopf bifurcation is obtained  as follows
\begin{equation}\label{normalform}
\begin{aligned}
&\dot{z}_1=i\omega_1\tau_1^*z_1+K_{11}\sigma_1z_1+K_{21}\sigma_2z_1+ K_{2100}z_1^2z_2+K_{1011}z_1z_3z_4, \\
&\dot{z}_2=-i\omega_1\tau_1^*z_2+\overline{K_{11}}\sigma_1z_2+\overline{K_{21}}\sigma_2z_2+ \overline{K_{2100}}z_1z_2^2+\overline{K_{1011}}z_2z_3z_4,\\
&\dot{z}_3=i\omega_2\tau_1^*z_3+K_{13}\sigma_1z_3+K_{23}\sigma_2z_3+ K_{0021}z_3^2z_4+K_{1110}z_1z_2z_3 ,\\
&\dot{z}_4=-i\omega_2\tau_1^*z_4+\overline{K_{13}}\sigma_1z_4+\overline{K_{23}}\sigma_2z_4+ \overline{K_{0021}}z_3z_4^2+\overline{K_{1110}}z_1z_2z_4 .
\end{aligned}
\end{equation}
Make the  polar  coordinate transformation
\[\begin{array}{ll}
z_1=\rho_1\cos \theta_1+i\rho_1\sin\theta_1,& z_2=\rho_1\cos \theta_1-i\rho_1\sin\theta_1,\\z_3=\rho_2\cos \theta_2+i\rho_2\sin\theta_2,& z_4=\rho_2\cos \theta_2-i\rho_2\sin\theta_2,
\end{array}\]
where $\rho_1,\rho_2>0$. Denote $\epsilon_1={\rm Sign}({\rm Re}K_{2100})$,  $\epsilon_2={\rm Sign}({\rm Re}K_{0021})$,  rescale as $\widehat{\rho}_1=\rho_1\sqrt{|{\rm}K_{2100}|}$, $\widehat{\rho}_2=\rho_2\sqrt{|{\rm}K_{0021}|}$, $\widehat{t}=t\epsilon_1$, and drop the hats, then we obtain the equivalent system of (\ref{normalform})
 \begin{equation}\label{normalformcylin}
 \begin{aligned}
 &\dot{\rho}_1=\rho_1(\nu_1+\rho_1^2+b\rho_2^2),\\
 &\dot{\rho}_2=\rho_2(\nu_2+c\rho_1^2+d\rho_2^2).
 \end{aligned}
 \end{equation}
Here
 \[\begin{aligned}
& \nu_1=\epsilon_1({\rm Re}K_{11}\sigma_1+{\rm Re}K_{21}\sigma_2),
 \nu_2=\epsilon_1({\rm Re}K_{13}\sigma_1+{\rm Re}K_{23}\sigma_2),  \\&b=\frac{\epsilon_1\epsilon_2{\rm Re}K_{1011}}{{\rm Re}K_{0021}}, c=\frac{{\rm Re}K_{1110}}{{\rm Re}K_{2100}},d=\epsilon_1\epsilon_2.\\
 \end{aligned}\]
As was discussed in  chapter 7.5 in Ref. \cite{Guckenheimer},  there are twelve distinct kinds of unfoldings for Eq. (\ref{normalformcylin}) (see Table 1).
 \begin{table}[tbp]
 \label{twelve}
 \caption{The twelve unfoldings of system (\ref{normalformcylin}).}
 \centering
 \begin{tabular}{lcccccccccccc}
 \hline
 Case & {I}a & {I}b  & {II} & {III}& {IV}a & {IV}b  & {V}& {VI}a& {VI}b& {VII}a& {VII}b& {VIII}     \\ \hline  
 $d$ &  +1 &+1& +1&+1&+1&+1&--1&--1&--1&--1&--1&--1\\         
 $b $ &  + &+& +&--&--&--&+&+&+&--&--&--\\        
 $c$ &  + &+& --&+&--&--&+&--&--&+&+&--\\
 $d-bc$ &  + &--& +&+&+&--&--&+&--&+&--&--\\
 \hline
 \end{tabular}
 \end{table}
 In section \ref{simulations}, case VIa arises, thus we draw bifurcation set and phase portraits for the unfolding of case VIa in Fig. \ref{fig:VIa}.

\section{Numerical simulations}
\label{simulations}

In this section,  we carry out some simulations. We choose
      \begin{equation*}
      \label{para1}
          r_1=0.8,r_2=1,a=1.3,K=0.7,\gamma=1,m=0.27,l=2,
      d_1=0.3,d_2=0.4,
        \end{equation*}
 then one  can get the unique positive constant equilibrium $E^*(0.4358,0.3181)$, which is globally asymptotically stable when $\tau_1=\tau_2=0$ according to Remark \ref{Du}.

\begin{figure}
     \centering
  a) \includegraphics[width=0.27\textwidth]{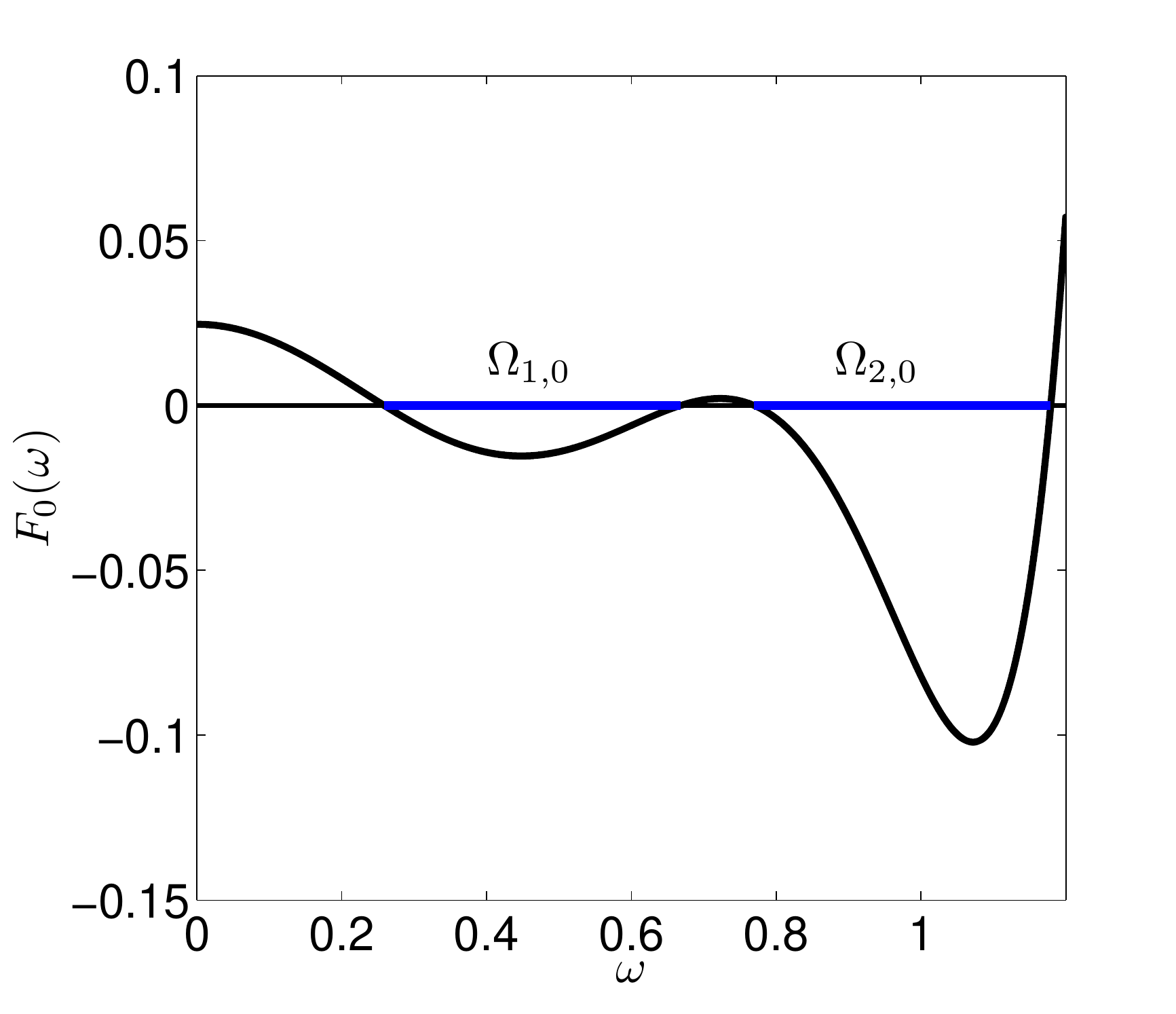}
    b)   \centering
  \includegraphics[width=0.27\textwidth]{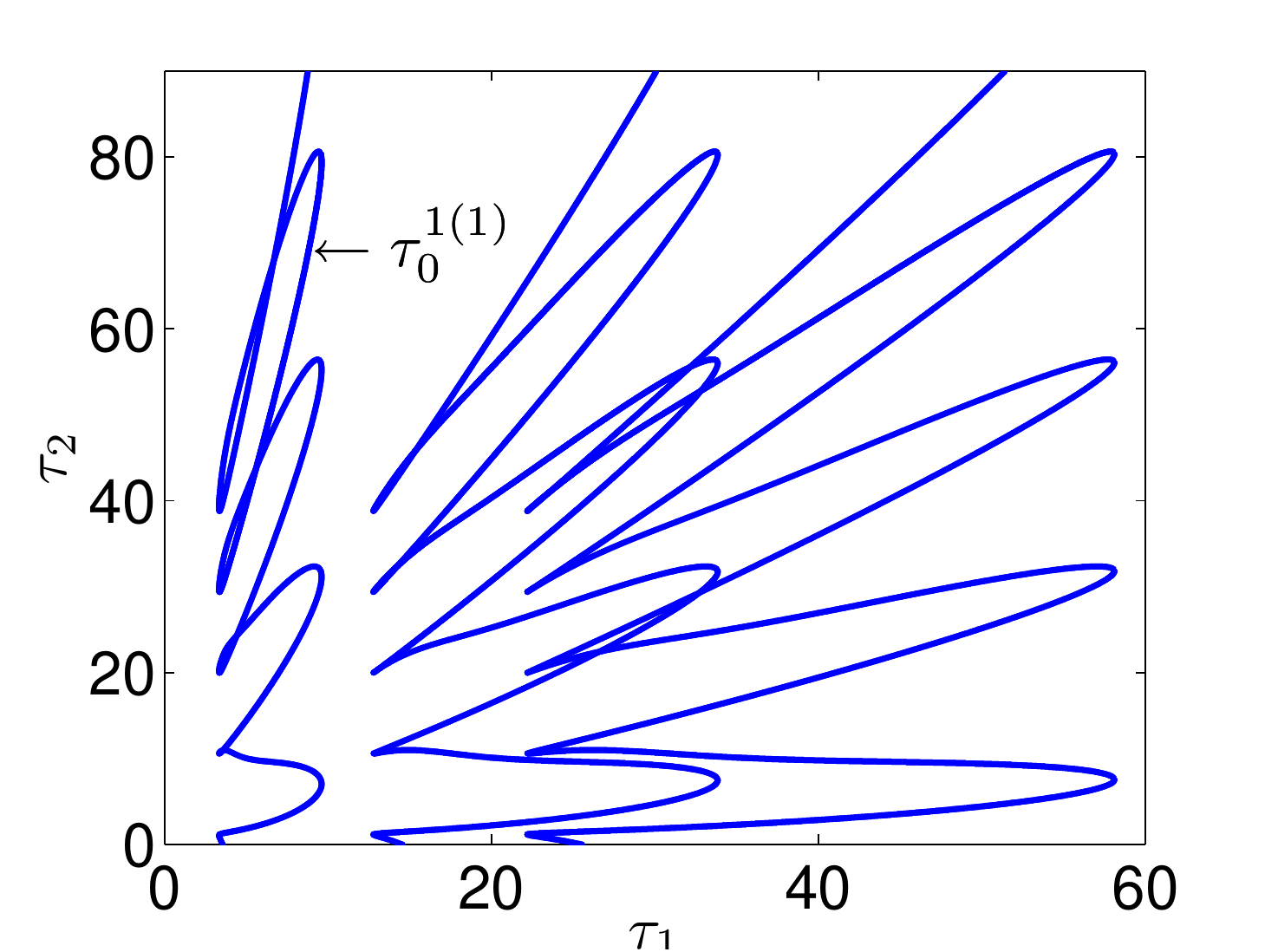}
  c) \includegraphics[width=0.27\textwidth]{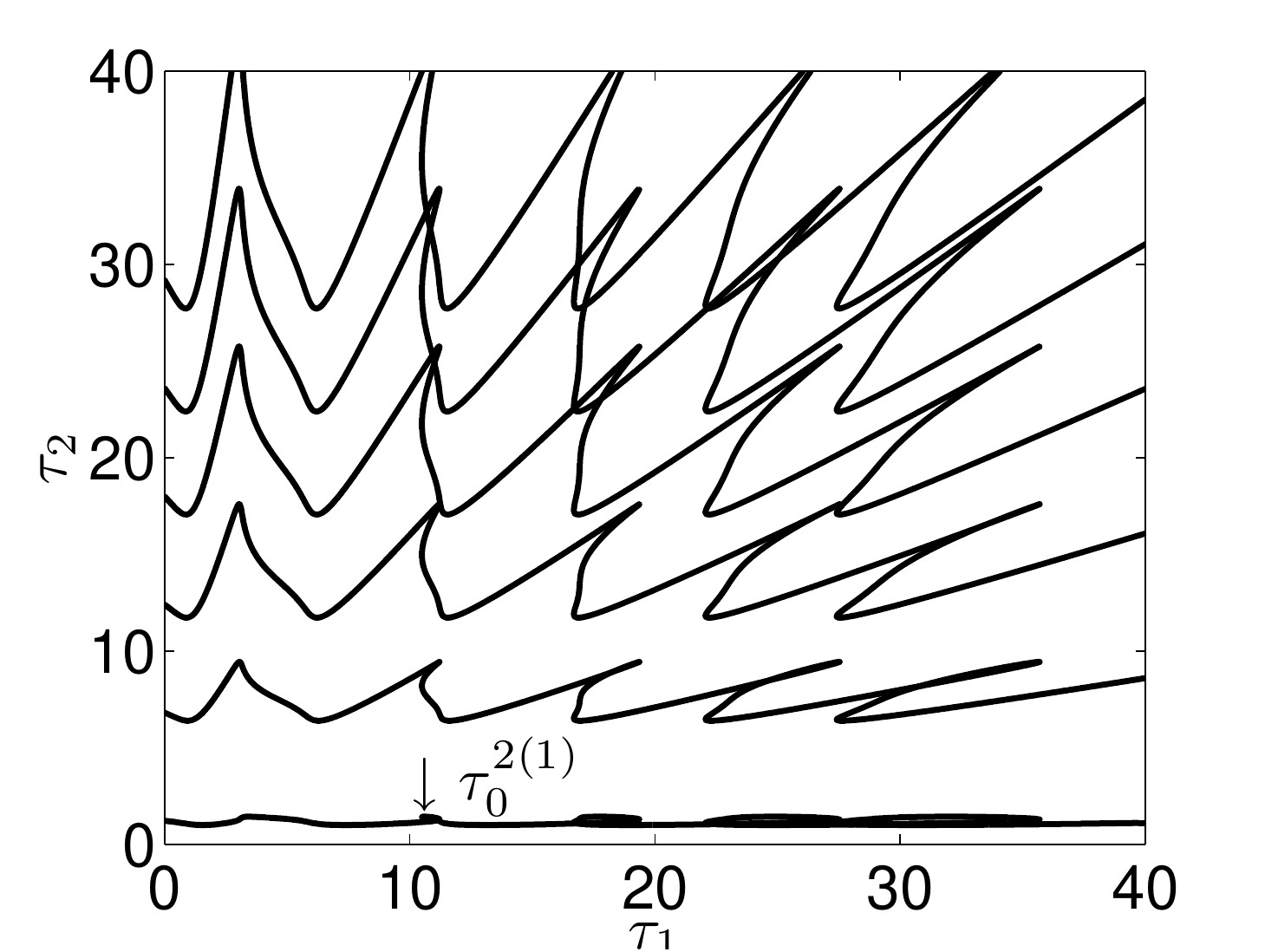}
 \caption{a) Graph of  $F_0(\omega)$. b) Stability switching curves $\mathcal{T}^1_0$. c) Stability switching curves $\mathcal{T}^2_0$.}
  \label{fig:F0T0}
 \end{figure}

\begin{figure}
     \centering
  a) \includegraphics[width=0.46\textwidth]{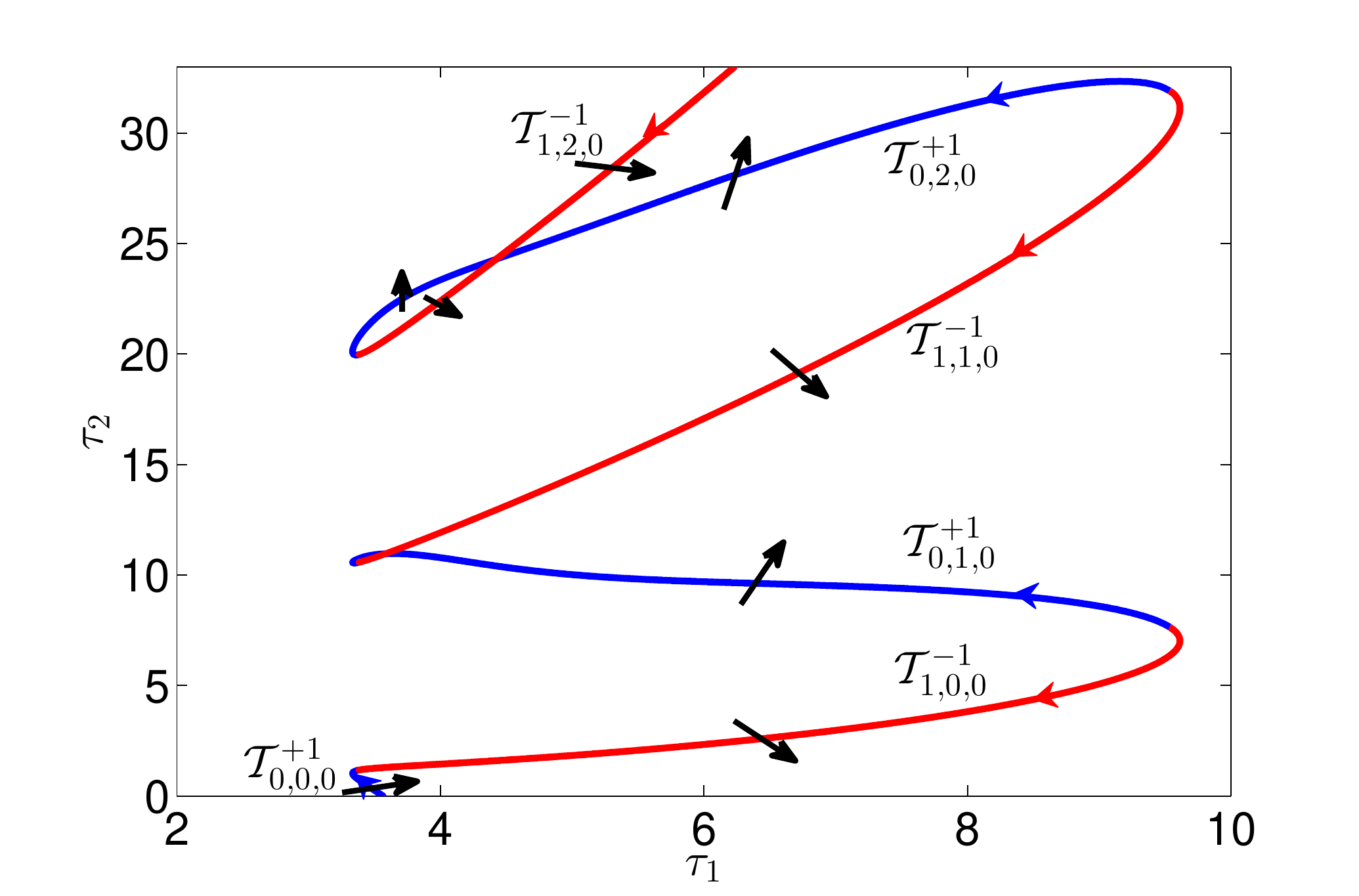}
  b) \includegraphics[width=0.46\textwidth]{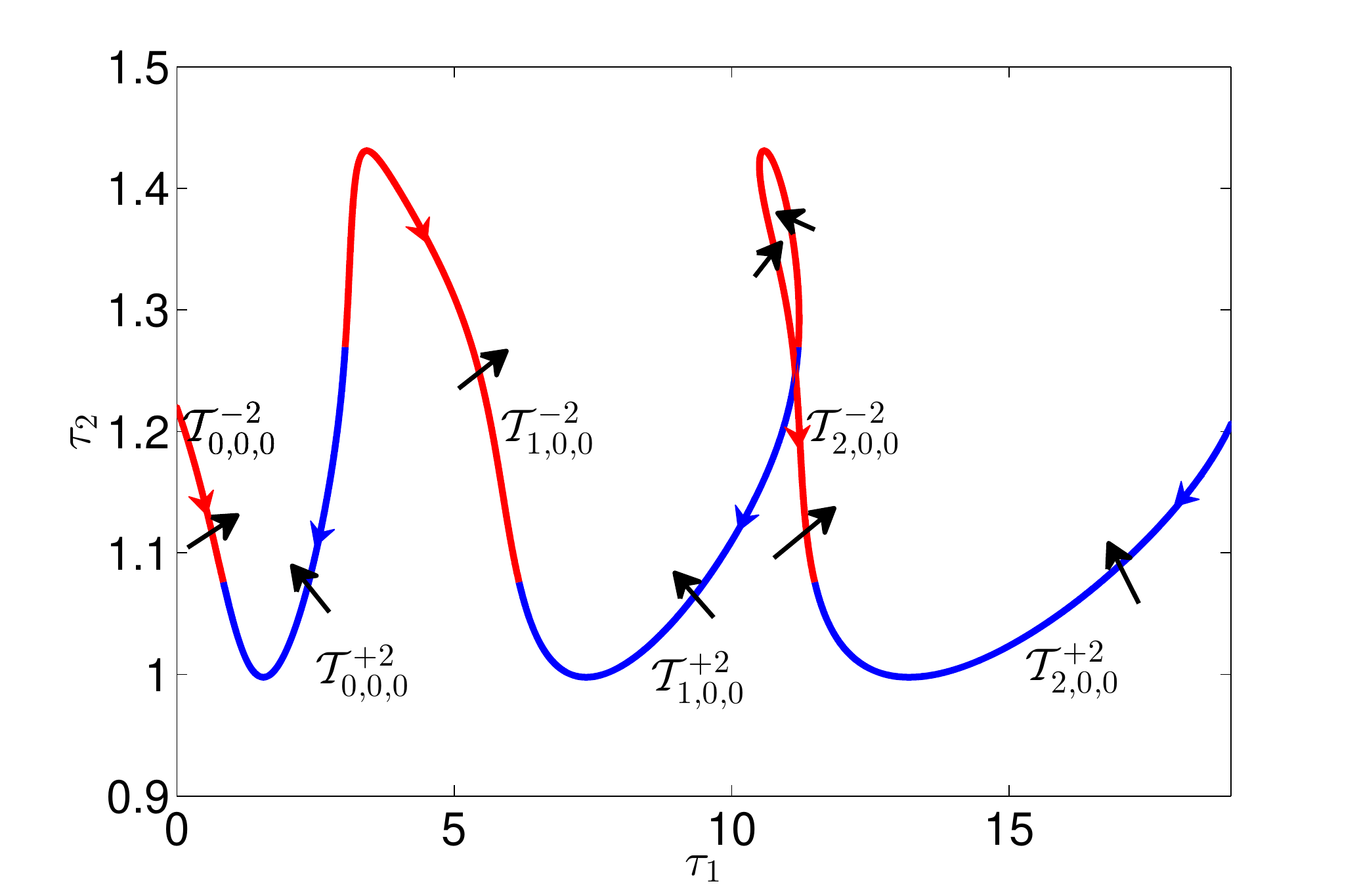}
 \caption{a) The detailed structure of the left-most curve of $\mathcal{T}_0^1$ (marked $\tau_0^{1(1)}$ in Fig. \ref{fig:F0T0} b)). The blue (red) arrow represents the positive direction of  $\mathcal{T}_{j_1,j_2,0}^{+1}$ ($\mathcal{T}_{j_1,j_2,0}^{-1}$). b) The detailed structure of the lowest curve of $\mathcal{T}_0^2$ (marked $\tau_0^{2(1)}$ in Fig. \ref{fig:F0T0} c)). The blue (red) arrow represents the positive direction of  $\mathcal{T}_{j_1,j_2,0}^{+2}$ ($\mathcal{T}_{j_1,j_2,0}^{-2}$). From Lemma \ref{direction}, we know that the regions on the right (left) of the blue (red) curves, which the black arrows point to, have two more characteristic roots with positive real parts.}
  \label{fig:T0direction}
 \end{figure}

To illustrate the dynamics in the presence of delays,
we follow the process given in section 2.1. As shown in Fig. \ref{fig:F0T0} a), $F_0(0)>0$, and $F_0(\omega)=0$ has four roots $a_{1,0}=0.2587,b_{1,0}=0.6682,a_{2,0}=0.7697$, and $b_{2,0}=1.1791$. The crossing set is  $\Omega_{1,0}\bigcup\Omega_{2,0}=[a_{1,0},b_{1,0}]\bigcup[a_{2,0},b_{2,0}]$.  For the two ends of the $\Omega_{1,0}$,
 we have $\theta_{1,0}(a_{1,0})=\pi$, $\theta_{2,0}(a_{1,0})=\pi$, $\theta_{1,0}(b_{1,0})=\pi$, and $\theta_{2,0}(b_{1,0})=0$, and $\delta_1^a=1$, $\delta_2^a=1$, $\delta_1^b=1$, $\delta_2^b=0$. From (\ref{Tzfk}) and (\ref{Tk}), we can get the stability switching curves $\mathcal{T}^1_0$ corresponding to $\Omega_{1,0}$ which is shown in Fig. \ref{fig:F0T0}  b).  From the previous discussion in (\ref{connect}), $\mathcal{T}_{j_1,j_2,0}^{+1}$   is connected to  $\mathcal{T}_{j_1+1,j_2-1,0}^{-1}$  at the left point $a_{1,0}$ and $\mathcal{T}_{j_1+1,j_2,0}^{-1}$ at the right point $b_{1,0}$ for any $j_1,j_2$. To show the structure of stability switching curves and the crossing direction clearly, we take the left-most curve of $\mathcal{T}_0^1$ (i.e. $\tau_0^{1(1)}$ in  Fig. \ref{fig:F0T0}  b)) as an example, and draw the figure in Fig. \ref{fig:T0direction} a). From bottom to top, it starts with a part of $\mathcal{T}_{0,0,0}^{+1}$, which is connected to $\mathcal{T}_{1,0,0}^{-1}$ at $b_{1,0}$.  $\mathcal{T}_{1,0,0}^{-1}$ is linked to $\mathcal{T}_{0,1,0}^{+1}$ at $a_{1,0}$, which is again connected to $\mathcal{T}_{1,1,0}^{-1}$ at $b_{1,0}$ $\cdots\cdots$.       The numerical results coincide with the analysis result in (\ref{connect}).
 In fact, the rest curves  of $\mathcal{T}_0^1$ in Fig. \ref{fig:F0T0}  b)) have  similar  structure as  $\tau_0^{1(1)}$.
 Similarly,   the stability switching curves $\mathcal{T}^2_0$ corresponding to $\Omega_{2,0}$  are shown in Fig. \ref{fig:F0T0}  c), and    the lowest curve of $\mathcal{T}_0^2$ (marked $\tau_0^{2(1)}$ in Fig. \ref{fig:F0T0} c)) is drawn in  Fig. \ref{fig:T0direction} b).  All the stability switching curves for $n=0$ are given by $\mathcal{T}_0=\mathcal{T}^1_0\cup\mathcal{T}^2_0$.

\begin{figure}
 \centering
 a) \includegraphics[width=0.27\textwidth]{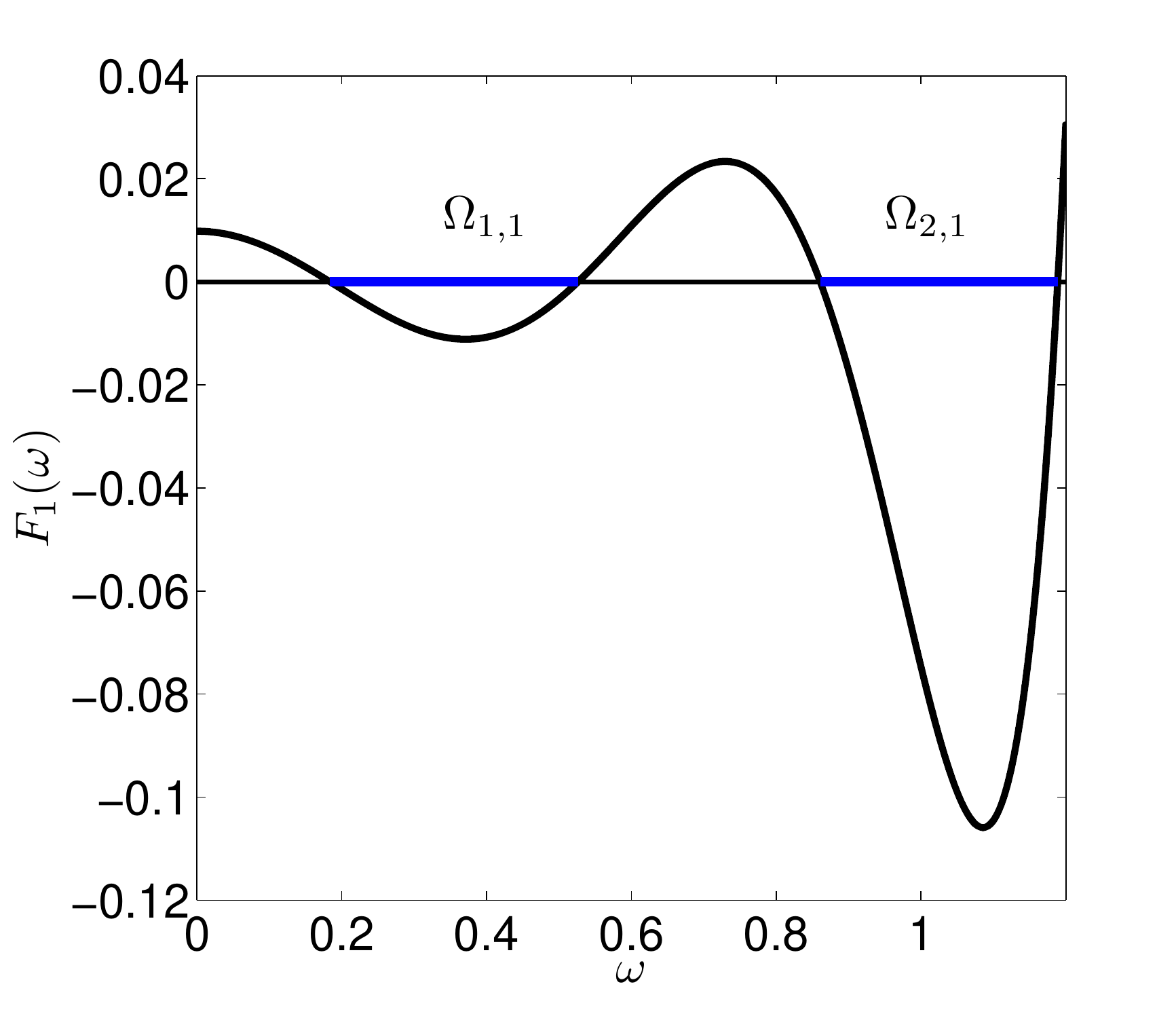}
   b)   \centering
  \includegraphics[width=0.27\textwidth]{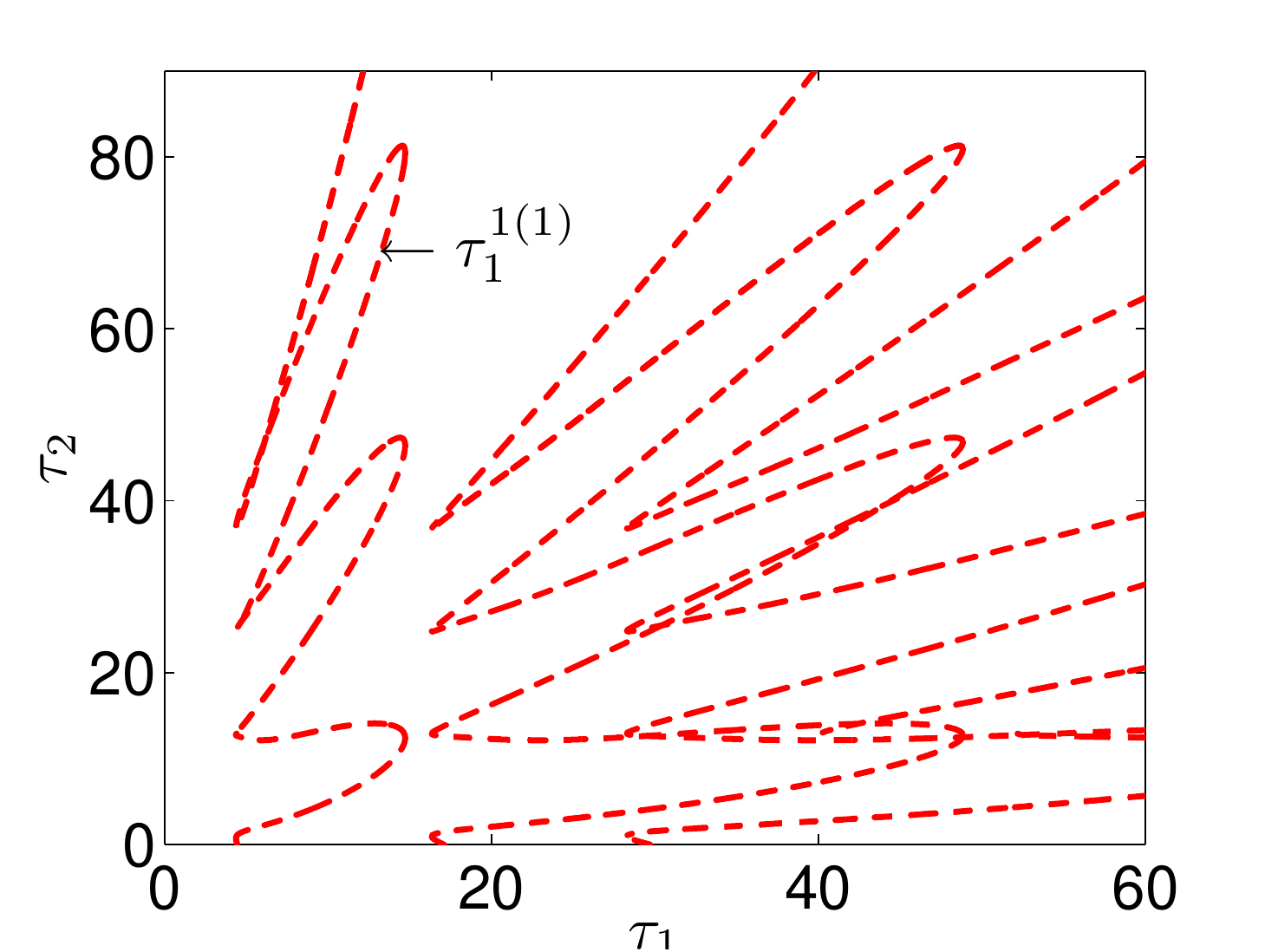}
  c) \includegraphics[width=0.27\textwidth]{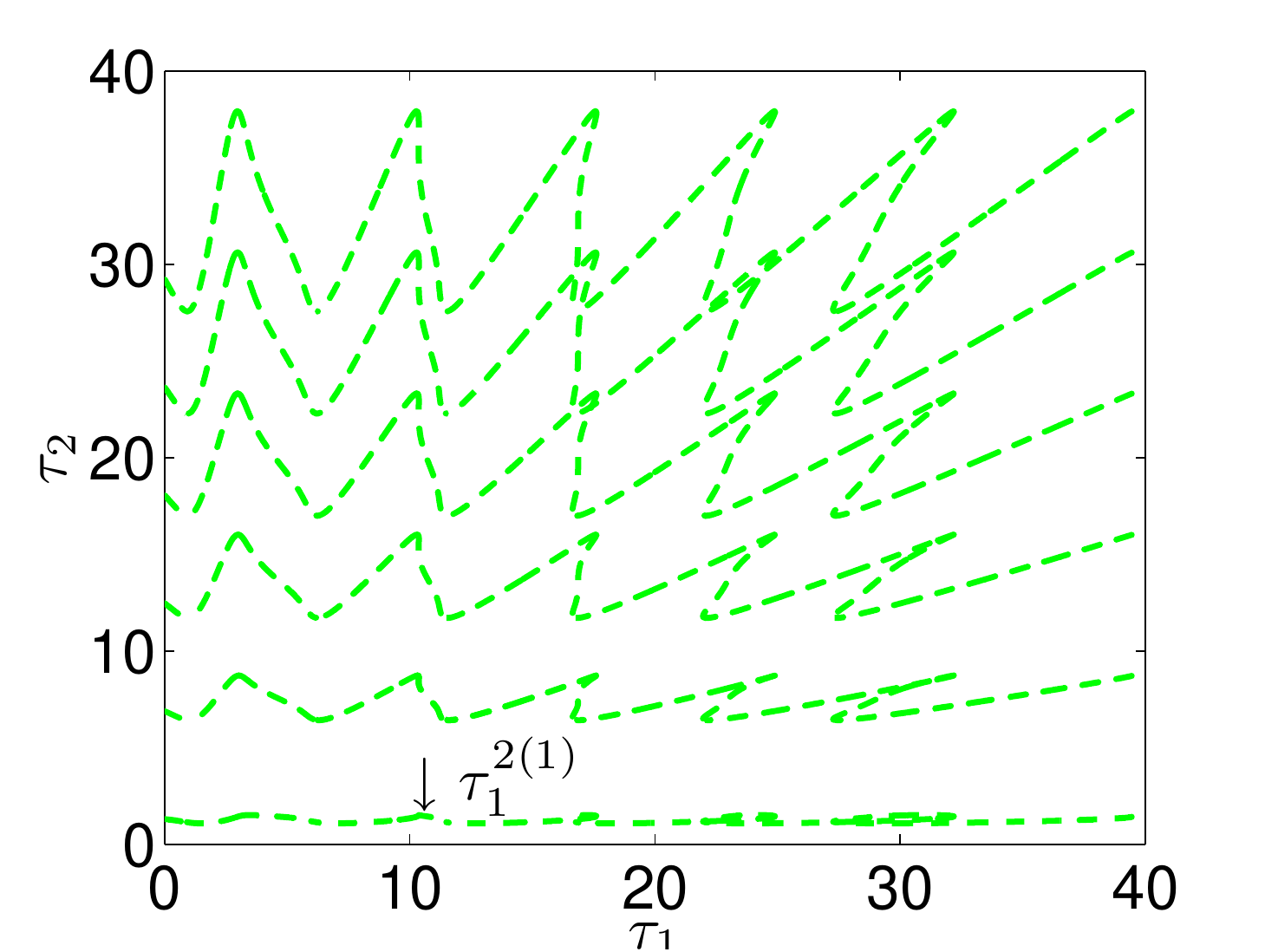}
 \caption{a) Graph of  $F_1(\omega)$. b) Stability switching curves $\mathcal{T}^1_1$. c) Stability switching curves $\mathcal{T}^2_1$.}
   \label{fig:F1T1}
       \end{figure}

When $n=1$, $F_1(0)>0$, and $F_1(\omega)=0$ has four roots $a_{1,1}=0.184,b_{1,1}=0.5264,a_{2,1}=0.8607,b_{2,1}=1.189$, which is shown in Fig. \ref{fig:F1T1} a). Thus, the crossing set is  $\Omega_{1,1}\bigcup\Omega_{2,1}=[a_{1,1},b_{1,1}]\bigcup[a_{2,1},b_{2,1}]$, and  we can get the stability switching curves $\mathcal{T}^1_1$ and $\mathcal{T}^2_1$, which is shown in Fig. \ref{fig:F1T1} b) and c).  Thus all the stability switching curves for $n=1$  are given by $\mathcal{T}_1=\mathcal{T}^1_1\cup\mathcal{T}^2_1$.

\begin{figure}
     \centering
 a) \includegraphics[width=0.27\textwidth]{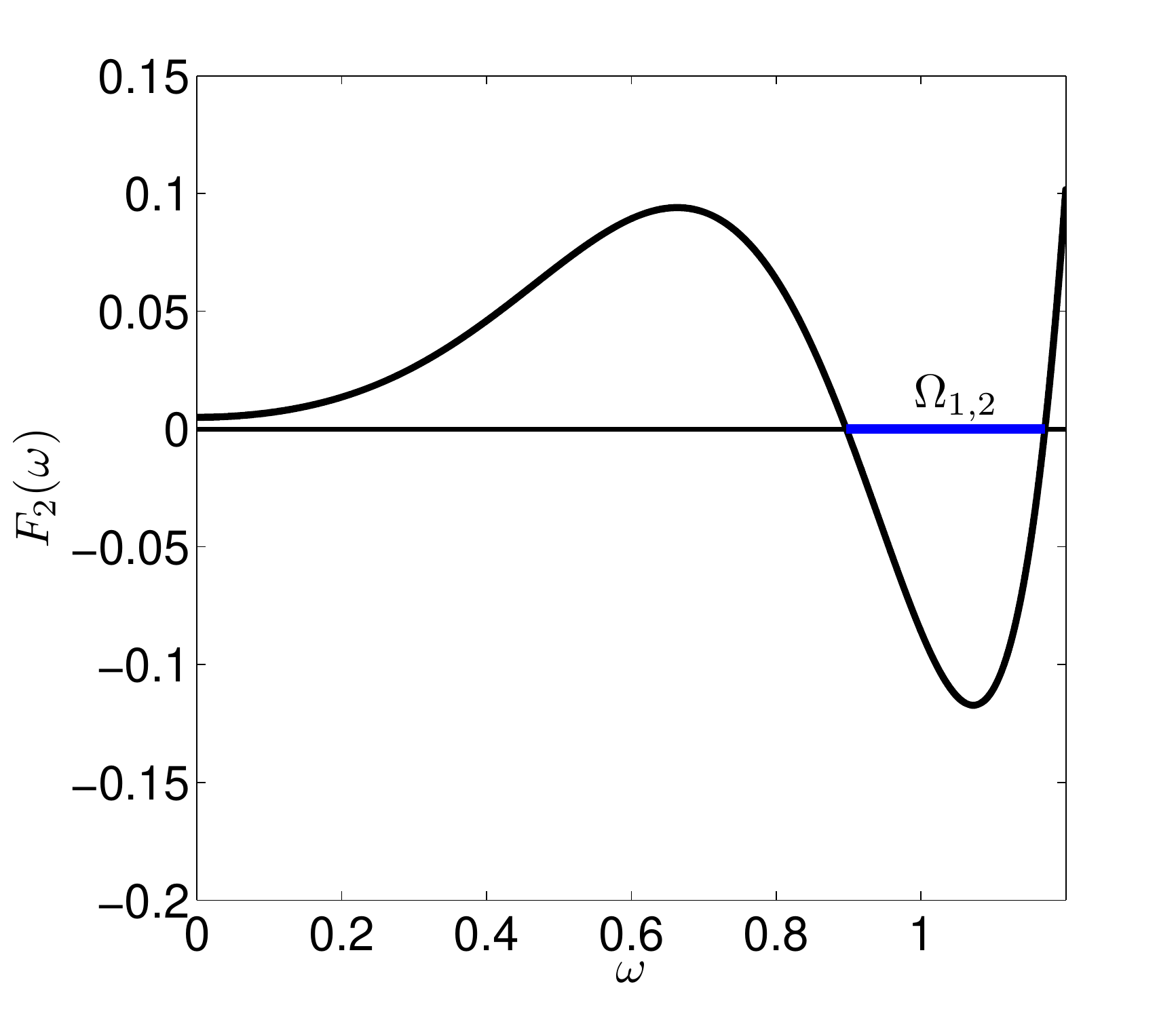}
  b)   \centering
  \includegraphics[width=0.27\textwidth]{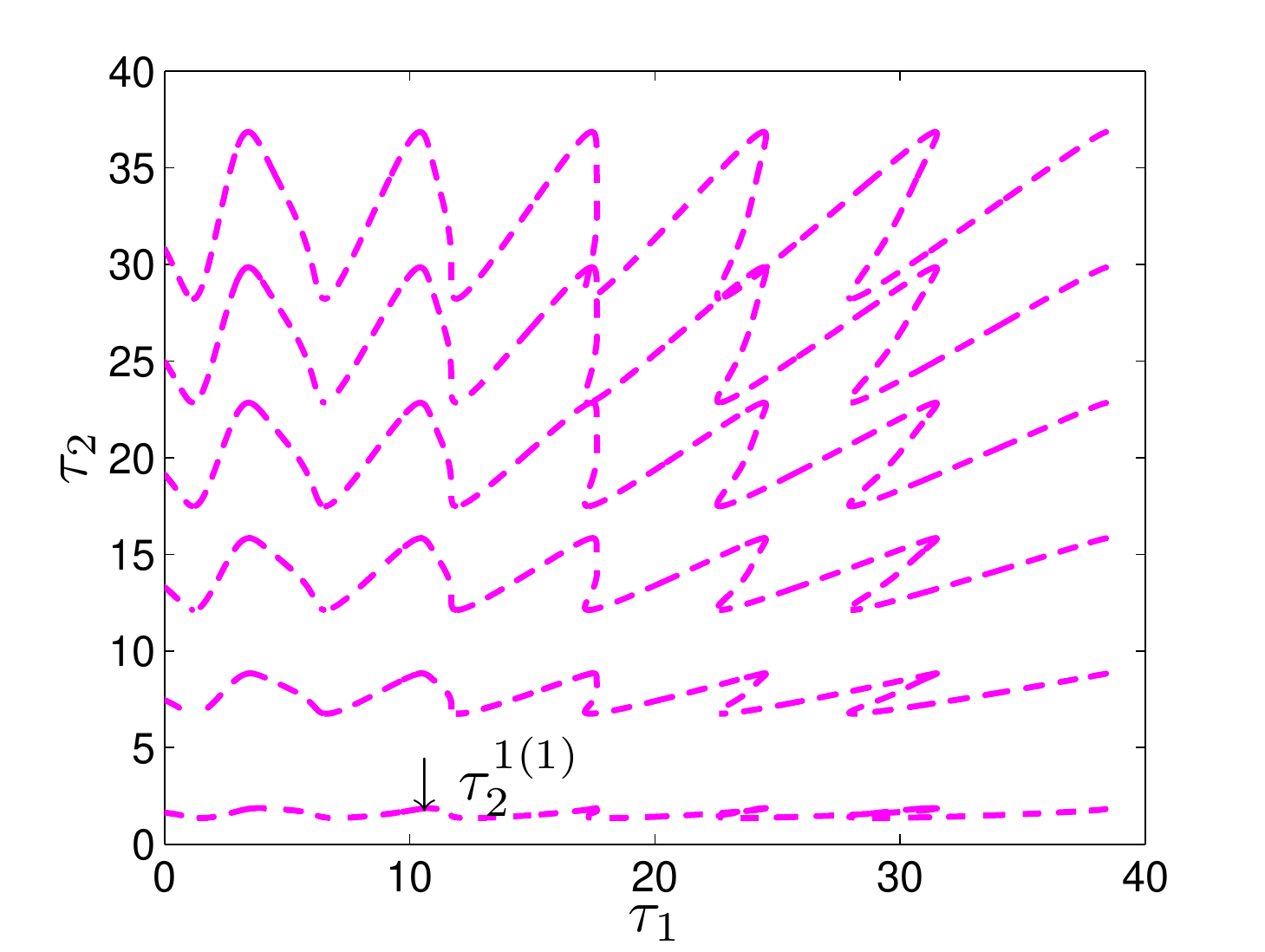}
  \caption{a) Graph of $F_2(\omega)$. b) Stability switching curves  $\mathcal{T}^1_2$. }
           \label{fig:F2T2}
        \end{figure}

When $n=2$, $F_2(0)>0$, and $F_2(\omega)=0$ has two roots $a_{1,2}=0.8968,b_{1,2}=1.171$, which is shown in Fig. \ref{fig:F2T2} a). The crossing set is  $\Omega_{1,2}=[a_{1,2},b_{1,2}]$. The stability switching curves $\mathcal{T}^1_2$ corresponding to $\Omega_{1,2}$ is shown in Fig. \ref{fig:F2T2} b).  All the stability switching curves for $n=2$ are given by $\mathcal{T}_2=\mathcal{T}^1_2$.

\begin{figure}
         \centering
  a) \includegraphics[width=0.27\textwidth]{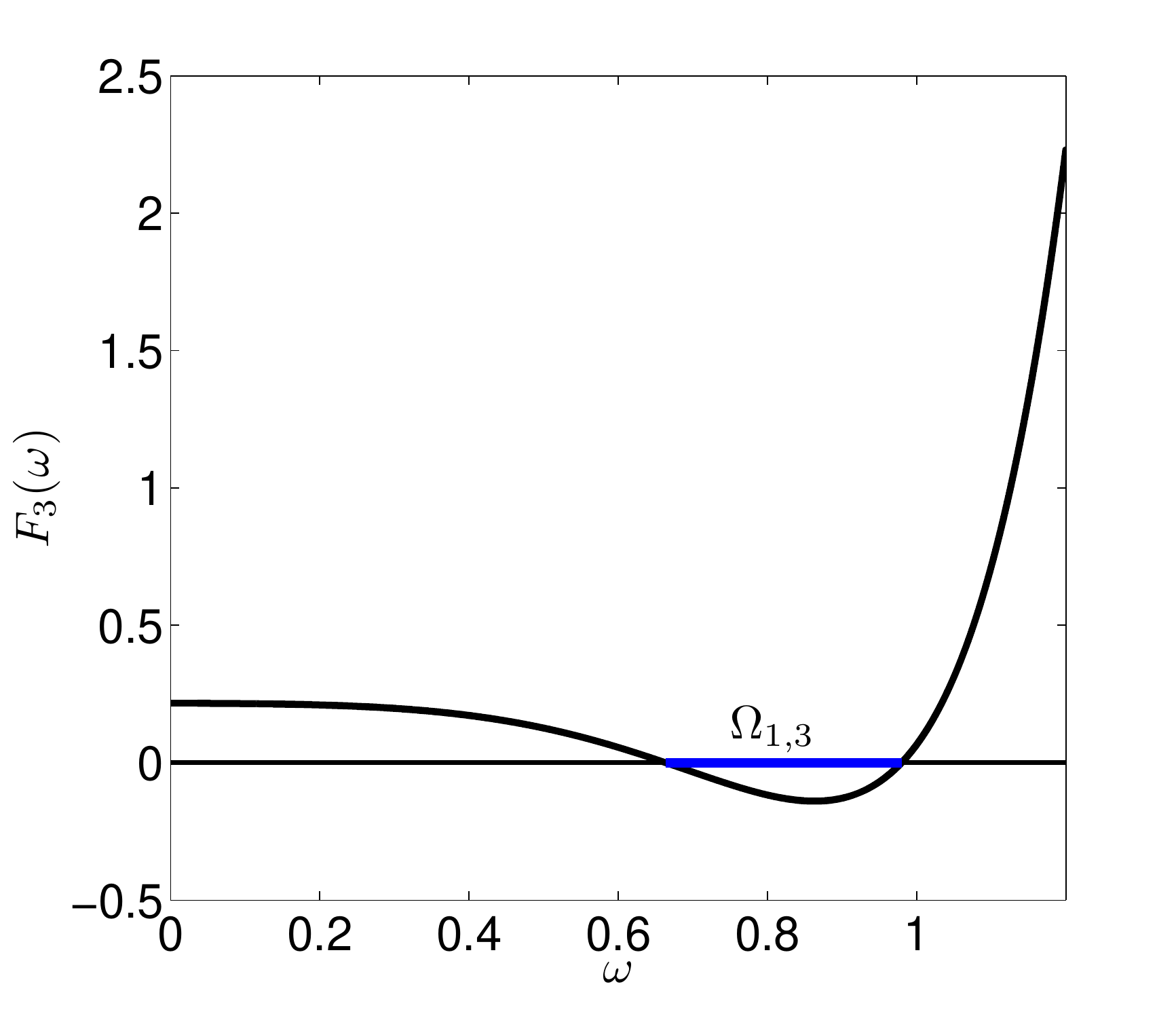}
      b)   \centering
  \includegraphics[width=0.27\textwidth]{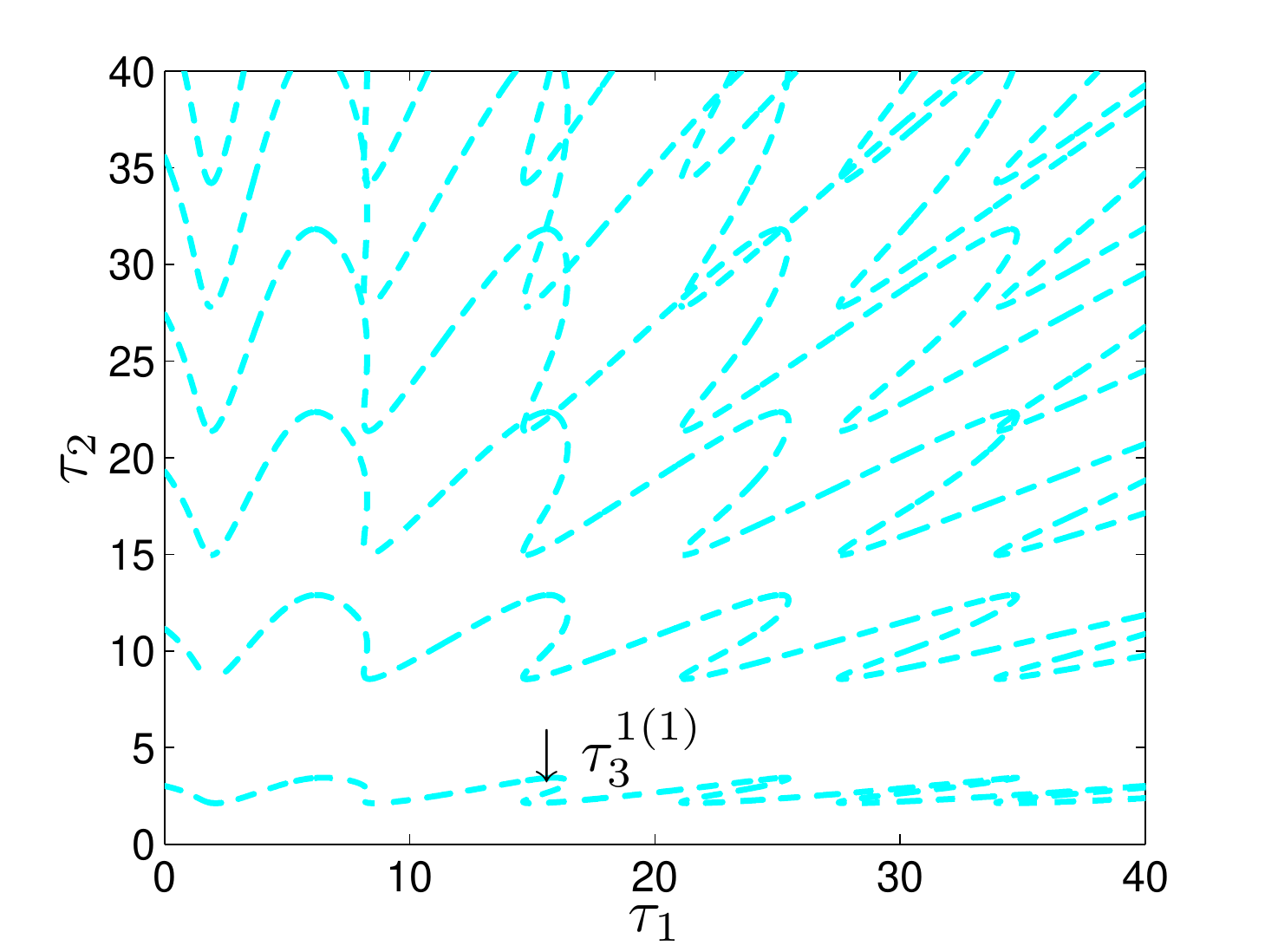}
 \caption{a) Graph of  $F_3(\omega)$. b) Stability switching curves $\mathcal{T}^1_3$. }
        \label{fig:F3T3}
         \end{figure}

When $n=3$, $F_3(0)>0$, and $F_3(\omega)=0$ has two roots $a_{1,3}=0.6638,b_{1,3}=0.9798$, which is shown in Fig. \ref{fig:F3T3} a). The crossing set is  $\Omega_{1,3}=[a_{1,3},b_{1,3}]$. And the stability switching curves $\mathcal{T}^1_3$ corresponding to $\Omega_{1,3}$ is shown in Fig. \ref{fig:F3T3} b).  Thus all the stability switching curves for $n=3$ are given by $\mathcal{T}_3=\mathcal{T}^1_3$.

When $n\geq 4$, numerical calculation indicates $F_n(\omega)>0$ for any $\omega$, thus there are no  stability switching curves on $(\tau_1,\tau_2)$ plane  for  $n\geq 4$.

 \begin{figure}
            a)   \centering
       \includegraphics[width=0.44\textwidth,height=0.32\textwidth]{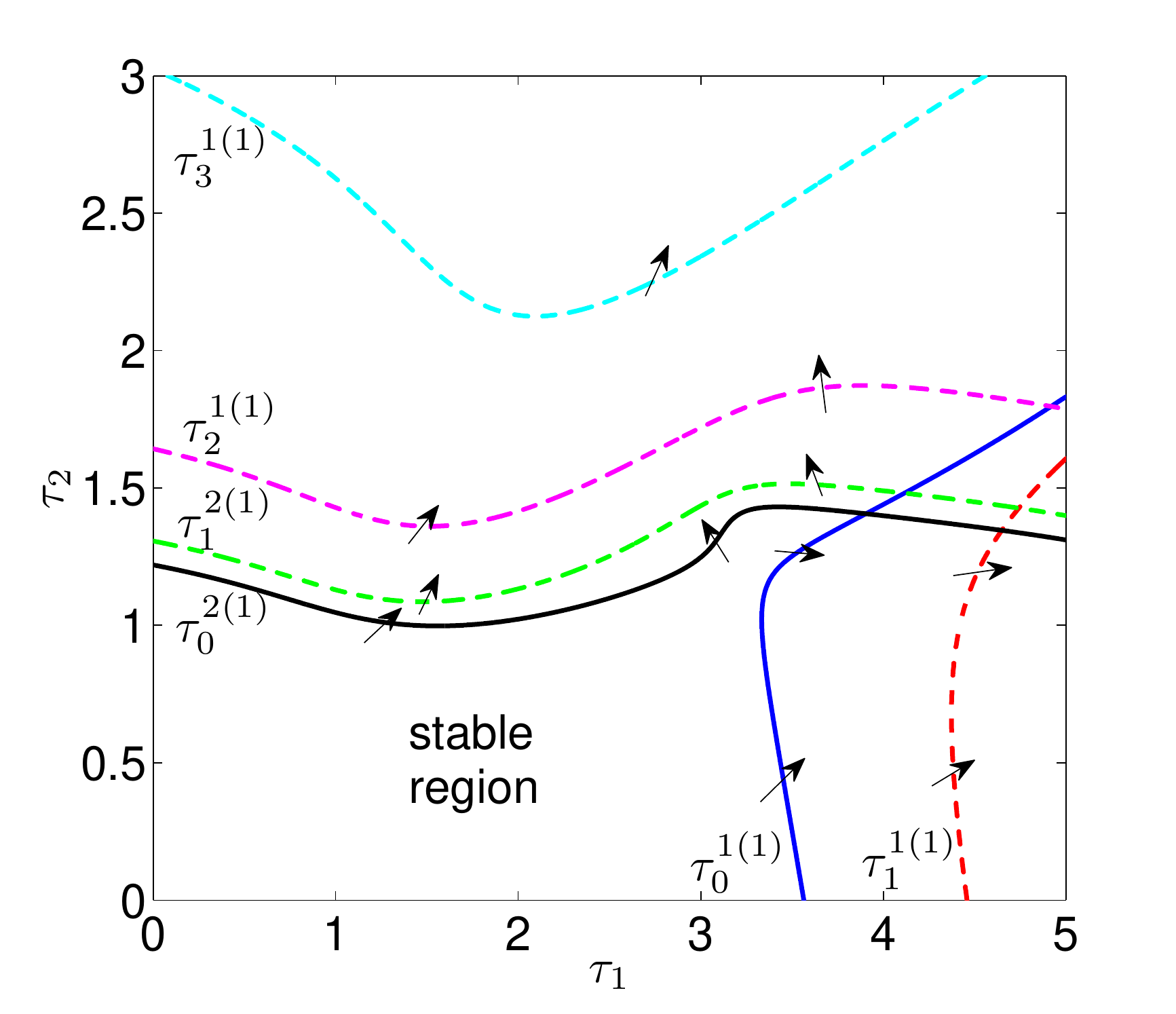}
       b) \includegraphics[width=0.44\textwidth,height=0.32\textwidth]{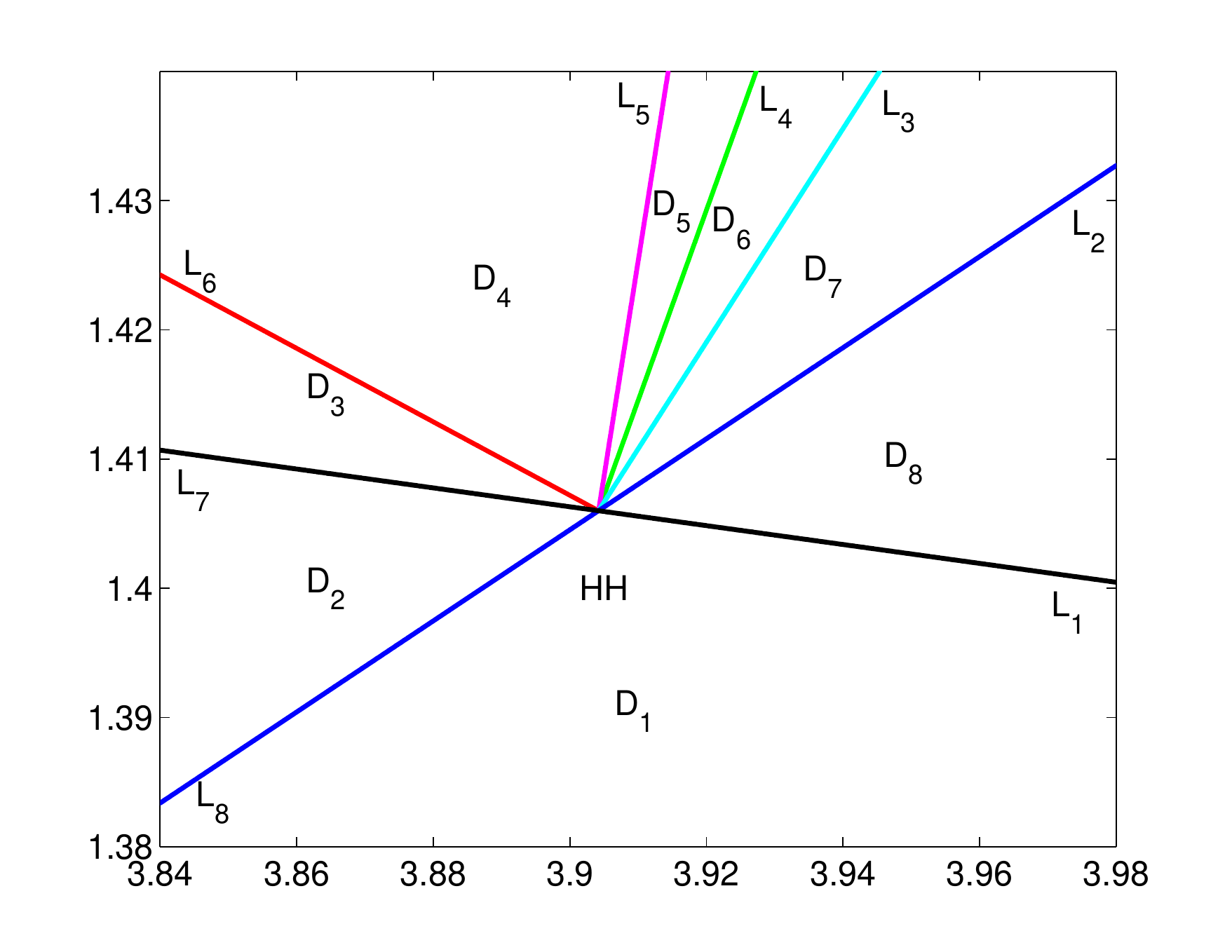}
   \caption  {a) The left-most curve and the lowest curve of $\mathcal{T}_0$ intersect at  $(\tau_1,\tau_2)=(3.9042,1.406)$, which is a double Hopf bifurcation point on the $\tau_1-\tau_2$ plane.  Crossing directions are marked by arrows.  b) The complete bifurcation sets near HH.}
    \label{fig:tau1tau2}
   \end{figure}

  Combining the stability switching curves shown in Fig. \ref{fig:F0T0}-\ref{fig:F3T3} together, and zooming in the part when $(\tau_1,\tau_2)\in [0,5]\times[0,3]$, we have the Hopf bifurcation curves shown in Fig. \ref{fig:tau1tau2} a). We focus on the  bottom left region bounded by left-most curve of $\mathcal{T}^1_0$ and the lowest curve of $\mathcal{T}^2_0$,  which is shown in  Fig. \ref{fig:tau1tau2} a). Notice that the  left-most and  the lowest curve  among all the stability switching curves are both part of $\mathcal{T}_0$.   By Theorem \ref{direction}, we can verify that the positive equilibrium  $E^*$ is stable in the bottom left region, since the crossing directions of the two switching curves (the black line and blue line) are all pointing outside of the region.  Moreover,  we can see that  the two stability switching curves  intersect at the point $(3.9042,1.406)$, and we denote the double Hopf bifurcation point by HH. For HH, using the normal form derivation process given in section \ref{normal form}, we have $\omega_1= 0.61081$, $\omega_2=0.94964    $,  $K_{11} =      0.0947 - 0.0071i, K_{21} =  -0.2689 + 0.4408i,K_{13} =
           0.1196 + 1.2137i,K_{23} =    1.6381 - 2.5531i,K_{2100} =   0.0154 - 0.0146i,
   K_{1011 }=        0.4878 + 0.2082i,K_{0021} =   -0.9861 - 0.9526i,K_{1110}=    -0.1778 - 0.1523i$.      Furthermore, we have the  normal form (\ref{normalformcylin}) with $\epsilon_1 =1, \epsilon_2 =   -1, b =    0.4946,c =     -11.5623, d =       -1$, and  $d-bc =       4.7192$.

This means that the unfolding system near the double Hopf bifurcation point HH is of type VIa.
 According to Guckenheimer and Holmes,
 \cite{Guckenheimer}  near the double Hopf bifurcation point there are
 eight different kinds of phase diagrams in eight different
 regions which are divided by semi-lines L1-L8 with
    \[\begin{aligned}
     &L_1:\tau_2=(\tau_1-3.9042)/(  -13.6972 )+1.406 ~~ ~ (\tau_1>3.9042);\\
    &L_2:\tau_2=(\tau_1-3.9042)/(    2.8383)+1.406 ~~ ~ (\tau_1>3.9042) ;\\
        &L_3:\tau_2=(\tau_1-3.9042)/(  1.2106
        )+1.406 ~~ ~(\tau_2>1.406) ;\\
    &L_4:\tau_2=(\tau_1-3.9042)/(   0.6790)+1.406+o(\tau_1-3.9042)~~ ~(\tau_2>1.406);\\
    &L_5:\tau_2=(\tau_1-3.9042)/(    0.6790)+1.406~~ ~(\tau_2>1.406);\\
    &L_6:\tau_2=(\tau_1-3.9042)/(  -3.5180)+1.406~~ ~(\tau_2>1.406);\\
       &L_7:\tau_2=(\tau_1-3.9042)/(  -13.6972)+1.406  ~~ ~(\tau_1<3.9042) ;\\
     &L_8:\tau_2=(\tau_1-3.9042)/(  2.8381)+1.406 ~~ ~(\tau_1<3.9042)  .\\
         \end{aligned}\]
 According to Fig. \ref{fig:VIa}, we have the bifurcation set near HH showing in Fig. \ref{fig:tau1tau2} b).  It is found that HH is the intersection of a supercritical Hopf bifurcation curve and a subcritical Hopf bifurcation curve.

\begin{figure}
               a) \centering
          \includegraphics[width=0.27\textwidth,height=0.3\textwidth]{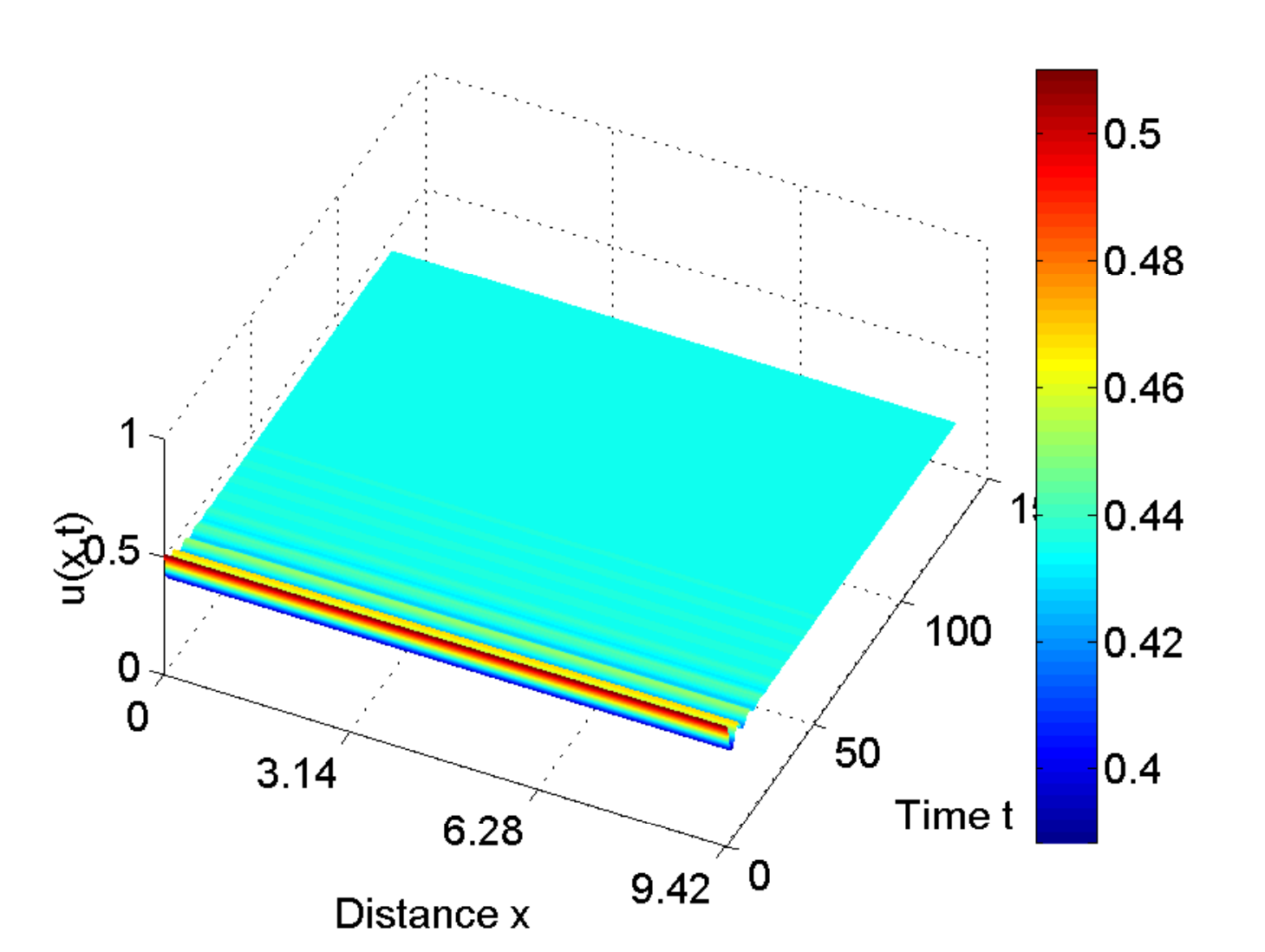}
          b)\includegraphics[width=0.27\textwidth,height=0.3\textwidth]{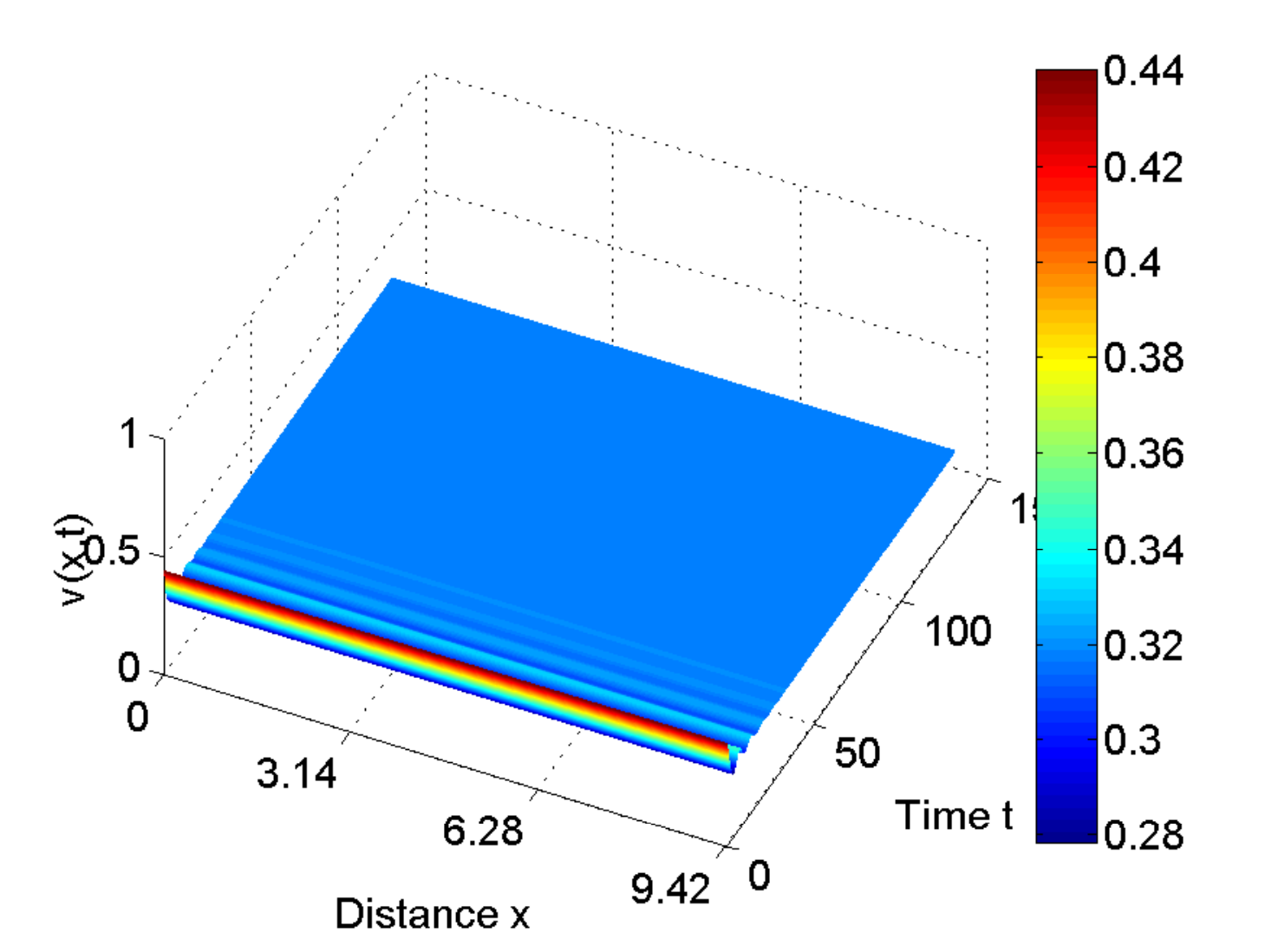}
 c) \includegraphics[width=0.27\textwidth,height=0.3\textwidth]{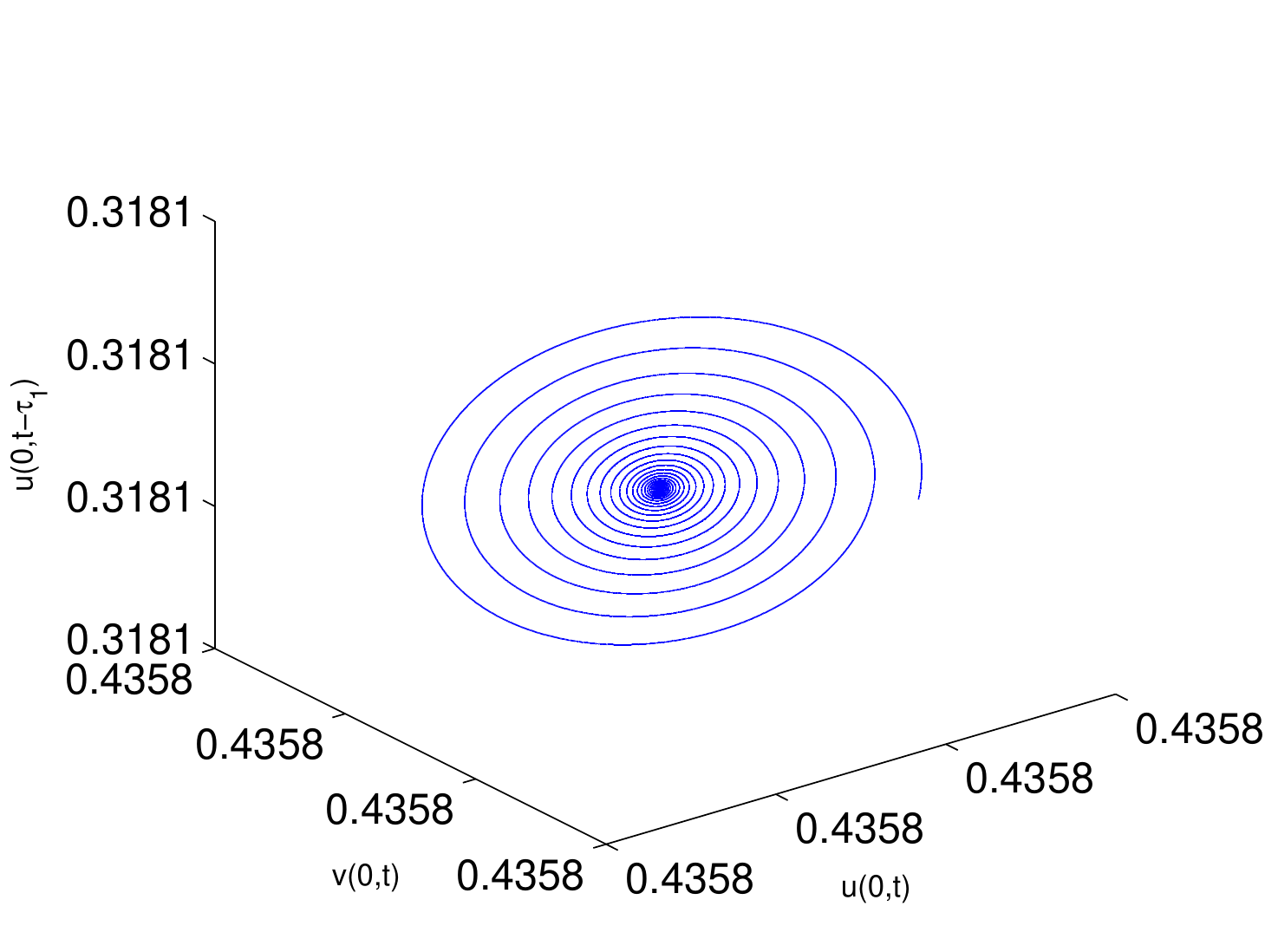}
      \caption  { When $\tau_1=1.74,\tau_2=0.67$ in $D_2$, the positive equilibrium is asymptotically stable.}
       \label{fig:E}
      \end{figure}

\begin{figure}
               a) \centering
          \includegraphics[width=0.27\textwidth,height=0.3\textwidth]{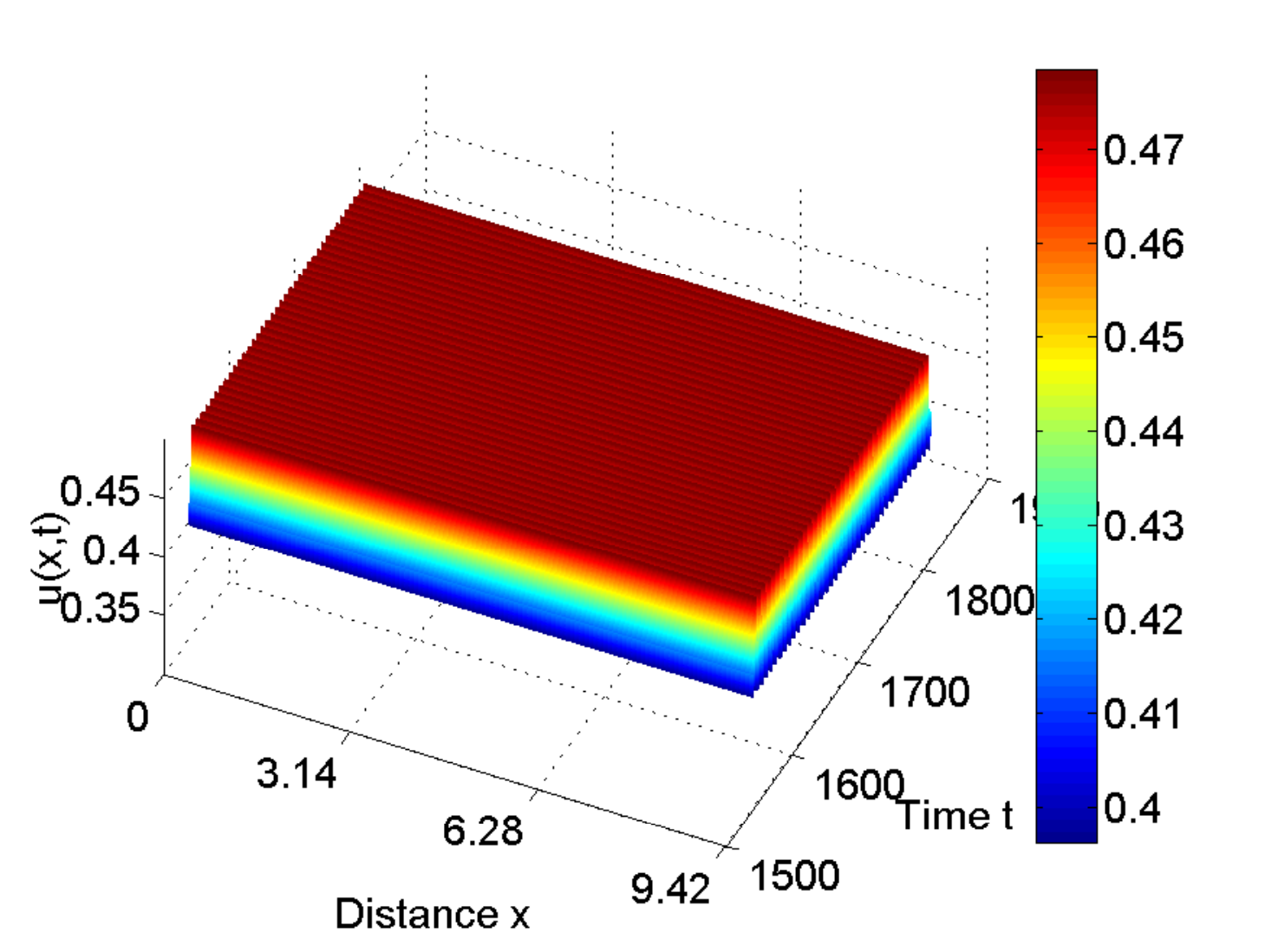}
          b)\includegraphics[width=0.27\textwidth,height=0.3\textwidth]{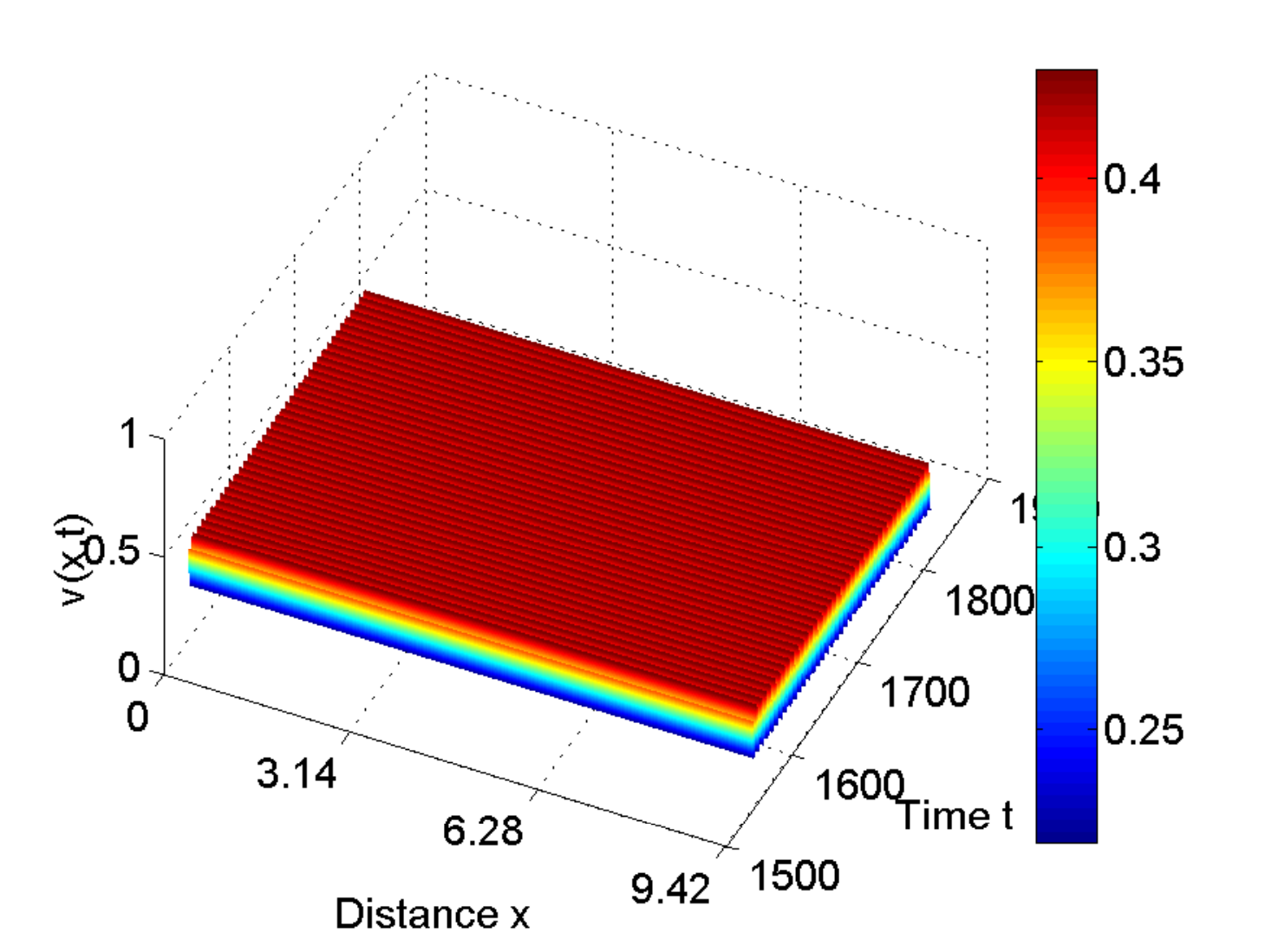}
  c) \includegraphics[width=0.27\textwidth,height=0.3\textwidth]{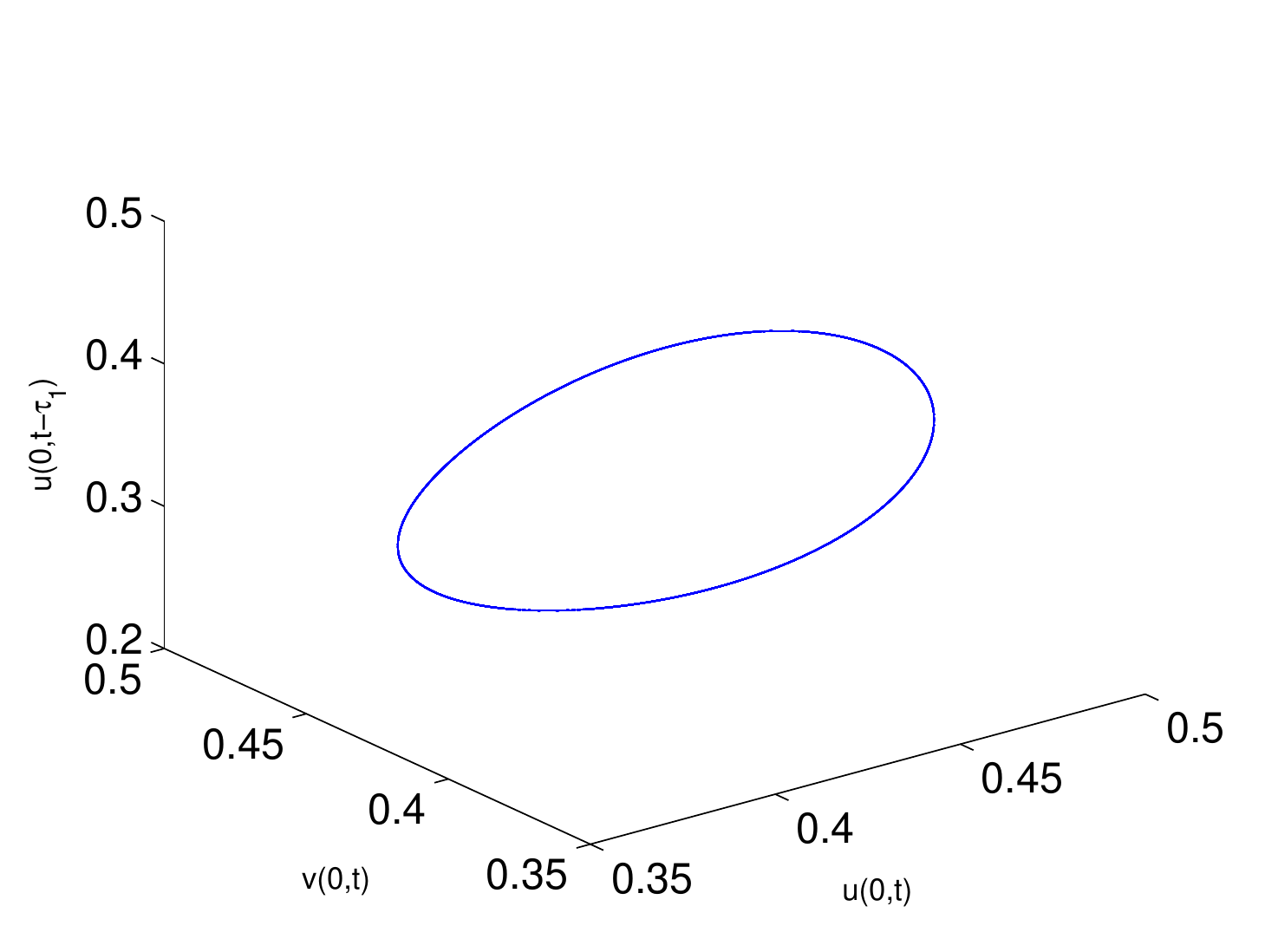}
      \caption  { When $\tau_1=3.62,\tau_2=1.435$ in $D_3$, there is a stable periodic solution.}
       \label{fig:periodic}
      \end{figure}

When $\tau_1=1.74,\tau_2=0.67$ in  $D_2$,  the  positive equilibrium is  a sink, which is shown in Fig. \ref{fig:E}. When $\tau_1=3.62,\tau_2=1.435$ in  region $D_3$, there is a stable periodic solutions originating from a supercritical Hopf bifurcation, which is shown in Fig. \ref{fig:periodic}.

\begin{figure}
             a)   \centering
        \includegraphics[width=0.275\textwidth,height=0.3\textwidth]{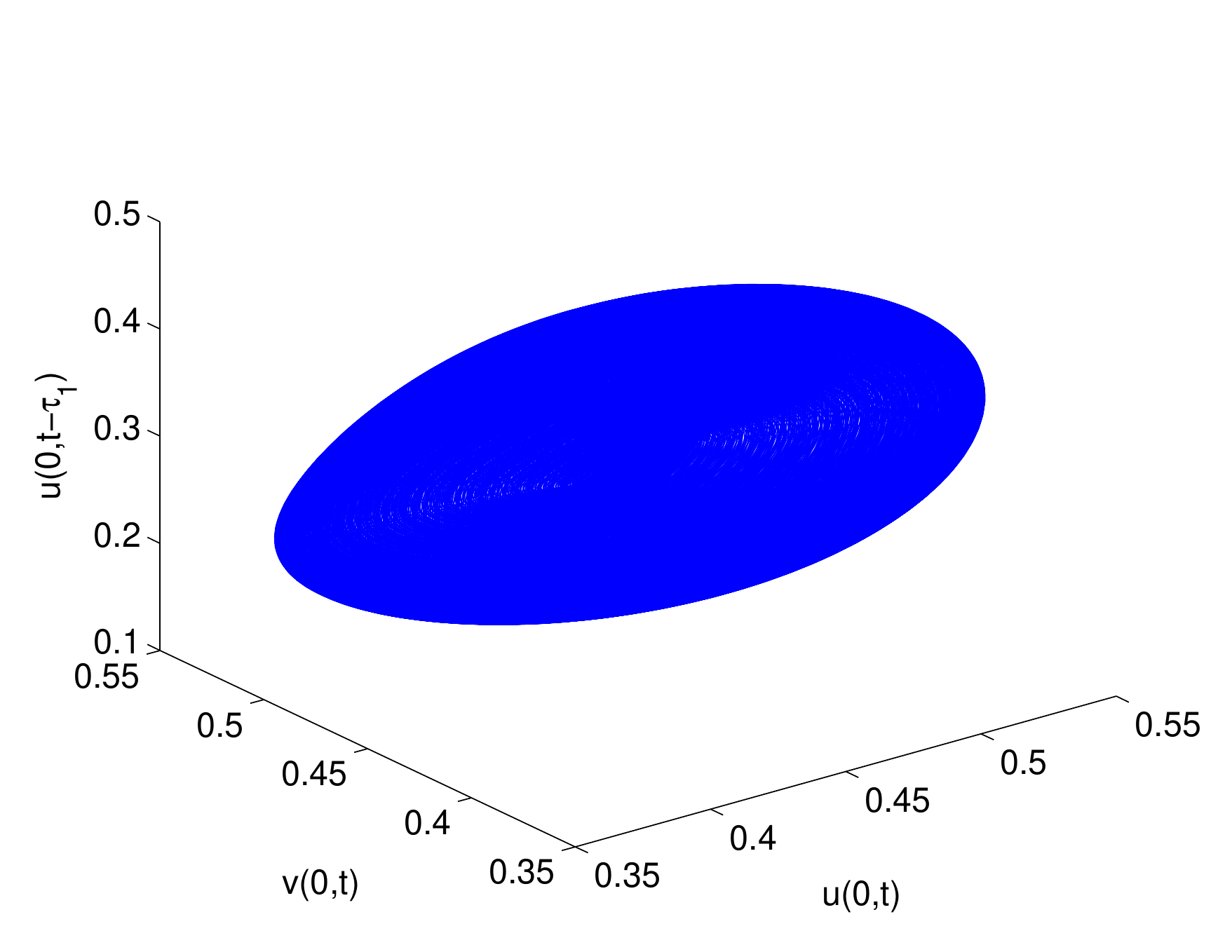}
       \includegraphics[width=0.275\textwidth,height=0.3\textwidth]{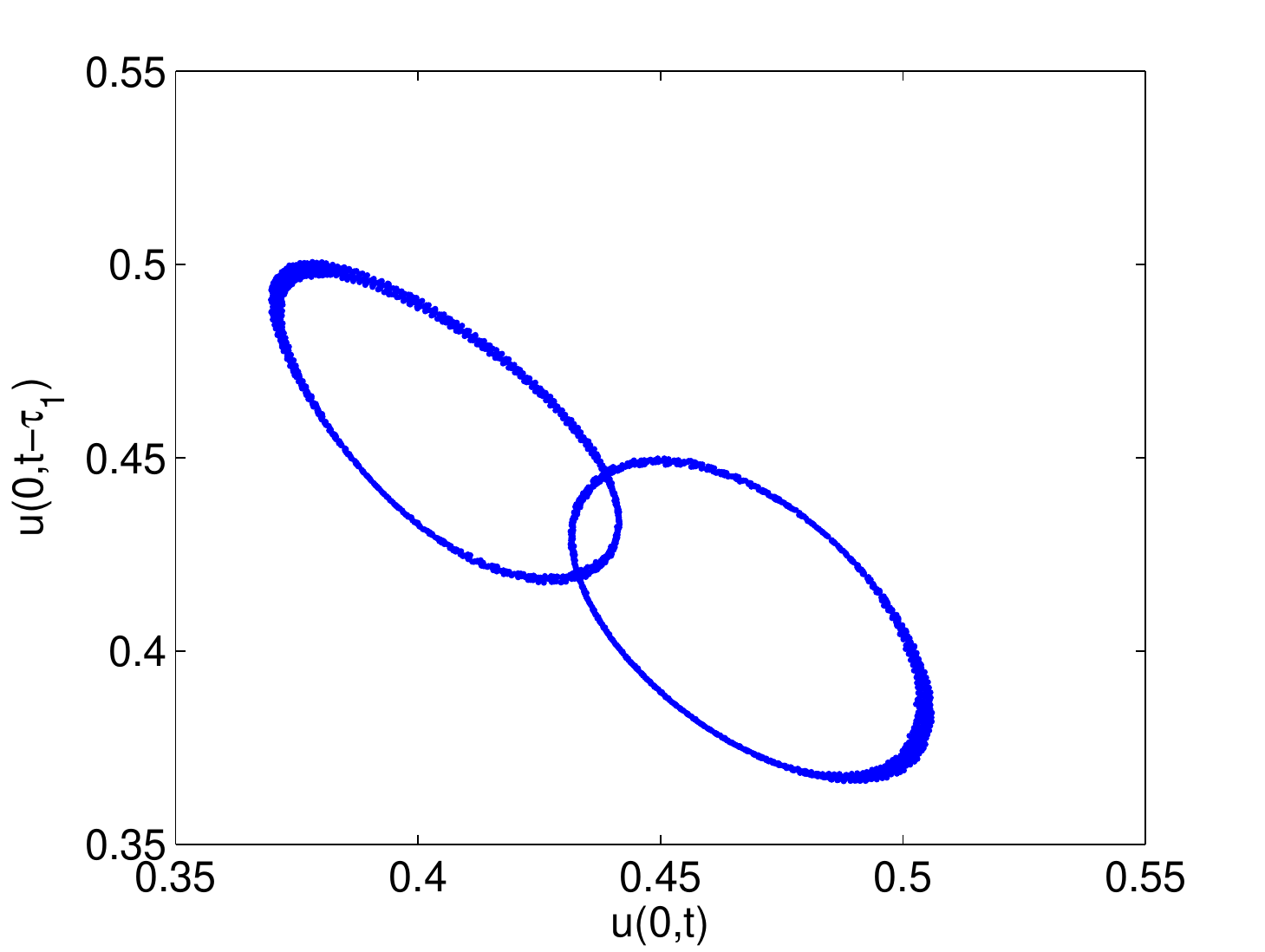}
\includegraphics[width=0.275\textwidth,height=0.3\textwidth]{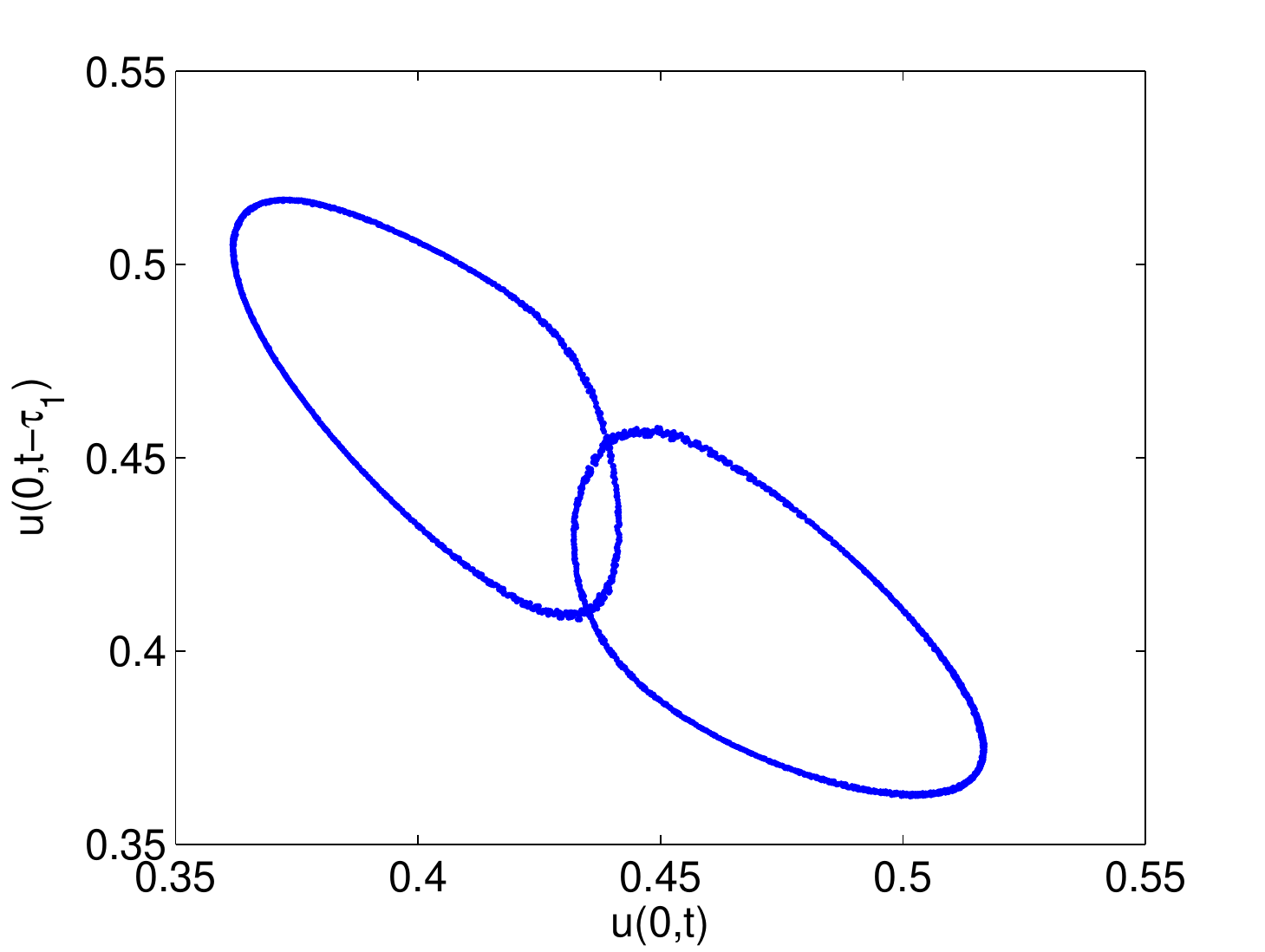}
               \\

           \centering
 b)                  \includegraphics[width=0.275\textwidth,height=0.3\textwidth]{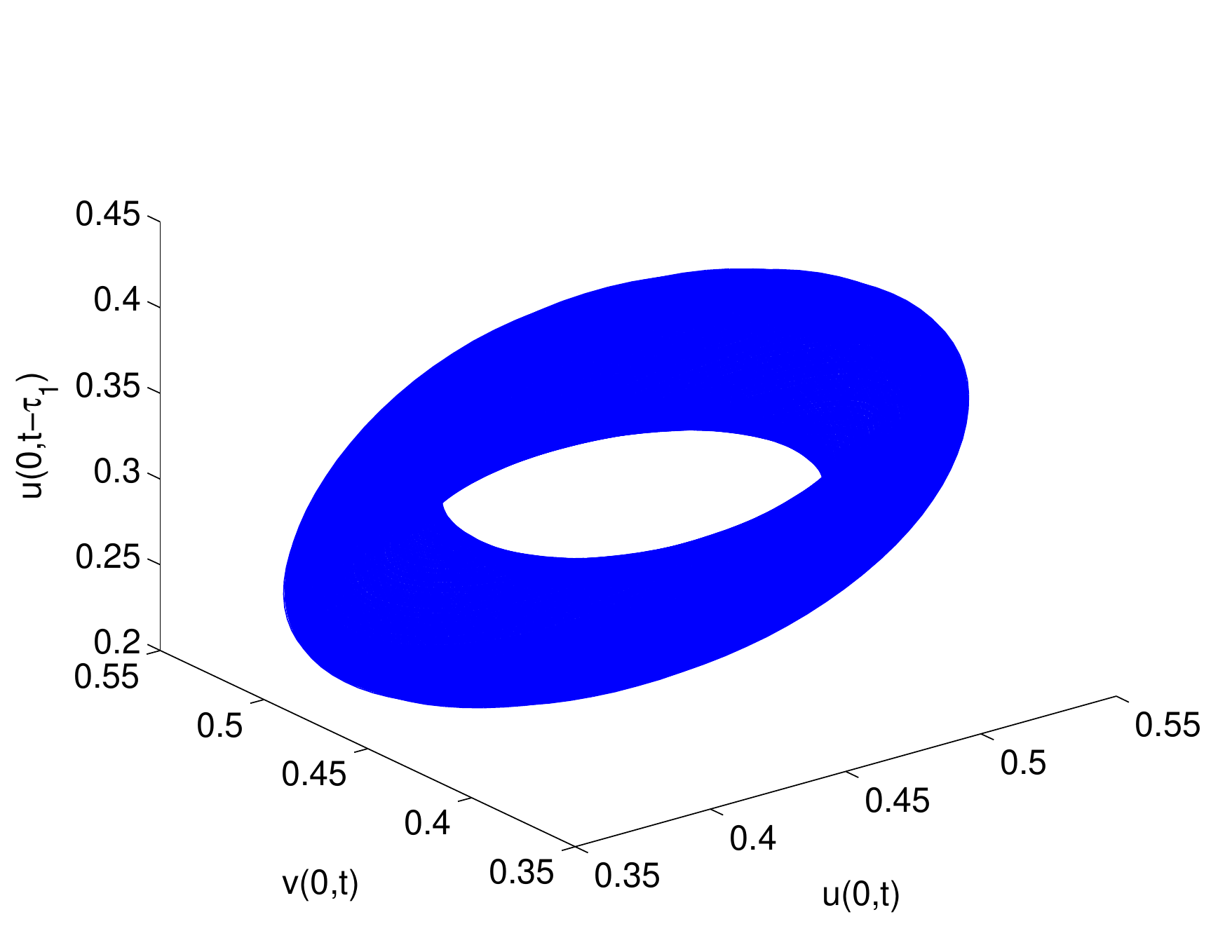}
                    \includegraphics[width=0.275\textwidth,height=0.3\textwidth]{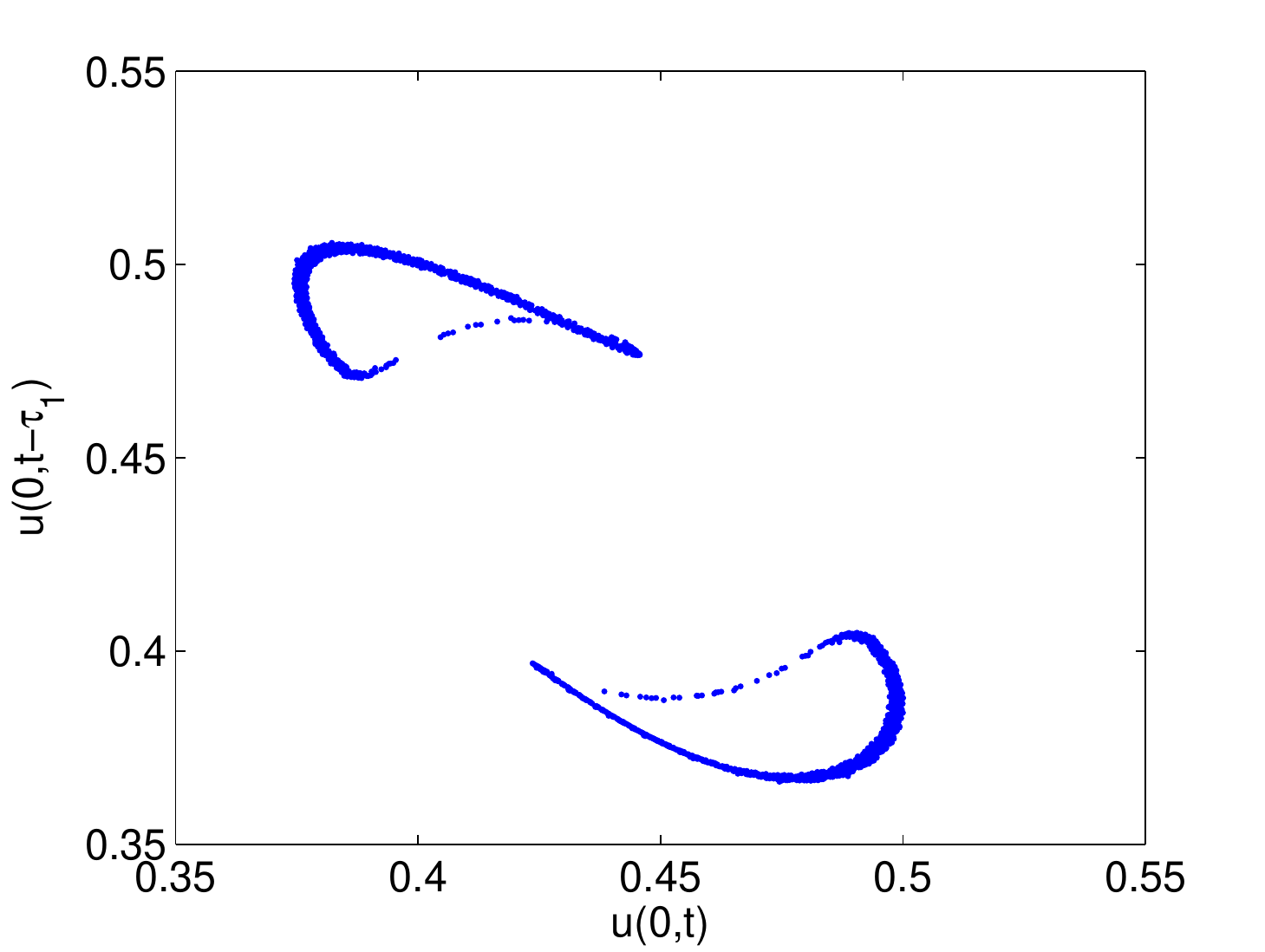}
             \includegraphics[width=0.275\textwidth,height=0.3\textwidth]{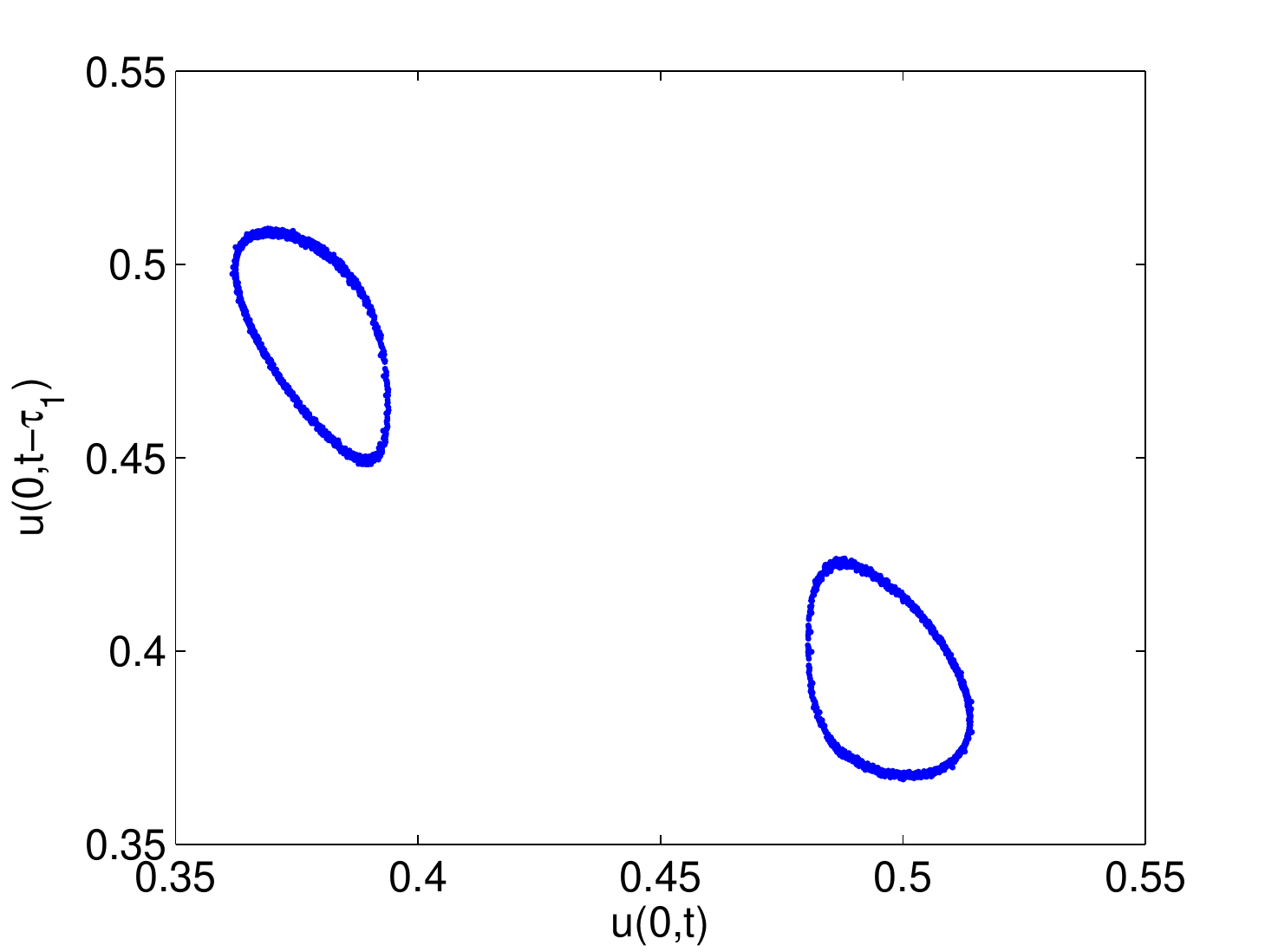}
                       \\
                        \centering
                     c)
          \includegraphics[width=0.275\textwidth,height=0.3\textwidth]{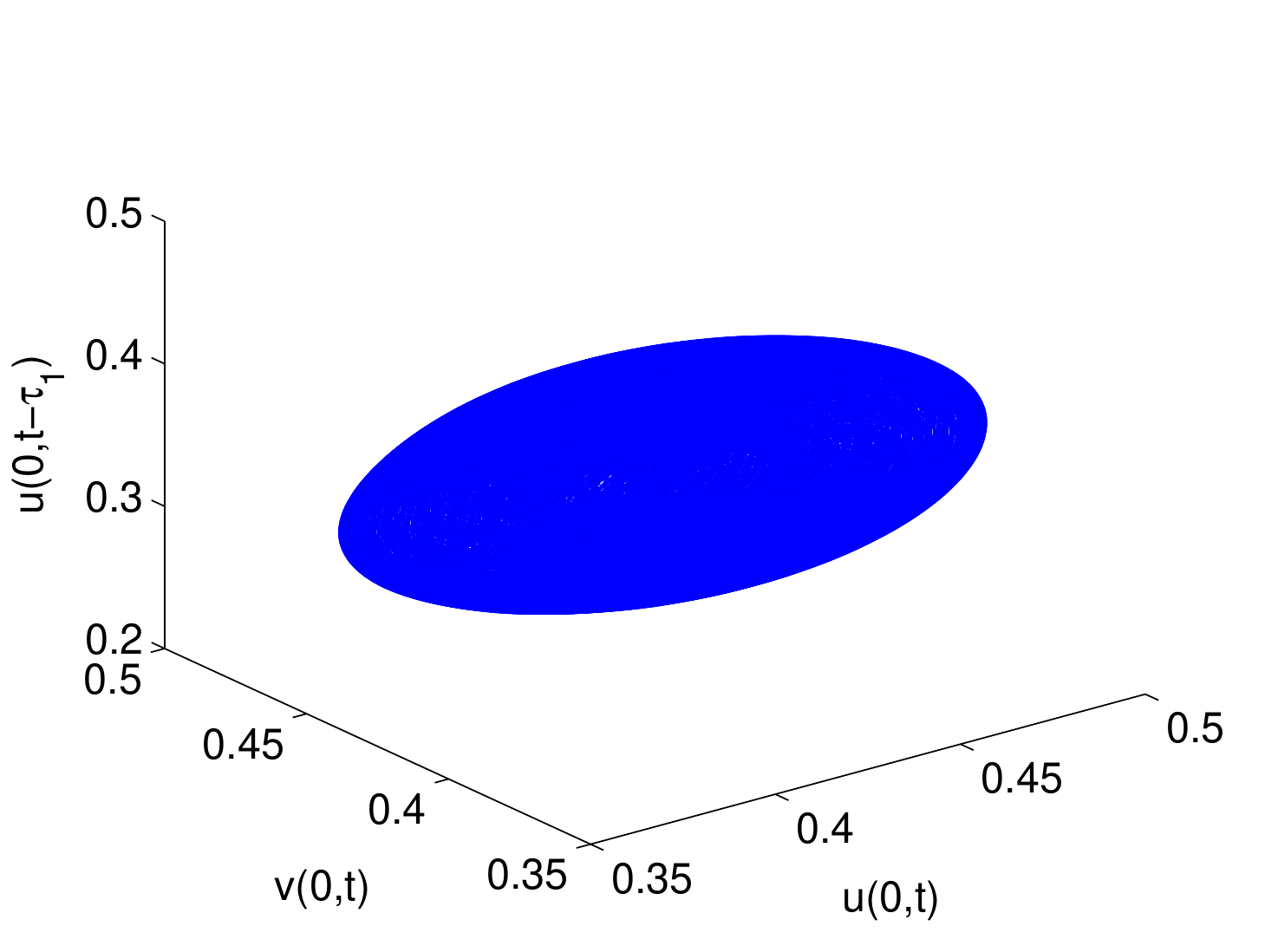}
         \includegraphics[width=0.275\textwidth,height=0.3\textwidth]{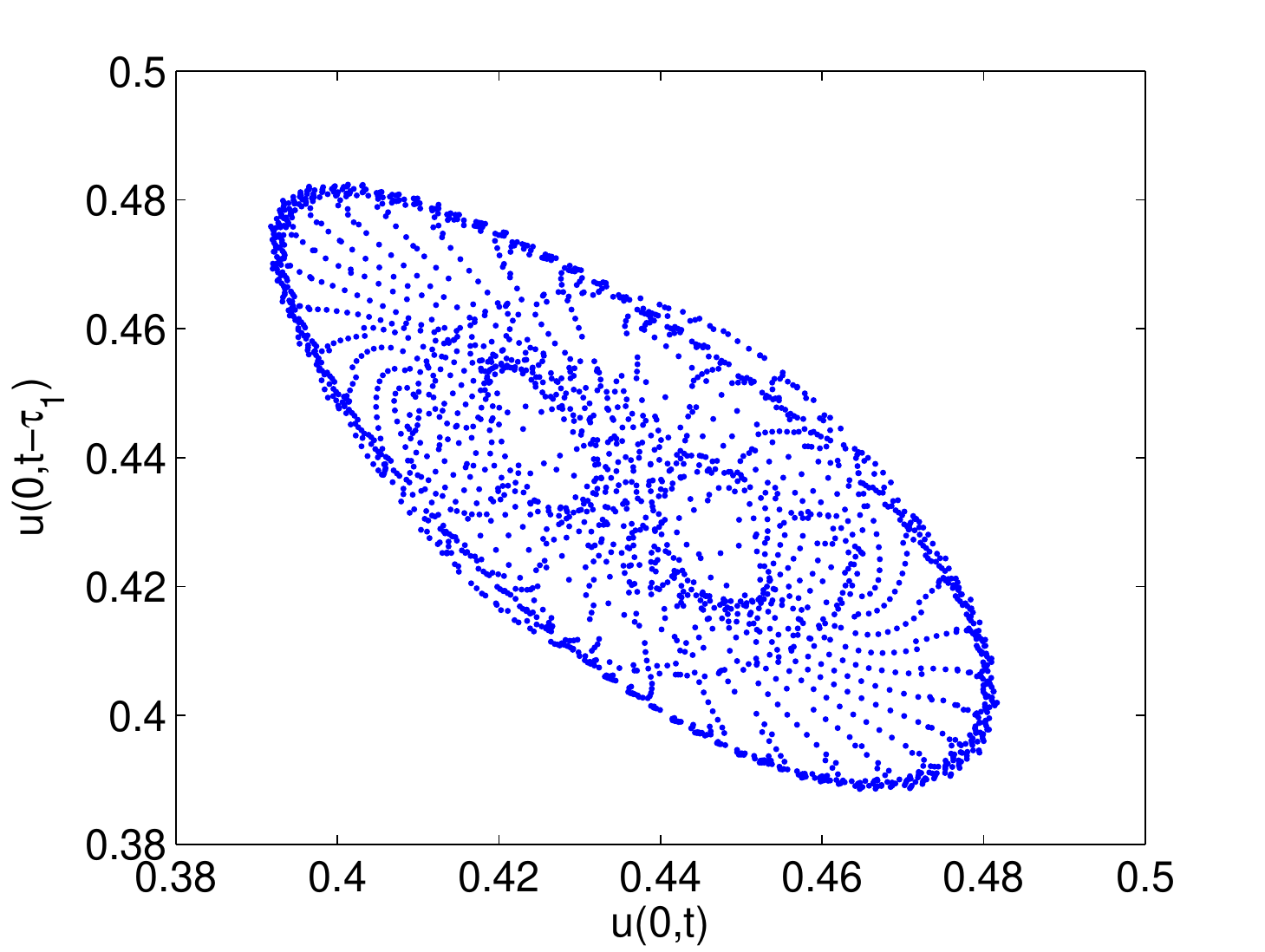}
  \includegraphics[width=0.275\textwidth,height=0.3\textwidth]{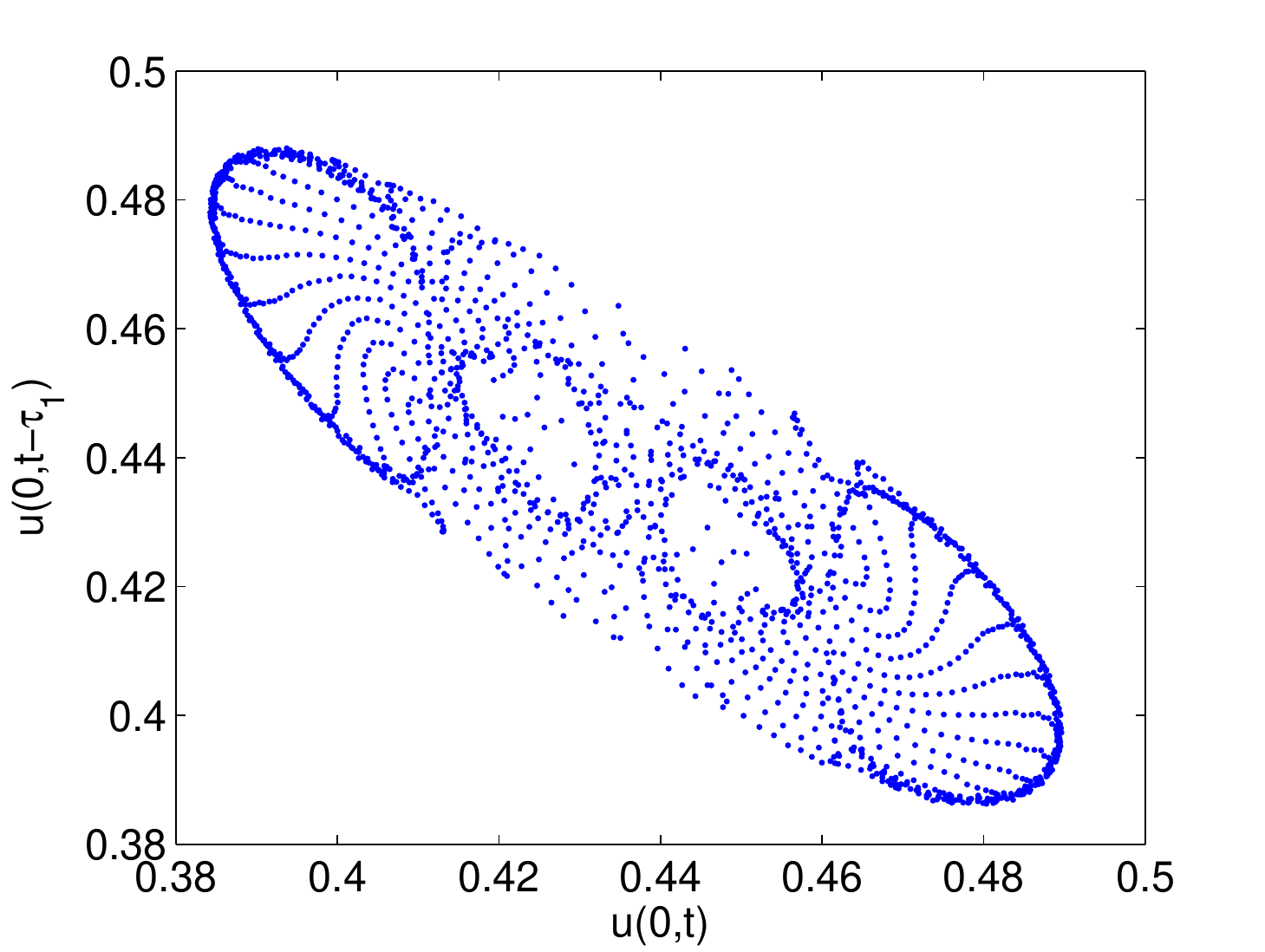}
      \caption  {  The phase portraits in $u(0,t)-v(0,t)-u(0,t-\tau_1)$, and the corresponding  Poincar\'e map on a  Poincar\'e section  $v(0,t)=v^*$ and $\dot{u}(x,t)=0$. The parameter is given as a) $\tau_1=3.82,\tau_2=1.4345$ in $D_4$,  b) $\tau_1=3.9043, \tau_2=1.418$ near $L_4$, c)  $\tau_1=3.905, \tau_2=1.4136$ in region $D_6$, respectively. Note that the transient states have been deleted for a clear expression.}
       \label{fig:23chaos}
      \end{figure}

  Finally, we  show the existence of quasi-periodic solutions and the results of the Poincar\'e map on a  Poincar\'e section. Due to the fact that the double Hopf bifurcation point is the intersection of two curves from $\mathcal{T}_0$, all periodic or quasi-periodic solutions nearby are  spatially homogeneous. Thus, we choose the solution curve  of $(u(0,t),v(0,t))$ at $x=0$ to show the rich dynamics. Since the periodic solutions oscillate in an infinite dimensional phase space, \cite{JWu} we give simulations near the double Hopf bifurcation point in the space $u(0,t)-v(0,t)-u(0,t-\tau_1)$ and choose the Poincar\'e sections $v(0,t)=v^*$, $\dot{u}(0,t)=0$, respectively,  in Fig. \ref{fig:23chaos}. The system exhibits rich dynamical behavior near the bifurcation point.  When $(\tau_1,\tau_2)=(3.82,1.4345)$ in $D_4$, system  (\ref{diffusion predator}) has a quasi-periodic solution on a 2-torus   (Fig. \ref{fig:23chaos} a)). It  becomes a quasi-periodic solution on a 3-torus, which breaks down through an orbit-connection bifurcation at  $(\tau_1,\tau_2)=(3.9043, 1.418)$ ( Fig. \ref{fig:23chaos} b)).     When the parameters vary and enter $D_6$, three-dimensional torus vanish.  Due to the fact that a vanishing 3-torus might accompany the phenomenon of chaos, \cite{P. Battelino,D. Ruelle,J.P. Eckmann}  near the double Hopf bifurcation point, we find   strange attractors exist.        We can see that when $\tau_1=3.905, \tau_2=1.4136$ in region $D_6$, system  (\ref{diffusion predator}) has a strange attractor,  which is shown  on the Poincar\'e section in Fig. \ref{fig:23chaos} c).

\section{Concluding remarks}
   This paper deals with a  modified Leslie-Gower predator-prey system with two delays and diffusion. We focus on the joint effect of two delays on the dynamical behavior of the system. Applying the method of stability switching curves, we find the stable region of the positive equilibrium and obtain Hopf bifurcation results.  By searching the intersection of  stability switching curves near the stable region, we get the double Hopf bifurcation point. Through the calculation of normal form of the system, we get the corresponding unfolding system and the bifurcation set. We theoretically prove and illustrate the existence of quasi-periodic solution on  two-torus, quasi-periodic solution on  three-torus, and even strange attractor.

   The theorem of Hopf bifurcation corresponding to systems with single parameter has been proposed for a long time. However, the theorem with two parameters has not been well stated. In this paper, we define   Hopf bifurcation curve on the plane $(\tau_1,\tau_2)$, and give the  sufficient condition of the existence of Hopf bifurcation in two-parameter plane.

The derivation process of normal form  for double Hopf bifurcation is very difficult, and the calculation is very long.  It is  even harder when we deal with systems with two  simultaneously varying delays, since the change of time scale $t\rightarrow\frac{t}{\tau_1}$ only transforms one delay to 1, and the other delay becomes $\frac{\tau_2}{\tau_1}$, which makes the Taylor expansion of the nonlinear term $G(\sigma, U^t)$   with $U^t(-\frac{\tau_2}{\tau_1})=U^t(-\frac{\tau_2^*+\sigma_2}{\tau_1^*+\sigma_1})$ very complicated. We should  notice that the method of  calculation of normal form used here can be also used in other systems with two delays, one delay or without delay by slight modifications.  The calculation formula of the  normal form  we give here is at a double Hopf bifurcation point  with $k_1=k_2=0$. The cases of $k_1=0,k_2\neq 0$ and  $k_1\neq 0,k_2\neq 0$ can be deduced
in a similar way, but in our model these double Hopf bifurcation points are not located at the boundary of stability region, hence unstable manifold always exists.
\section*{Supplementary Material}
In the supplementary material, we give the detailed  calculation process of  second and third order normal forms  near double Hopf bifurcation.

\section*{Acknowledgments}
The authors are grateful to the handling editor and anonymous  referees for their careful reading of the manuscript and  valuable comments, which improve the exposition of the paper very much.  Y. Du is supported by Education Department of Shaanxi Province (grant No. 18JK0123). B. Niu and J. Wei are supported by  National Natural Science Foundation of China (grant Nos.11701120 and No.11771109) and the Foundation for Innovation at HIT(WH).

\end{document}